\input amssym.def
\input amssym  
\vphantom{}

\def\R{{\bf R}} 
\def\C{{\bf C}}

\def\D{{\bf D}}

\def\N{{\bf N}}

\def\n{\overline{\nu}}
\catcode`\|=13 
\font\helvb=cmssbx10
\font\eightrm=cmr8
\font\eighti=cmmi8
\font\eightsy=cmsy8
\font\eightbf=cmbx8
\font\eighttt=cmtt8
\font\eightit=cmti8
\font\eightsl=cmsl8
\font\sixrm=cmr6
\font\sixi=cmmi6
\font\sixsy=cmsy6
\font\sixbf=cmbx6
\skewchar\eighti='177 \skewchar\sixi='177
\skewchar\eightsy='60 \skewchar\sixsy='60

\def\tenpoint{%
  \textfont0=\tenrm \scriptfont0=\sevenrm \scriptscriptfont0=\fiverm
  \def\rm{\fam0\tenrm}%
  \textfont1=\teni \scriptfont1=\seveni \scriptscriptfont1=\fivei
  \def\mit{\fam\@ne}\def\oldstyle{\fam1\teni}%
  \textfont2=\tensy \scriptfont2=\sevensy \scriptscriptfont2=\fivesy
    \def\itfam{4}\textfont\itfam=\tenit
  \def\it{\fam\itfam\tenit}%
  \def\slfam{5}\textfont\slfam=\tensl
  \def\sl{\fam\slfam\tensl}%
  \def\bffam{6}\textfont\bffam=\tenbf \scriptfont\bffam=\sevenbf
  \scriptscriptfont\bffam=\fivebf
  \def\bf{\fam\bffam\tenbf}%
  \def\ttfam{7}\textfont\ttfam=\tentt
  \def\tt{\fam\ttfam\tentt}%
  \abovedisplayskip=6pt plus 2pt minus 6pt
  \abovedisplayshortskip=0pt plus 3pt
  \belowdisplayskip=6pt plus 2pt minus 6pt
  \belowdisplayshortskip=7pt plus 3pt minus 4pt
  \smallskipamount=3pt plus 1pt minus 1pt
  \medskipamount=6pt plus 2pt minus 2pt
  \bigskipamount=12pt plus 4pt minus 4pt
  \normalbaselineskip=12pt
  \setbox\strutbox=\hbox{\vrule height8.5pt depth3.5pt width0pt}%
  \normalbaselines\rm}

\def\eightpoint{%
  \textfont0=\eightrm \scriptfont0=\sixrm \scriptscriptfont0=\fiverm
  \def\rm{\fam0\eightrm}%
  \textfont1=\eighti \scriptfont1=\sixi \scriptscriptfont1=\fivei
  \def\oldstyle{\fam1\eighti}%
  \textfont2=\eightsy \scriptfont2=\sixsy \scriptscriptfont2=\fivesy
  \textfont\slfam=\eightit
  \def\sl{\fam\itfam\eightit}%
  \textfont\slfam=\eightsl
  \def\sl{\fam\slfam\eightsl}%
  \textfont\bffam=\eightbf \scriptfont\bffam=\sixbf
  \scriptscriptfont\bffam=\fivebf
  \def\bf{\fam\bffam\eightbf}%
  \textfont\ttfam=\eighttt
  \def\tt{\fam\ttfam\eighttt}%
  \abovedisplayskip=9pt plus 2pt minus 6pt
  \abovedisplayshortskip=0pt plus 2pt
  \belowdisplayskip=9pt plus 2pt minus 6pt
  \belowdisplayshortskip=5pt plus 2pt minus 3pt
  \smallskipamount=2pt plus 1pt minus 1pt
  \medskipamount=4pt plus 2pt minus 1pt
  \bigskipamount=9pt plus 3pt minus 3pt
  \normalbaselineskip=9pt
  \setbox\strutbox=\hbox{\vrule height7pt depth2pt width0pt}%
  \normalbaselines\rm}

\font\petcap=cmcsc10

\tenpoint
\hsize=12.5cm
\vsize=19cm
\parskip 5pt plus 1pt
\parindent=1cm
\baselineskip=13pt
\hoffset=-0.1cm 
\def\footnoterule{\kern-6pt
  \hrule width 2truein \kern 5.6pt} 

\def\ie{{\sl i.e.\ }}
\def\cf{{\rm cf.\ }}

\def\omini{\raise 1ex\hbox{\ept o}}
\def\emini{\raise 1ex\hbox{\ept e}}
\def\ermini{\raise 1ex\hbox{\ept er}}
\def\remini{\raise 1ex\hbox{\ept re}}

\def\lead{\leaders\hbox to 10pt{\hss.\hss}\hfill}
\def\somt#1|#2|{\vskip 8pt plus1pt minus 1pt
                \line{#1\lead #2}}
\def\soms#1|#2|{\vskip 2pt
                \line{\qquad #1\lead #2}}  
\def\somss#1|#2|{\vskip 1pt   
                 \line{\qquad\qquad #1\lead #2}}

\def\aujour{\ifnum\day=1 1\ermini\else\number\day\fi\
\ifcase\month\or janvier\or f\'evrier\or mars\or avril\or mai\or juin\or
juillet\or aout\or septembre\or octobre\or novembre\or d\'ecembre\fi\
\number\year}
\def\today{\ifcase\month\or january \or february \or march \or april
\or may \or june\or july\or august \or september\or october\or november\or
december\fi\ \number\day , \number\year}

\newskip\afterskip
\catcode`\@=11
\def\p@int{.\par\vskip\afterskip\penalty100} 
\def\p@intir{\discretionary{.}{}{.\kern.35em---\kern.7em}}
\def\pointir{\afterassignment\pointir@\global\let\next=}
\def\pointir@{\ifx\next\par\p@int\else\p@intir\fi\egroup\next}
\catcode`\@=12
\def|{\relax\ifmmode\vert\else\findef\fi}
\def\findef{\errhelp{Cette barre verticale ne correspond ni a un \vert
mathematique
                        ni a une fin de definition, le contexte doit vous
indiquer ce qui manque.
                        Si vous vouliez inserer un long tiret, le codage
recommande est ---,
                        dans tous les cas, la barre fautive a ete supprimee.}%
                        \errmessage{Une barre verticale a ete trouvee en
mode texte}}

\def\TITR#1|{\null{\mss\baselineskip=17pt
                           \vskip 3.25ex plus 1ex minus .2ex
                           \leftskip=0pt plus \hsize
                           \rightskip=\leftskip
                           \parfillskip=0pt
                           \noindent #1
                           \par\vskip 2.3ex plus .2ex}}
 
\def\auteur#1|{\penalty 500
               \vbox{\centerline{\si
                 \iffrance par \else par \fi #1}
                \vskip 10pt}\penalty 500}

\def\resume#1|{\penalty 100
                           {\leftskip=\parindent
                            \rightskip=\leftskip
                            \eightpoint\bgroup\petcap \skip\afterskip=0pt
                             \iffrance R\'esum\'e \else Abstract \fi\pointir
                            #1 \par}
                           \penalty -100}

\def\titre#1|{\null\baselineskip14pt
                           {\helvb
                           \vskip 3.25ex plus 1ex minus .2ex
                           \leftskip=0pt plus \hsize
                           \rightskip=\leftskip
                           \parfillskip=0pt
                           \noindent #1
                           \par\vskip 2.3ex plus .2ex}}


\def\section#1|{
                                \bgroup\bf
                                 \par\penalty -500
                                 \vskip 3.25ex plus 1ex minus .2ex
                                 \skip\afterskip=1.5ex plus .2ex
                                  #1\pointir}

\def\ssection#1|{
                                 \bgroup\petcap
                                  \par\penalty -200
                                  \vskip 3.25ex plus 1ex minus .2ex
                          \skip\afterskip=1.5ex plus .2ex
                                   #1\pointir}
 

\def\th#1|{ 
                   \bgroup \sl
                        \def\findef{\egroup\par}
                        \bgroup\petcap
                         \par\vskip 2ex plus 1ex minus .2ex
                         \skip\afterskip=0pt
                           #1\pointir}

\def\defi#1|{ 
                   \bgroup \rm
                        \def\findef{\egroup\par}
                        \bgroup\petcap
                         \par\vskip 2ex plus 1ex minus .2ex
                         \skip\afterskip=0pt
                           #1\pointir}

\def\rque#1|{\bgroup \sl
                          \par\vskip 2ex plus 1ex minus .2ex\skip\afterskip=0pt
                          #1\pointir}

\def\dem{\bgroup \sl
                  \par\vskip 2ex plus 1ex minus .2ex\skip\afterskip=0pt
                  \iffrance D\'emonstration\else Preuve\fi\pointir}

\def\preuve{\bgroup \sl
                  \par\vskip 2ex plus 1ex minus .2ex\skip\afterskip=0pt
                  \iffrance Preuve\else Preuve\fi\pointir}


\def\_#1{_{\baselineskip=.7 \baselineskip
                                                                       
\vtop{\halign{\hfil$\scriptstyle{##}$\hfil\cr #1\crcr}}}}

\def\build#1#2\fin{\mathrel{\mathop{\kern0pt#1}\limits#2}}

\def\frac#1/#2{\leavevmode\kern.1em
   \raise.5ex\hbox{$\scriptstyle #1$}\kern-.1em
      /\kern-.15em\lower.25ex\hbox{$\scriptstyle #2$}}

{\obeylines
}
 
\def\decale#1{\par\noindent\hskip 4em
               \llap{#1\enspace}\ignorespaces}
\newif\ifchiffre
\def\chiffre{\chiffretrue}
\chiffre
\newdimen\laenge
\def\lettre#1|{\setbox3=\hbox{#1}\laenge=\wd3\advance\laenge by 3mm
\chiffrefalse}
\def\article#1|#2|#3|#4|#5|#6|#7|%
    {{\ifchiffre\leftskip=7mm\noindent
     \hangindent=2mm\hangafter=1
\llap{[#1]\hskip1.35em}\bgroup\petcap #2\pointir {\sl #3}, {\rm #4}
\nobreak{\bf #5}
({\oldstyle #6}), \nobreak #7.\par\else\noindent \advance\laenge by 4mm
\hangindent=\laenge\advance\laenge by -4mm\hangafter=1
\rlap{[#1]}\hskip\laenge\bgroup\petcap #2\pointir {\sl #3}, #4 {\bf #5}
({\oldstyle
#6}), #7.\par\fi}} 
\def\livre#1|#2|#3|#4|#5|%
    {{\ifchiffre\leftskip=7mm\noindent
    \hangindent=2mm\hangafter=1
\llap{[#1]\hskip1.35em}\bgroup\petcap #2\pointir{\sl #3}, #4, {\oldstyle
#5}.\par
\else\noindent
\advance\laenge by 4mm \hangindent=\laenge\advance\laenge by -4mm
\hangafter=1
\rlap{[#1]}\hskip\laenge\bgroup\petcap #2\pointir
{\sl  #3}, #4, {\oldstyle #5}.\par\fi}}
\def\divers#1|#2|#3|#4|%
    {{\ifchiffre\leftskip=7mm\noindent
    \hangindent=2mm\hangafter=1
     \llap{[#1]\hskip1.35em}\bgroup\petcap #2\pointir #3, {\oldstyle #4}.\par
\else\noindent
\advance\laenge by 4mm \hangindent=\laenge\advance\laenge by -4mm
\hangafter=1
\rlap{[#1]}\hskip\laenge\bgroup\petcap #2\pointir #3,{\oldstyle #4}.\par\fi}}
\def\div#1|#2|#3|#4|
{{\ifchiffre\leftskip=7mm\noindent
\hangindent=2mm\hangafter=1
\llap{[#1]\hskip1.35em}\bgroup\petcap #2\pointir {\sl  #3},{\oldstyle #4}.\par
\else\noindent
\advance\laenge by 4mm \hangindent=\laenge\advance\laenge by -4mm
\hangafter=1
\rlap{[#1]}\hskip\laenge\bgroup\petcap #2\pointir {\sl  #3},{\oldstyle
#4}.\par\fi}}


\font\si=cmssi10

\font\mss=cmss12 scaled \magstep1 

\hsize=12.5cm
\vsize=19cm
\parskip 5pt plus 1pt
\parindent=1cm
\baselineskip=13pt

  \hoffset=0.7cm
 \voffset=0.8cm

\font\cms=cmbsy10 at 5pt
\font\gros=cmex10 at 17.28pt

\newif\iffrance
\def\francais{\francetrue 
\frenchspacing
\emergencystretch3em
}

\def\TITR#1|{\null{\mss\baselineskip=17pt
                           \vskip 3.25ex plus 1ex minus .2ex
                           \leftskip=0pt plus \hsize
                           \rightskip=\leftskip
                           \parfillskip=0pt
                           \noindent #1
                           \par\vskip 2.3ex plus .2ex}}
 
\def\auteur#1|{\penalty 500
               \vbox{\centerline{\si
                 \iffrance par \else par \fi #1}
                \vskip 10pt}\penalty 500}

\def\resume#1|{\penalty 100
                           {\leftskip=\parindent
                            \rightskip=\leftskip
                            \eightpoint\bgroup\petcap \skip\afterskip=0pt
                             \iffrance R\'esum\'e \else Abstract \fi\pointir
                            #1 \par}
                           \penalty -100}

\def\titre#1|{\null\penalty-500\baselineskip14pt
                           {\helvb
                           \vskip 3.25ex plus 1ex minus .2ex
                           \leftskip=0pt plus \hsize
                           \rightskip=\leftskip
                           \parfillskip=0pt
                           \noindent #1
                           \par\nobreak\vskip 2.3ex plus .2ex}\penalty 5000}


\def\summ{ \mathop{\sum}\limits}
\def\cupp{\mathop{\bigcup}\limits}
\def\capp{\mathop{\bigcap}\limits}
\def\opp{\mathop{\bigoplus}\limits}
\def\ott{\mathop{\bigotimes}\limits}

\def\hfl{{\hbox to 12mm{\rightarrowfill}}}

\def\egal{{{\lower5pt\hbox to 0.5cm{$\hrulefill$}}\atop 
{\raise10pt\hbox to 0.5cm{$\hrulefill$}}}}

\def\upuparrow{\eqalign{&\big\uparrow\cr\noalign{\vskip-14pt}&\kern-1.8pt
               \uparrow\cr}}
\def\ddarrow{\eqalign{&\raise15pt\hbox{$\big\downarrow$}\cr
             \noalign{\vskip-25pt} &\big\downarrow\cr}}

\def\cu{\hbox{\cms\char'133}}
\def\ca{\hbox{\cms\char'134}}  
 \def\hookup{{\lower6.3pt\hbox{\cu}\kern-4.32pt\big\uparrow}}
 \def\hookdown{\raise5.3pt\hbox{\ca}\kern-4.32pt\lower3pt\hbox{\big\downarrow}}

\def\grtilde{\kern10pt\hbox{\gros\char'147}}

\let\le=\leq 
\let\ge=\geq 

\def\S{\mathhexbox278\kern.15em}
\def\\ {\smallsetminus}
\def\cad{c'est-\`a-dire~}
\def\wt{\widetilde}
\def\wh{\widehat}
\def\ept{\eightpoint}

\def\la{\longrightarrow}
\def\ld {, \ldots,}
\def\simto{\build \la^ {\vbox to 0pt{\hbox{$_{\,\displaystyle\sim}$}}}\fin}

\def\coin{\mathrel{\raise2pt\hbox{$\scriptstyle
|$}\kern-2pt\lower2pt\hbox{$-$}}}
\def\rect{\mathrel{\lower2pt\hbox{$-$}\kern-2pt\raise2pt\hbox{$\scriptstyle
|$}}} 


 \def\num{n$^\circ$~}

\def\im{\mathop{\rm Im}\nolimits}
\def\resp{\mathop{\rm resp.}\nolimits}
 
{\obeylines
}


\long\def\nomm#1|#2|#3|{\line{$\vtop{\hsize2.5cm #1~}\vtop{\hsize1cm #2~}
\vtop{\hsize=9cm\normalbaselines\parshape 1 0cm 9cm #3.}$}
\medskip}

\long\def\nom#1|#2|{\centerline{$\vtop{\hsize=3cm{\bf #1}~: }
\vtop{\hsize=9.5cm\normalbaselines\parshape 1 0cm 9.5cm #2.}$}
\medskip}

\def\decale#1{\par\noindent\hskip 4em
               \llap{#1\enspace}\ignorespaces}

\def\bsmf{~Bull. Soc. Math. France~}

\def\boxit#1#2{\hbox{\vrule
 \vbox{\hrule\kern#1
  \vtop{\hbox{\kern#1 #2\kern#1}%
   \kern#1\hrule}}%
 \vrule}}

\newbox\texte
\def\texteencadre#1|#2|#3|{\setbox1=\vbox{#3}
\setbox\texte=\vbox{\hrule height#1pt%
\hbox{\vrule width#1pt\kern#2pt\vbox{\kern#2pt \hbox{\box1}\kern#2pt}%
\kern#2pt\vrule width#1pt}\hrule height#1pt}}

\def\blanc{\hbox to 0pt{\vrule height13pt depth5pt width0pt}}
\def\blancs{\hbox to 0pt{\vrule height9pt depth5pt width0pt}}
\def\blan{\hbox to 0pt{\vrule height6pt depth5pt width0pt}}
\def\bla{\hbox to 0pt{\vrule height7pt depth5pt width0pt}}
\def\hrulefill{\leaders\hrule height0.2pt\hfill}

\def\bb#1&#2\cr{
\hbox to 1cm{\hfil\strut#1~\vrule}
\hbox to 1cm{\hfil\strut#2~\vrule}
}
\def\aa#1&#2&#3&#4&#5&#6&#7&#8&#9\cr{\line{$
\hbox to 3cm{\vrule width 1pt \strut~#1\hfil\vrule}
\hbox to 1cm{\hfil\strut#2~\vrule}
\hbox to 1cm{\hfil\strut#3~\vrule}
\hbox to 1cm{\hfil\strut#4~\vrule}
\hbox to 1cm{\hfil\strut#5~\vrule}
\hbox to 1cm{\hfil\strut#6~\vrule}
\hbox to 1cm{\hfil\strut#7~\vrule}
\hbox to 1cm{\hfil\strut#8~\vrule width 1pt}
\bb#9\cr
$\hfil}}

\def\ac{{\cal A}}                       
                      
\def\cc{{\cal C}}           \def\cb{{\bf C}}            
           \def\db{{\bf D}}

\def\ic{{\cal I}}           \def\ib{{\bf I}}            
\def\jc{{\cal J}}                       
\def\kc{{\cal K}}                       
                       
\def\mc{{\cal M}}                       
                       \def\nnf{{\bf N}}
\def\oc{{\cal O}}                      
\def\pcc{{\cal P}}          \def\pb{{\bf P}}            
           \def\qb{{\bf Q}}            
           \def\rb{{\bf R}}

\def\xc{{\cal X}}                       
                       
           \def\zb{{\bf Z}}            


\font\tenmib=cmmib10 
\font\sevenmib=cmmib7
\font\fivemib=cmmib5
\expandafter\chardef\csname pre boldmath.tex at\endcsname=\the\catcode`\@
\catcode`\@=11
\def\hexanumber@#1{\ifcase#1 0\or 1\or 2\or 3\or 4\or 5\or 6\or 7\or 8\or
 9\or A\or B\or C\or D\or E\or F\fi}

\skewchar\tenmib='177 \skewchar\sevenmib='177 \skewchar\fivemib='177

\newfam\mibfam   
\textfont\mibfam=\tenmib \scriptfont\mibfam=\sevenmib
\scriptscriptfont\mibfam=\fivemib
\def\mib@hex{\hexanumber@\mibfam}
\mathchardef\bfalpha="0\mib@hex 0B
\mathchardef\bfbeta="0\mib@hex 0C
\mathchardef\bfgamma="0\mib@hex 0D
\mathchardef\bfdelta="0\mib@hex 0E
\mathchardef\bfepsilon="0\mib@hex 0F
\mathchardef\bfzeta="0\mib@hex 10
\mathchardef\bfeta="0\mib@hex 11
\mathchardef\bftheta="0\mib@hex 12
\mathchardef\bfiota="0\mib@hex 13
\mathchardef\bfkappa="0\mib@hex 14
\mathchardef\bflambda="0\mib@hex 15
\mathchardef\bfmu="0\mib@hex 16
\mathchardef\bfnu="0\mib@hex 17
\mathchardef\bfxi="0\mib@hex 18
\mathchardef\bfpi="0\mib@hex 19
\mathchardef\bfrho="0\mib@hex 1A
\mathchardef\bfsigma="0\mib@hex 1B
\mathchardef\bftau="0\mib@hex 1C
\mathchardef\bfupsilon="0\mib@hex 1D
\mathchardef\bfphi="0\mib@hex 1E
\mathchardef\bfchi="0\mib@hex 1F
\mathchardef\bfpsi="0\mib@hex 20
\mathchardef\bfomega="0\mib@hex 21
\mathchardef\bfvarepsilon="0\mib@hex 22
\mathchardef\bfvartheta="0\mib@hex 23
\mathchardef\bfvarpi="0\mib@hex 24
\mathchardef\bfvarrho="0\mib@hex 25
\mathchardef\bfvarsigma="0\mib@hex 26
\mathchardef\bfvarphi="0\mib@hex 27
\mathchardef\bfimath="0\mib@hex 7B
\mathchardef\bfjmath="0\mib@hex 7C
\mathchardef\bfell="0\mib@hex 60
\mathchardef\bfwp="0\mib@hex 7D
\mathchardef\bfpartial="0\mib@hex 40
\mathchardef\bfflat="0\mib@hex 5B
\mathchardef\bfnatural="0\mib@hex 5C
\mathchardef\bfsharp="0\mib@hex 5D
\mathchardef\bftriangleleft="2\mib@hex 2F
\mathchardef\bftriangleright="2\mib@hex 2E
\mathchardef\bfstar="2\mib@hex 3F
\mathchardef\bfsmile="3\mib@hex 5E
\mathchardef\bffrown="3\mib@hex 5F
\mathchardef\bfleftharpoonup="3\mib@hex 28
\mathchardef\bfleftharpoondown="3\mib@hex 29
\mathchardef\bfrightharpoonup="3\mib@hex 2A
\mathchardef\bfrightharpoondown="3\mib@hex 2B
\mathchardef\bflhook="3\mib@hex 2C 
\mathchardef\bfrhook="3\mib@hex 2D 
\mathchardef\bfldotp="6\mib@hex 3A 
\catcode`\@=\csname pre boldmath.tex at\endcsname

\def\bnu{\bar\nu}
\def\oro{\overline\rb_0}
\def\ta{\wh A}
\def\oxp{\overline X'}
\def\oa{\overline A}
\def\ox{\overline X}

\def\oi{\overline I}
\def\oj{\overline J}
\def\ov{\overline V}

\def\gr{\mathop{\rm gr}\nolimits}
\def\ogr{\overline{\gr}}
\def\inn{\mathop{\rm in}\nolimits}
\def\mod{\mathop{\rm mod}\nolimits}
\def\tot{\mathop{\rm Tot}\nolimits}
\def\spe{\mathop{\rm Specan}\nolimits}
\def\pro{\mathop{\rm Projan}\nolimits}
\def\red{\mathop{\rm red}\nolimits}
\def\reg{\mathop{\rm reg}\nolimits}

\def\oinn{\overline{\inn}}
\def\bu{\bar u}
\def\uu{\underline u}
\def\ol{\overline}
\def\ul{\underline}
\def\inff{\mathop{\inf}\limits}
\def\supp{\mathop{\sup}\limits}
\def\opp{\mathop{\oplus}\limits}
\def\ott{\mathop{\otimes}\limits}
\def\wtx{\wt X}
\def\wtu{\wt U}
\def\tf{\tilde f}
\def\tit{\tilde\theta }
\def\tic{\tilde\ic}
\def\tig{\tilde g}
\def\tx{\tilde x}
\def\tih{\tilde h}
\def\oic{\overline{\ic}}
\def\ojc{\overline{\jc}}
\def\hfll{{\hbox to 25mm{\rightarrowfill}}}
\def\trdot{{\triangle\kern-6pt\cdot}\kern6pt}
\def\trdott{{\triangle\kern-5pt\cdot}\kern5pt}

\font\cms=cmss17  

\francais
\vphantom{}
\vskip1cm

\TITR CL\^OTURE INT\'EGRALE DES ID\'EAUX ET \'EQUISINGULARIT\'E|
\auteur Monique LEJEUNE-JALABERT et Bernard TEISSIER|\par
\centerline{\it Avec un appendice de Jean-Jacques RISLER\rm}\par\bigskip
\noindent

\centerline{\bf English summary}
This text has two parts. The first one is the essentially unmodified text of our 1973-74 seminar on integral dependence in complex analytic geometry at the Ecole Polytechnique with J-J. Risler's appendix on the \L ojasiewicz exponents in the real-analytic framework. The second part is a short survey of more recent results directly related to the content of the seminar.\par\noindent
The first part begins with the definition and elementary properties of the $\bnu$ order function associated to an ideal $I$ of a reduced analytic algebra $A$. Denoting by $\nu_I (x)$ the largest power of $I$ containing the element $x\in A$, one defines $\bnu_I (x)=\hbox{\rm lim}_{i\to \infty}{\nu_I(x^k)/k}$. The second paragraph is devoted to the equivalent definitions of the integral closure of an ideal in complex analytic geometry, one of them being $\overline I=\{x\in A/\bnu_I (x)\geq 1\}$. The third paragraph describes the normalized blowing-up of an ideal and the fourth explains how to compute $\bnu_I (x)$ with the help of the normalized blowing-up of the ideal $I$. It contains the basic finiteness results of the seminar, such as the rationality of $\bnu_I (x)$ (which had been proved by Nagata in algebraic geometry, a fact of which we were not aware at the time), the definitions of the fractional powers of coherent sheaves of ideals and the proof of their coherency. Given a coherent sheaf ${\cal I}$ of ${\cal O}_X$-ideals on a reduced analytic space $X$ one can define for each open set $U$ of $X$ and $f\in \Gamma (U,{\cal O}_X)$ the number $\bnu_{\cal I}^U(f)$ as the infimum of the $\bnu_{{\cal I}_y}(f_y)$ for $y\in U$.\par Then one defines for each positive real number $\nu$ the sheaf $\overline{{\cal I}^\nu}$ (resp. $\overline{{\cal I}^{\nu +}}$) associated to the presheaf $$ U\mapsto \{f\in \Gamma (U,{\cal O}_X)/\bnu_{\cal I}^U(f)\geq \nu\}$$(resp. $$ U\mapsto \{f\in \Gamma (U,{\cal O}_X)/\bnu_{\cal I}^U(f)>\nu\}).$$
Finally one has the graded ${\cal O}_X/{\cal I}$-algebra $$\overline{\hbox{\rm gr}}_{\cal I}{\cal O}_X=\bigoplus_{\nu\in \R_0}\overline{{\cal I}^\nu}/\overline{{\cal I}^{\nu +}}.$$ One important result is then that this algebra is locally finitely generated and that locally there is a universal denominator $q$ in the sense that all nonzero homogeneous components of the graded algebra have degree in ${{1}\over{q}}\N$. \par
In \S 5 it is shown that one can compute $\bnu$ using analytic arcs $h\colon (\C,0)\to (X,x)$, and \S 6 shows that \L ojasiewicz exponents are the inverses of $\bnu$, which implies that they are rational.\par
Risler's appendix shows how to use blowing-ups to compute \L ojasiewicz exponents and prove their rationality in the real analytic case.\par
The complements, added for this publication, point to some developments directly related to the subject of the seminar:\par The first one is the proof in the spirit of the seminar of the classical \L ojasiewicz inequality $\vert \hbox{\rm grad}(f(z))\vert \geq C_1\vert f(z)\vert^{\theta}$  with $\theta <1$. \par Then we point to later work which shows that in fact given an ideal $I$ and an element $f\in A$ the rational number $\bnu_I(f)$ can be seen as the slope of one of the sides of a natural Newton polygon associated to $I$ and $f$, which is in several ways a better indicator of the relations of the powers of $f$ with the powers of $I$ and has some useful incarnations. The third complement points to results of Izumi using $\bnu$ to characterize the Gabrielov rank condition for a morphism of analytic algebras, the fourth is a presentation of a generalization due to Ciuper\v ca, Enescu and Spiroff of the rationality of $\bnu$ to the case of several ideals, where it becomes the rationality of a certain polyhedral cone. \par The fifth comment presents the connection of $\bnu$ with the {\it type} of ideals, which was introduced by D'Angelo in complex analysis and used recently by Heier for the proof of an effective Nullstellensatz. In the middle 1980's, A. P\l oski, J. Chadzy\'nski and T. Krasi\'nski found methods of evaluation for the local and global \L ojasiewicz exponents in inequalities of the form $\vert P(z)\vert \geq C\vert z\vert^\theta$ where either $P=(P_1,\ldots ,P_k)$ is a collection of analytic functions on $\C^n$ having an isolated zero at the origin and the inequality should be true for $\vert z\vert$ small enough, or $P$ is a collection of polynomials with finitely many common zeroes and the inequality should be true for $\vert z\vert$ large enough. The results on the type are of the same nature, because it follows from the seminar that the type is in fact a \L ojasiewicz exponent.\par
The sixth comment points to results of Morales and others about the Hilbert function associated to the integrally closed powers $\overline{I^n}$ of a primary ideal in an excellent local ring and the associated graded algebra.\par Finally we point to two different but not unrelated uses of what is in fact the main object of study in the seminar: the reduced graded ring $\overline{\hbox{\rm gr}}_IA$ defined and studied in \S 4. In [T5] the second author uses the fact that for the local algebra ${\cal O}$ of a plane analytic branch the algebra $\overline{\hbox{\rm gr}}_m{\cal O}$ is the algebra of the semigroup associated to the singularity and is a complete intersection (a result due to the first author) to revisit the local moduli problem. The key is that the local analytic algebra ${\cal O}$ of every plane branch in the same equisingularity class has the same $\overline{\hbox{\rm gr}}_m{\cal O}$ because it has the same semigroup, so that the branch is a deformation of the monomial curve corresponding to that algebra. In [Kn], Allen Knutson uses the same specialization to the ''balanced normal cone" corresponding to $\overline{\hbox{\rm gr}}_IA$ in intersection theory.\par\noindent
Each paragraph has its own bibliography. Unfortunately at the time of the seminar we were unaware of the beautiful results of Samuel, Rees and Nagata (see [Sa], [N], [R1], [R2], [R3] in the bibliography of the complements), of which it appears \it a posteriori \rm that some parts of the seminar are translations into the complex analytic framework. The demand for this text over the years, however, and the fact that some mathematicians are led to rediscover some of its results, indicate that its publication is probably of some use.

\par\medskip

\tenrm{\centerline{\bf Pr\'eambule} Ceci est le texte du s\'eminaire tenu \`a l'Ecole Polytechnique en 1973-74, r\'edig\'e presque aussit\^ot et paru comme pr\'epublication de l'Institut Fourier, auquel nous avons ajout\'e \`a la fin quelques indications sur des r\'esultats plus r\'ecents qui sont en rapport direct avec le texte et une bibliographie compl\'ementaire. Chaque paragraphe a par ailleurs sa propre bibliographie.} 
\titre Introduction|
La motivation originelle des r\'esultats du chapitre~I venait d'une
d\'emonstration du th\'eor\`eme de ``continuit\'e du contact" de Hironaka,
point important de sa d\'emonstration de la r\'esolution des singularit\'es des
espaces analytiques complexes.

Il s'est av\'er\'e que ces r\'esultats \'etaient aussi fort utiles dans
l'\'etude des probl\`emes d'\'equisingularit\'e, et c'est cela qui a \'et\'e
expos\'e dans la seconde partie du s\'eminaire qui sera r\'edig\'ee
ult\'erieurement. Une des id\'ees essentielles de ces applications se trouve
d'ailleurs d\'ej\`a au chapitre~I, dans le lien entre $\bnu$ et l'exposant de
\L ojasiewicz, lien qui permet d'alg\'ebriser des conditions de nature
transcendante dans certains cas, et en particulier d'appliquer \`a l'\'etude
de ``conditions d'incidences" du type conditions de Whitney la puissance de
l'alg\`ebre. (Voir en particulier Ast\'erisque \num 7--8, 1973).

Signalons pour terminer qu'\`a l'\'epoque du s\'eminaire nous ignorions
l'existence des travaux de Samuel (Some asymptotic properties of powers of
ideals, Annals of Math., Series 2, t. 56) et de Nagata (Note on a paper of
Samuel, Mem. Call. of Sciences Univ. of Kyoto, t.~30, 1956--1957) o\`u la
rationalit\'e de $\bnu$ en g\'eom\'etrie alg\'ebrique est d\'emontr\'ee par
une m\'ethode qui est essentiellement celle donn\'ee ici au \S4. Ces articles
nous ont \'et\'e signal\'es par L.~Szpiro~; nous l'en remercions. Apr\`es
r\'eflexion, il nous a sembl\'e que notre travail d'\'etude syst\'ematique de
la filtration par le $\bnu$ et de ses propri\'et\'es de finitude en
g\'eom\'etrie analytique complexe, pr\'ecisait suffisamment les r\'esultats
de Samuel et Nagata, et en montrait des applications assez nouvelles, pour
pr\'esenter un certain int\'er\^et par lui-m\^eme.
\par\medskip\noindent

\titre 0. Fonction d'ordre, gradu\'e associ\'e. D\'efinition de $\bnu$|

Tous les anneaux consid\'er\'es ici sont unitaires et commutatifs et les
homomorphismes d'anneaux transforment l'\'el\'ement unit\'e en l'\'el\'ement
unit\'e.

$\zb$ d\'esigne l'anneau des entiers relatifs, $\nnf$ les entiers non
n\'egatifs, $\rb$ le corps des r\'eels, $\rb_0$ les r\'eels non n\'egatifs,
$\rb_+$ les r\'eels positifs et $\oro =\rb_0\cup \infty$.

\section 0.1. Quelques rappels|

\th 0.1.1. D\'efinition|Soit $A$ un anneau. Soit $\mu :A\to\oro$ une
application. On dit que c'est une fonction d'ordre, si elle v\'erifie les
propri\'et\'es suivantes pour tout $(x,y)\in A\times A$~:

\decale{i)} $\mu (x+y)\ge\inf\big(\mu (x),\mu (y)\big)$

\decale{ii)} $\mu (x\cdot y)\ge \mu (x)+\mu (y)$

\decale{iii)} $\mu (0)=\infty\quad \mu (1)=0$.|

(Pour donner un sens pr\'ecis \`a ces conditions, on utilise les conventions
habituelles sur le symbole $\infty$, \`a savoir~: $\infty$ est le seul
\'el\'ement de $\oro$ \`a \^etre sup\'erieur \`a tout $d\in\rb_0$, $\infty +
d=d+\infty=\infty$ pour tout $d\in\oro$).

\th 0.1.2. D\'efinition|Soit $A$ un anneau. On appelle filtration
(d\'ecrois\-sante) sur $A$, une suite d\'ecroissante $(A_d)_{d\in\rb_0}$ de
sous-groupes de $A$ v\'erifiant~:

\decale{i)} $A_d\cdot A_e\subset A_{d+e}$ pour $d\in\rb_0, e\in\rb_0$

\decale{ii)} $A_0=A$

\decale{iii)} Il existe $d\in\rb_+$ tel que $A_d\ne A$.

\noindent $A$ est alors dit filtr\'e.|

\rque 0.1.3. Remarques|Si $\mu $ est une fonction d'ordre, on a toujours $\mu
(x)=\mu (-x)$, et si $\mu (x)<\mu (y),\mu (x+y)=\mu (x)$.

Si $A$ est un anneau filtr\'e, $A_d$ est un id\'eal de $A$ pour tout
$d\in\rb_0$ et pour tout $d\in\rb_+,A_d\ne A$.

\rque 0.1.4. Remarque|Il y a \'equivalence entre les notions de fonction
d'ordre et d'anneau filtr\'e.

La fonction d'ordre \'etant donn\'ee, on d\'efinit
$$A_d=\{x\in A:\mu (x)\ge d\},\quad d\in\rb_0~.$$
R\'eciproquement, la filtration \'etant donn\'ee, on pose~:
$$\mu (x)=\{\sup d:x\in A_d\}~.$$

\ssection 0.1.5|Dans ces conditions, on pose
$$d(x,y)=e^{-\mu (x-y)}\quad\hbox{pour}\quad x,y\quad\hbox{dans}\quad A~.$$
On v\'erifie que~:

$\bullet\quad d(x,x)=0$

$\bullet\quad d(x,y)=d(y,x)$

$\bullet\quad d(x,y)\le\sup (d(x,z), d(z,y))$

\noindent $d$ est donc un \'ecart ultram\'etrique sur $A$, invariant par les
translations et $A_d=\{x,d(0,x)\le e^{-d}\}$. La topologie d\'efinie par $d$
est compatible avec la structure d'anneau de $A$. C'est celle pour laquelle
les $A_d$ constituent un syst\`eme fondamental de voisinages de 0 dans $A$.
On obtiendrait d'ailleurs la m\^eme topologie en choisissant comme syst\`eme
fondamental de voisinages de 0 les $A_d$ o\`u $d$ parcourt $\nnf$. $A$ admet
alors un s\'epar\'e compl\'et\'e $\ta$ qui a lui-m\^eme une structure
d'anneau topologique filtr\'e et on sait que le morphisme canonique
$i:A\to\ta$ est une injection si et seulement si $A$ est s\'epar\'e. Ceci a
lieu si et seulement si l'une des conditions suivantes est satisfaite~:

$\bullet\quad \mu (x)=\infty\Rightarrow x=0$

$\bullet\quad \capp_{d\in\rb_0} A_d=\{0\}$.

\th 0.1.6. D\'efinition|Soient $A$ un anneau et $G$ une $A$-alg\`ebre. On dit que
$G$ est une $A$-alg\`ebre gradu\'ee si l'on s'est donn\'e une d\'ecomposition~:
$$G=\bigoplus_{d\in\rb_0}G_d$$
o\`u

1) $G_d$ est un $A$-module pour $d\in\rb_0$

2) $G_0=A$

3) $G_{d_1}\cdot G_{d_2}\subset G_{d_1+d_2}$ pour $d_1,d_2$ dans $\rb_0$.

\noindent On appelle $G_d$ la composante homog\`ene de degr\'e $d$ de $G$.|

\ssection 0.1.7|Soient $A$ un anneau et $\mu $ une fonction d'ordre. On pose
$$\eqalign{
A_d(\resp A^+_d)&=\{x\in A\quad \mu (x)\ge d(\resp > d)\},\quad d\in\rb_0\cr
\gr A&=\oplus_{d\in\rb_0}A_d/A^+_d\cr}$$
est une $A/A^+_0$-alg\`ebre gradu\'ee.

\ssection 0.1.8|La construction de 0.1.7 est fonctorielle.

\rque 0.1.9. Un exemple fondamental|Soient $A$ un anneau, $I$ un id\'eal de
$A$ ne contenant pas 1, et consid\'erons la suite d\'ecroissante d'id\'eaux de
$A$, $(I^n)_{n\in\nnf}$. Par convention $I^\circ =A$. C'est la {\it filtration}
$I$-{\it adique\/} de $A$. Pour \^etre coh\'erent avec 0.1.2, nous poserons
naturellement pour $d$ r\'eel non n\'egatif quelconque
$$I^d=I^{n(d)}$$
o\`u $n(d)$ est le plus petit entier sup\'erieur ou \'egal \`a $d$.

Dans ce cas particulier, nous noterons la fonction d'ordre $\nu _I$. Ainsi~:
$$\nu _I(x)=\sup\{n:n\in\nnf,x\in I^n\}$$
et elle est en fait \`a valeur enti\`ere. Quant au gradu\'e associ\'e, nous le
noterons $\gr_I A$. Soit $x$ un \'el\'ement de $A$ tel que $\nu _I(x)\in\rb_0$.
On note $\inn_I(x)$ l'image canonique de $x$ dans la composante homog\`ene de
degr\'e $\nu _I(x)$ de $\gr_I A$. Si $A$ est un anneau n\oe th\'erien et si $I$
est contenu dans le radical de $A$, la topologie $I$-adique sur $A$ est
s\'epar\'ee [1]. De plus $\ta$ est un $A$-module fid\`element plat [1].

Par exemple si $A=\cb\{x_1\ld x_r;y_1\ld y_s\}$ est l'anneau de s\'eries
convergentes \`a $r+s$ variables et si $I=(y_1\ld y_s)$, $\ta$ s'identifie \`a
$\cb\{x_1\ld x_r\}$ $[[y_1\ld y_0]]$ anneau de s\'eries formelles \`a $s$
variables sur $\cb\{x_1\ld x_r\}$ et $\gr_I A$ s'identifie \`a $\cb\{x_1\ld
x_r\}$ $[y_1\ld y_s]$ anneau de polyn\^omes \`a $s$ variables sur $\cb\{x_1\ld
x_n\}$.

\ssection 0.1.10|De m\^eme si $J$ est un id\'eal de $A$, on pose~:
$$\nu _I(J)=\{\sup n:n\in\nnf,J\subset I^n\}~.$$
On a \'evidemment
$$\nu _I(J)=\inf_{x\in J}\nu _I(x)$$
et si $(x_1\ld x_n)$ est un syst\`eme de g\'en\'erateurs de $J$
$$\nu _I(J)=\inf_{i=1\cdots n}\nu _I(x_i)~.$$

\section 0.2. Une autre fonction d'ordre $\bnu$. Quelques propri\'et\'es
g\'en\'erales|

Revenons \`a l'exemple 0.1.9. $A$ \'etant un anneau et $I$ un id\'eal de $A$,
$\gr_I A$ a en g\'en\'eral des \'el\'ements nilpotents. Ceci vient du fait
qu'on peut tr\`es bien avoir $\nu _I(x^k)>k\nu _I(x)$ avec $k$ entier positif.
Par exemple si $A=\cb\{x,y\}/x^2-y^3$ et si $I=(x,y)$, on a
$$\eqalign{
\gr_I A&=\cb[X,Y]/X^2\cr
\nu _I(x^2)&=3>2\nu _I(x)=2~.\cr}\leqno{\rm et}$$

\th 0.2.1. Lemme|Soient $A$ un anneau, $I$ un id\'eal de $A$ ne contenant pas
1. Soit $J$ un id\'eal de $A$. La suite
$$u_k={\nu _I(J^k)\over k},\quad k\in\nnf$$
est convergente dans $\oro$.|

\dem Soit $\bu (\resp \uu)$ la limite sup\'erieure (resp. inf\'erieure) de la
suite $u_k$. Il s'agit de montrer que $\uu=\bu$. Si $\uu=\infty$ ou $\bu=0$,
c'est clair. Nous supposerons donc $\uu$ fini et $\bu>0$. Par d\'efinition,
$$\matrix{
&\hfill *  
\hfill&\forall \varepsilon 
\hfill&\kern-5pt>0,\quad  \forall i 
\hfill&\kern-5pt\in\nnf,\quad  \exists j
\hfill&\kern-5pt\ge i 
\hfill&:\quad u_j
\hfill&\kern-5pt\le\uu+\varepsilon
\hfill&\cr 
&\hfill **  
\hfill&\forall \varepsilon
\hfill&\kern-5pt>0,\quad  \forall i
\hfill&\kern-5pt\in\nnf,\quad \exists j
\hfill&\kern-5pt \ge i  
\hfill&:\quad u_j
\hfill&\kern-5pt\ge \bu-\varepsilon  
\hfill&(\hbox{si~} \bu <\infty)\hfill\cr
&\hfill ***  
\hfill&\forall N
\hfill&\kern-5pt\in\nnf,\quad  \forall i
\hfill&\kern-5pt\in\nnf,\quad  \exists j
\hfill&\kern-5pt\ge i  
\hfill&:\quad u_j
\hfill&\kern-5pt\ge N 
\hfill&(\hbox{si~} \bu =\infty )~.\hfill\cr}$$  

Fixons un $\varepsilon >0$ (et si $\bu=\infty$ un $N$) et choisissons un
indice $i$ assez grand pour que ${1\over i}<{\varepsilon \over
\bu-\varepsilon }(\resp {\varepsilon \over N})$. D'apr\`es ** (resp. ***) il
existe $j\ge i$ tel que $u_j>\bu-\varepsilon (\resp N)$ et d'apr\`es *, il
existe $k\ge ji$ tel que $u_k<\uu +\varepsilon $.

Divisons $k$ par $j$, $k=j\ell + q$ o\`u $q<j$. On a alors 
$$
\uu +\varepsilon >{\nu _I(J^{j\ell+q})\over j\ell+q}\ge{\ell\nu
_I(J^j)\over \ell j+q}={\nu _I(J^j)\over j}(1-{q\over k})
\ge (\bu-\varepsilon )(1-{1\over i})\ge \bu-2\varepsilon$$
$\big(\resp N(1-{1\over i})\ge N-\varepsilon \big)$.

Nous avons ainsi montr\'e que pour tout $\varepsilon >0$, $\uu+\varepsilon
\ge\bu-2\varepsilon $ (resp. pour tout $\varepsilon >0$ et tout $N$,
$\uu+\varepsilon \ge N-\varepsilon $) et donc $\uu=\bu$.

\rque 0.2.2. Remarque|La suite $(u_k)_{k\in\nnf}$ n'est pas croissante en
g\'en\'eral. N\'ean\-moins si on fixe $i$ la sous-suite $(u_{i^n})_{n\in\nnf}$
est croissante. Ceci montre que~:
$$\lim_{k\to\infty}(u_k)=\sup_{k\in\nnf}u_k~.$$

\th 0.2.3. D\'efinition|Soient $A$ un anneau et $I$ un id\'eal de $A$ ne
contenant pas 1. Soient $x$ un \'el\'ement de $A$ et $J$ un id\'eal de $A$. On
pose
$$
\eqalign{
\bnu_I(x)&=\lim_{k\to\infty}\nu _I(x^k)/k\cr
\bnu_I(J)&=\lim_{k\to\infty}\nu _I(J^k)/k~.\cr}$$|

\rque 0.2.4. Un exemple|Limite de rationnels, $\bnu_I(x)$ n'est pas en
g\'en\'eral un rationnel comme le montre l'exemple suivant~:

Soit $A$ l'alg\`ebre du mono\"{\i}de additif $\rb_0$, \cad l'anneau de
polyn\^ome \`a une variable $X$ \`a c\oe fficients dans $\zb$ dont les
exposants prennent leurs valeurs dans $\rb_0$. Soit $I$ l'id\'eal engendr\'e
par $X$
$$\bnu_I(X^\lambda )=\lambda \quad\hbox{alors que}\quad \nu _I(X^\lambda
)=[\lambda ]$$
o\`u $[]$ d\'esigne la partie enti\`ere.

\th 0.2.5. Proposition|Si $(x_1\ld x_m)$ est un syst\`eme de g\'en\'erateurs
de~$J$,
$$\bnu_I(J)=\inf_{1\le i\le m}\bnu_I(x_i)~.$$|

\dem Tout d'abord, il est clair que $\bnu_I(J)\le \inf\bnu_I(x_i)$. En effet
$x^k_i\in J^k$, $k\in\nnf$, donc $\nu _I(J^k)\le\nu _I(x^k_i)$.

R\'eciproquement $J^k$ \'etant engendr\'e par les $x^{a_1}_1\cdots x^{a_m}_m$
tels que $a_1+\cdots + a_m=k$,
$$\nu _I(J^k)=\inf_{\Sigma a_i=k}\nu _I(x^{a_1}_1\cdots x^{a_m}_m)~.$$

Fixons un $\varepsilon >0$ et soit $k'_0$ le plus petit entier tel que pour
tout $k>k'_0$ on ait, pour tout $1\le i\le m$,
$$\bnu_I(x_i)-\varepsilon \le {\nu _I(x^k_i)\over k}~.$$

Soit $k_0=mk'_0$ et posons $N=\supp_{1\le i\le m,a\le k'_0}[a\bnu_I(x_i)-\nu
_I(x^a_i)]$. On a~:
$$\nu _I(x^{a_1}_1\cdots x^{a_m}_m)\ge \summ_{1\le i\le m}\nu
_I(x^{a_i}_i)=\summ_{a_i>k'_0}\nu _I(x^{a_i}_i)+\summ_{a_i\le k'_0}\nu
_I(x^{a_i}_i)~.$$
(Si $a_1+\cdots +a_m=k>k_0$, la 1\`ere sommation n'est certainement pas vide).\par\noindent
Donc~:
$$\nu _I(x^{a_1}_1\cdots
x^{a_m}_m)\ge\summ_{a_i>k'_0}\big(a_i\bnu_I(x_i)-a_i\varepsilon
\big)+\summ_{a_i\le k'_0}\big(a_i\bnu_I(x_i)-N\big)~.$$
Or puisque $a_1+\cdots +a_m=k$, on a $\summ_{a_i>k'_0}a_i\le k$. Ainsi
$$\nu _I(x^{a_1}_1\cdots x^{a_m}_m)\ge\summ_{1\le i\le
m}a_i\bnu_I(x_i)-\varepsilon k-mN\ge k(\inf_{1\le i\le m}\bnu_I(x_i)-\varepsilon
)-mN~.$$
Finalement $\nu _I(J^k)\ge k\big[\inff_{1\le i\le m}\bnu_I(x_i)-\varepsilon
-{mN\over k}\big]$. Faisant tendre $k$ vers l'infini, $N$ ne d\'ependant que de
$\varepsilon $, il vient 
$$\bnu_I(J)\ge\inf_{1\le i\le m}\bnu_I(x_i)-\varepsilon $$
et ceci \'etant vrai pour tout $\varepsilon $, $\bnu_I(J)\ge\inff_{1\le i\le
m}\bnu_I(x_i)$.

\th 0.2.6. Corollaire|$\bnu_I$ est une fonction d'ordre.|

\dem Soient $x,y$ dans $A$ et soit $J$ l'id\'eal engendr\'e par $x$ et $y$.
Alors $\bnu_I(J)=\inf\big(\bnu_I(x),\bnu_I(y)\big)$. Mais $x+y\in J$. Donc
$\bnu_I(x+y)\ge\bnu_I(J)$. C'est i). D'autre part $\nu _I(x^k\cdot
y^k)\ge\nu _I(x^k)+\nu _I(y^k)$. Ceci donne ii). Finalement puisque $\nu
_I(0)=\infty$, a fortiori $\bnu_I(0)=\infty$. De plus, puisque $1\notin I$,
$\nu _I(1)=0$. Mais $1^k=1$. Donc $\bnu_I(1)=0$.

\rque 0.2.7. Remarque|Si $x$ est un \'el\'ement nilpotent de $A$,
$\bnu_I(x)=\infty$. On a en fait le r\'esultat plus pr\'ecis suivant~:

\th 0.2.8. Proposition|Soient $A$ un anneau, $N$ son nilradical. Soient $I$
un id\'eal de $A$ ne contenant pas 1, $J$ un id\'eal de $A$. Soient
$A_1=A/N$, $I_1$, $J_1$ les images respectives de $I$, $J$. Si $J$ est de
type fini, alors
$$\bnu_I(J)=\bnu_{I_1}(J_1)~.$$|

\dem Soit $(x_1\ld x_m)$ un syst\`eme de g\'en\'erateurs de $J$. D'apr\`es
2.5, $\bnu_I(J)=\inff_{1\le i\le m}\bnu_I(x_i)$ et
$\bnu_{I_1}(J_1)=\inff_{1\le i\le m}\bnu_{I_1}$ (cl $x_i\mod N$). Il
suffit donc de montrer que si $y$ est un \'el\'ement de $A$, $y_1$ son
image dans $A_1$, $\bnu_I(y)=\bnu_{I_1}(y_1)$. Il est clair que
$\bnu_I(y)\le\bnu_{I_1}(y_1)$. D'autre part si $\nu _{I_1}(y^k_1)=\nu
(k)$, alors $y^k\in I^{\nu (k)}+N$ ou encore $y^k=z_k+u_k$ o\`u $\nu
_I(z_k)\ge\nu (k)$ et $u_k\in N$.

Soit $i$ le plus petit entier tel que $u^{i+1}_k=0$.

Pout tout $n\in \N$, on a
$$y^{kn}=z^n_k+{n\choose 1}z^{n-1}_k u_k+\cdots {n-i\choose i}z^{n-i}_k
u^i_k, $$ $$ \nu _I(y^{kn})\ge (n-i)\nu (k)~.$$
Faisons tendre $n$ vers l'infini en laissant $k$ fixe.
$$\eqalign{
\lim_{n\to\infty}{\nu _I(y^{kn})\over
n}&=\bnu_I(y^k)\ge\lim_{n\to\infty}{n-i\over n}\cdot \nu (k)=\nu (k).\cr
\bnu_I(y^k)&=k\bnu_I(y)\quad\hbox{et}\quad\bnu_I(y)\ge{\nu (k)\over
k}.\cr}\leqno{\rm Or}$$
Faisant maintenant tendre $k$ vers l'infini, on obtient
$\bnu_I(y)\ge\bnu_{I_1}(y_1)$.

\th 0.2.9. Proposition|On a~:
$$\eqalign{
\bnu_I(J^k)&=k\bnu_I(J)\quad k\in\nnf\cr
\bnu_{I^k}(J)&={1\over k}\bnu_I(J)\quad k\in\nnf.\cr}$$|

\dem La premi\`ere assertion est une cons\'equence imm\'ediate de 0.2.3.
D'autre part $k\cdot\nu _{I^k}(J^n)\le\nu _I(J^n)$ et donc
$\bnu_{I^k}(J)\le{1\over k}\bnu_I(J)$. Soit $\nu (n)=\nu _I(J^n)$. Divisant
$\nu (n)$ par $k$, on obtient $\nu (n)=qk+r$ o\`u $0\le r<k$ et $\nu
_{I^k}(J^n)\ge q={\nu (n)-r\over k}\ge{\nu (n)-k+1\over k}$. Faisant tendre $n$
vers l'infini, on obtient
$$\bnu_{I^k}(J)\ge{1\over k}\lim_{n\to\infty}{\nu (n)\over n}={1\over
k}\cdot\bnu_I(J)~.$$

\th 0.2.10. Proposition|Soit $A$ un anneau n\oe th\'erien r\'eduit et soit
$\oa$ le normalis\'e de $A$. Soit $I$ un id\'eal de $A$ ne contenant pas 1,
$J$ un id\'eal. On pose $J'=J\oa$, $I'=I\oa$.

Si $\oa$ est un $A$-module de type fini (par exemple si $A$ est un anneau
local excellent), $\bnu_I(J)=\bnu_{I'}(J')$.|

\dem Comme ci-dessus, $\bnu_I(J)\le\bnu_{I'}(J')$ est imm\'ediat. D'apr\`es
Artin-Rees, $A$ \'etant n\oe th\'erien et $\oa$ fini sur $A$, il existe
$n_0\in \nnf$ tel que si $n\ge n_0$
$$I'^n\cap A=I^{n-n_0}\cdot(I'^{n_0}\cap A).$$
Soit $\nu (n)=\nu _{I'}(J'^n)$
$$J^n\subset J'^n\cap A\subset I'^{\nu (n)}\cap A\subset I^{\nu (n)-n_0}$$
$${\nu _I(J^n)\over n}\ge{\nu (n)-n_0\over
n}\quad\hbox{et}\quad\bnu_I(J)\ge\bnu_I'(J')~.$$

\rque 0.2.11. Notations|Soient $A$ un anneau, $I$ un id\'eal ne contenant pas
1. $\bnu_I$ la fonction d'ordre correspondante. On note $\ogr_I A$ l'anneau
gradu\'e obtenu par la construction de 0.1.7.

Soit $x$ un \'el\'ement de $A$ tel que $\bnu_I(x)\in\rb_0$. On note $\oinn_I
x$ (ou $\oinn x$ lorsqu'aucune confusion n'est possible) l'image canonique de
$x$ dans la composante homog\`ene de degr\'e $\bnu_I(x)$ de $\ogr_I A$. C'est
un \'el\'ement non nul.

Soit $J$ un id\'eal de $A$ non inclus dans $A_\infty=\{x\in
A~\bnu_I(x)=\infty\}$. On note $\oinn_I(J,A)$ l'id\'eal de $\ogr_I A$
engendr\'e par les $\oinn_I x$ pour $x$ parcourant $I$ tels que
$\bnu_I(x)\in\rb_0$.

\rque 0.2.12. Remarques|$\ogr_I A$ est un anneau r\'eduit. Il existe un
homomorphisme canonique (d'alg\`ebres gradu\'ees)
$$\tau :\gr_I A\la\ogr_I A$$
dont le noyau est le nilradical de $\gr_I A$. C'est un isomorphisme si et
seulement si $\gr_I A$ est un anneau r\'eduit.

Il suffit de remarquer que tout \'el\'ement homog\`ene non nul de $\ogr_I A$
est de la forme $\oinn_I x$ et que $(\oinn_I x)^k=\oinn_I\cdot x^k\ne 0$ car
$\bnu_I(x^k)=k\bnu_I(x)$. De m\^eme tout \'el\'ement homog\`ene non nul de
$\gr_I A$ est de la forme $\inn_I x$ et $\tau (\inn_I x)=\oinn_I x$ ou 0 selon
que $\bnu_I(x)=\nu _I(x)$ ou $\bnu_I(x)>\nu _I(x)$. Cette derni\`ere condition
signifie qu'il existe $k$ tel que $\nu _I(x^k)>k\bnu_I(x)$ ou encore $(\inn_I
x)^k=0$. Si maintenant $\gr_I A$ est r\'eduit, on a toujours $\bnu_I(x)=\nu
_I(x)$.

\ssection 0.2.13|Soient $A$ un anneau, $I$, $J$ des id\'eaux tels que $1\notin
I+J$. La suite
$$0\la\sqrt{\oinn_I(J,A)}\la\ogr_I A\build\la^{\alpha }\fin\ogr_{I/J}
A/J$$
\vskip-24pt
{\leftskip-2cm {\cms Z}\par}
\noindent 
est exacte. (On prendra garde \`a ne pas ajouter un 0 \`a droite).

\dem Soit $H=\oinn_I x$ un \'el\'ement homog\`ene non nul de degr\'e $\bnu$ de
$\ogr_I A$. $H\in\ker \alpha $ si et seulement si $\bnu_{I/J}(x\mod J)>\bnu$.
Ceci se produit encore si et seulement s'il existe $k\in\nnf$ tel que
$$\nu _{I/J}(x^k\mod J)>k\bnu$$
ou encore $x^k=y+z$ o\`u $\nu _I(y)>k\bnu$ et $z\in J$. Or,
$$\bnu_I(x^k)=k\bnu,\bnu_I(y)>k\bnu;\bnu_I(z)=k\bnu$$
et
$$\oinn_I x^k=(\oinn_I x)^k=\oinn_I z\in\oinn_I(J,A)~.$$

\titre R\'ef\'erences|

\div 1|Bourbaki n|Alg\`ebre commutative, chapitres 3 et 4|{\rm~Hermann}|
\par\medskip\noindent
\titre 1. La notion de cl\^oture int\'egrale d'un id\'eal|

Dans tout ce paragraphe, $A$ est un anneau commutatif et unitaire et $I$ un
id\'eal propre.

\th 1.1. D\'efinition|On dit qu'un \'el\'ement $f\in A$ est entier sur $I$ (ou
satisfait une relation de d\'ependance int\'egrale) s'il existe une relation~:
$$f^k+a_1f^{k-1}+\cdots +a_k=0\leqno (1.1.1)$$ o\`u $a_i\in A$ et $\nu _I(a_i)\ge i$, \ie, $a_i\in I^i$ pour  $i=1\cdots
k$.|

\rque 1.2. Notation|Soit $A[T]$ l'anneau de polyn\^ome \`a une ind\'etermin\'ee
$T$ sur $A$. On note $\pcc(I)$ le sous-anneau $\opp_{n\in\nnf}I^nT^n$ de
$A[T]$.

Le lemme suivant relie la notion de d\'ependance int\'egrale sur un id\'eal,
avec la notion classique de d\'ependance int\'egrale sur un anneau.

\th 1.3. Lemme|L'\'el\'ement $f$ est entier sur $I$ si et seulement si $fT$ est
entier sur l'anneau $\pcc(I)$ au sens usuel ([4], Vol.1, chap. V).|

\dem Soit $f^k+a_1f^{k-1}+\cdots +a_k=0$, $a_i\in I^i$, une relation de
d\'ependance int\'egrale de $f$ sur $I$. Alors~:
$$(fT)^k+(a_1T)(fT)^{k-1}+\cdots +(a_kT^k)=0$$
est une relation de d\'ependance int\'egrale de $fT$ sur $\pcc(I)$.

R\'eciproquement, soit
$$(fT)^k+b_1(fT)^{k-1}+\cdots +b_k=0,\quad b_i\in \pcc(I)$$
une relation de d\'ependance int\'egrale de $fT$ sur $\pcc(I)$.

On en d\'eduit, en annulant les termes de degr\'e $k$ dans $A[T]$, la relation
de d\'ependance int\'egrale de $f$ sur $I$ cherch\'ee.

\th 1.4. Corollaire-D\'efinition|L'ensemble des \'el\'ements $f$ de $A$
entiers sur $I$ est un id\'eal qu'on appelle la cl\^oture int\'egrale de $I$
dans $A$ et qu'on note $\oi$. On dit que $I$ est int\'egralement clos si
$I=\oi$.|

\dem Il suffit de voir que si $f$ et $g$ sont entiers sur $I$, $f+g$ l'est
aussi. Or la fermeture int\'egrale de $\pcc(I)$ dans $A[T]$ est un sous-anneau de
$A[T]$.

\th 1.5. Lemme|Si $I$ et $J$ sont des id\'eaux de $A$, tels que $I\subset J$,
alors $\oi\subset\oj$.|\par\noindent C'est \'evident d'apr\`es 1.1.

\th 1.6. Lemme|$I\subset\oi\subset\sqrt{I}$. En particulier si $I$ est un
id\'eal \'egal \`a sa racine, $I$ est int\'egralement clos. (On se gardera de
croire que les puissances de $I$ sont alors \'egalement des id\'eaux
int\'egralement clos.)|

\dem 1.1. montre que si $f$ est entier sur $I$, $f^k\in I$.

\th 1.7. Lemme|La fermeture int\'egrale de $\pcc(I)$ dans $A[T]$ est
$\opp_{n\in\nnf}\overline{I^n} T^n$.|

\dem Tout d'abord [1] p.~30, la fermeture int\'egrale de $\pcc(I)$ dans $A[T]$
est un sous-anneau gradu\'e de $A[T]$. Il suffit donc d'en d\'eterminer les
composantes homog\`enes. Comme en 1.3, si $fT^n$ est entier sur $\pcc(I)$, la
relation de d\'ependance int\'egrale~:
$$(fT^n)^k+\Sigma b_i(fT^n)^{k-i}=0,\quad b_i\in \pcc(I)$$
fournit, en annulant les termes de degr\'e $nk$ dans $A[T]$, la relation de
d\'ependance int\'egrale de $f$ sur $I^n$ cherch\'ee.

\th 1.8. Corollaire|

\decale{i)} $(\oi)^n\subset\overline{I^n}$

\decale{ii)} $I\cdot \overline{I^n}\subset\ol{I^{n+1}}$

\decale{iii)} $\build I^{\ =}\fin =\oi$.|

\dem

\decale{i)} La fermeture int\'egrale de $\pcc(I)$ dans $A[T]$ est une
sous-alg\`ebre de $A[T]$ contenant $\oi T$. Elle contient donc $\oplus
(\oi)^nT^n$.

\decale{ii)} La fermeture int\'egrale de $\pcc(I)$ dans $A[T]$ est une
sous-alg\`ebre de $A[T]$ contenant $I\cdot T$ et $\overline{I^n} T^n$. Elle contient donc
$I\cdot\overline{I^n} T^{n+1}$.

\decale{iii)} Soit maintenant $f$ entier sur $\oi$. Alors d'apr\`es 1.3. $fT$ est
entier sur $\pcc(\oi)$ qui est une sous-alg\`ebre de $\opp_{n\in\nnf}\overline{I^n} T^n$ et
est de ce fait enti\`ere sur $\pcc(I)$.

\th 1.9. Proposition|$f$ est entier sur $I$, si et seulement s'il existe un
$A$-module de type fini $M$ tel que~:

\decale{i)} $fM\subset I\cdot M$

\decale{ii)} Si $aM=0$, il existe $k\ge 0$ tel que $af^k=0$.|

\dem Supposons $f$ entier sur $I$ et consid\'erons une relation de d\'ependance
int\'egrale (1.1.1) satisfaite par $f$. Soit $I_0$ un id\'eal de type fini de
$A$ inclus dans $I$ tel que $\nu _{I_0}(a_i)\ge i$, $i=1\ld k$. Alors $f^k\in
I_0\cdot (I_0+fA)^{k-1}$ et donc $(I_0+fA)^k\subset I_0\cdot (I_0+fA)^{k-1}$.
Ainsi, nous pouvons choisir $M=(I_0+fA)^{k-1}$. Le $A$-module $M$ est de type fini et poss\`ede la propri\'et\'e
i). D'autre part, si $a(I_0+fA)^{k-1}=0$, alors $af^k\in aI_0\cdot (I_0+fA)^{k-1}=0$.

R\'eciproquement, soient $m_1\ld m_s$ des g\'en\'erateurs de $M$. i) nous
fournit un syst\`eme lin\'eaire~:
$$fm_i=\summ_{j=1\ld s}b_{ij}m_j,\quad i=1\ld s\quad\hbox{o\`u}\quad b_{ij}\in
I,\quad i,j=1\ld s$$
et pour tout $i=1\ld s$, on a~:
$$\left\vert\matrix{
b_{11}-f & b_{12}&\cdots &b_{1s}\cr
\noalign{\vskip3pt}
\vdots &&&\vdots\cr
\noalign{\vskip3pt}
b_{s1}&\cdots&\cdots &b_{ss}-f\cr}\right\vert m_i=0~.$$
Le d\'eterminant annule donc $M$. ii) nous permet d'obtenir la relation de
d\'ependance int\'egrale cherch\'ee.

\rque 1.10. Remarque|Si $I$ est de type fini et contient un \'el\'ement non
diviseur de 0, la condition ii) de 1.9 peut \^etre remplac\'ee par~:

ii') $M$ est un $A$-module fid\`ele ($aM=0\Leftrightarrow a=0$).

\noindent En effet d'une part ii') entra\^{\i}ne ii). D'autre part, si $f$ est
entier sur $I$, on peut choisir $M=(I+fA)^{k-1}$, qui est un $A$-module
fid\`ele puisqu'il contient l'\'el\'ement non diviseur de 0 de $I$.

\th 1.11. Corollaire|Soient $I$ et $J$ des id\'eaux de $A$. Alors~:
$$\oi\cdot\oj\subset\ol{I\cdot J}~.$$|

\dem Il suffit de montrer que si $f\in\oi$ et $g\in\oj$ alors
$fg\in\overline{I\cdot J}$. Comme dans le d\'ebut de la d\'emonstration de
1.9, choisissons pour modules v\'erifiant les conditions i) et ii) de 1.9
relativement \`a $f$ et $g$ des id\'eaux $M$ et $N$ de type fini de $A$. Et,
soit $R=MN$. $R$ est un id\'eal de type fini de $A$. De plus~:

\decale{i)} $fgR=fgMN=fM\cdot gN\subset IJMN=IJR$

\decale{ii)} Si $aR=0$, d\'esignant par $m_1\ld m_s$ un syst\`eme de
g\'en\'erateurs de $M$ on obtient que $am_iN=0$, $i=1\ld s$. Il existe donc
$k_i$, $i=1\ld s$, tels que $am_i g^{k_i}=0$ et si $k=\sup k_i$, $ag^kM=0$. Il
existe alors $\ell$, tel que $ag^kf^\ell=0$ et $a(fg)^{\sup k,\ell}=0$.

\th 1.12. Lemme|Soient $A$ un anneau normal \footnote{(*)}{un anneau est dit
normal s'il est r\'eduit et int\'egralement clos dans son anneau total de
fractions.} et $f$ un \'el\'ement non diviseur de z\'ero dans $A$~; alors $fA$
est un id\'eal int\'egralement clos.|

\dem Soit $g$ entier sur $fA$. Une relation (1.1.1) s'\'ecrit~:
$$g^k+b_1fg^{k-1}+\cdots +b_kf^k=0~.$$
L'\'el\'ement $g/f$ de $\tot A$ est donc entier sur $A$. Ainsi, il
appartient en fait \`a $A$ et $g$ appartient \`a $fA$.

\th 1.13. Lemme|Soit $A$ un anneau n\oe th\'erien. Pour un id\'eal $I$ de $A$,
les conditions suivantes sont \'equivalentes~:

\decale{i)} $\opp_{n\in\nnf}\overline{I^n} T^n[\resp\opp_{n\in\nnf}(\oi)^nT^n]$ est un
$\opp_{n\in\nnf}I^nT^n$-module de type fini.

\decale{ii)} Il existe un entier $N$ tel que si $n\ge N$ on ait l'\'egalit\'e
$I\cdot\overline{I^n}=\overline{I^{n+1}}$ $[\resp I\cdot(\oi)^n=(\oi)^{n+1}]$.|

\dem i) $\Rightarrow$ ii). Soit $E_1\ld E_s$ un syst\`eme de g\'en\'erateurs de
$\oplus \overline{I^n} T^n$. On peut supposer les $E_i$ homog\`enes. Posons $n_i=\deg
E_i$, $i=1\ld s$~; et $N=\sup n_i$. Soit $n\ge N$ et soit
$f\in\overline{I^{n+1}}$. Alors~:
$$fT^{n+1}=\summ_{i=1}^s A_i E_i\quad A_i\in\oplus I^k T^k$$
et on peut supposer que $A_i$ est homog\`ene de degr\'e $n+1-n_i$. Ceci
entra\^{\i}ne que~:
$$f\in\summ_{i=1}^s I^{n+1-n_i}\overline{I^{n_i}}\subset I. \overline{I^n}$$
en utilisant (1.8, ii)).

ii) $\Rightarrow$ i): L'hypoth\`ese permet d'\'ecrire que~:
$$\oplus_{n\in\nnf}\overline{I^n}=\oplus_{n\le N}\overline{I^n}+\ol{I^N} \pcc(I),$$
resp:
$$\oplus_{n\in\nnf}(\overline{I})^n=\oplus_{n\le N}(\overline{I})^n+(\ol{I})^N \pcc(I).$$
L'anneau $A$ \'etant n\oe th\'erien, chaque $\overline{I^n}$ (resp. $(\overline{I})^n$)  est donc un $A$-module de type
fini et a fortiori un $\pcc(I)$-module de type fini. $\opp_{n\in\nnf}\overline{I^n}T^n$ (resp. $\opp_{n\in\nnf}(\overline{I})^nT^n$) est
donc lui-m\^eme un $\pcc(I)$-module de type fini.

\th 1.14. Proposition|Soient $A$ un anneau excellent r\'eduit
\footnote{(*)}{Voir [2] pour la notion d'anneau excellent. Un anneau local
complet est un anneau excellent, de m\^eme que les anneaux locaux de la
g\'eom\'etrie analytique.} et $I$ un id\'eal contenant un \'el\'ement non
diviseur de z\'ero dans $A$. Alors il existe un entier $N$ tel que si $n\ge N$

1) $I\cdot\overline{I^n}=\overline{I^{n+1}}$

2) $I\cdot(\oi)^n=(\oi)^{n+1}$~.|

\dem $A$ \'etant un anneau n\oe th\'erien, les assertions 1) et 2) sont
respectivement \'equivalentes \`a 1')  $\opp_{n\in\nnf}\overline{I^n} T^n$ est un $\pcc(I)$-module de type fini et 2') $\pcc (\oi)$ est un $\pcc(I)$-module de type
fini; voir (1.13).
De plus, $\pcc(I)$ \'etant alors lui-m\^eme un anneau n\oe th\'erien, d'apr\`es
1.8 i), 1') entra\^{\i}ne 2'). En effet, $\opp_{n\in\nnf}(\overline{I})^n T^n$ est un sous
$\pcc(I)$-module de $\opp_{n\in\nnf}\overline{I^n} T^n$, il est donc lui-m\^eme de type fini si ce dernier l'est. On remarque aussi que $\pcc(I)$ a m\^eme anneau
total de fractions que $A[T]$. En effet, si $H\in A[T]$ et si $g$ est
l'\'el\'ement de $I$ non diviseur de z\'ero dans $A$, $g^{\deg H}\cdot H\in
\pcc(I)$. Par suite si ${F(T)\over G(T)}\in\tot A[T]$, ${F(T)\over
G(T)}={g^mF(T)\over g^mG(T)}\in\tot \pcc(I)$ si $m\ge\sup (\deg F,
 \deg G)$.
Utilisant maintenant que $A$ est excellent, la fermeture int\'egrale $\opp_{n\in\nnf}\overline{I^n} T^n$ de
$\pcc(I)$, (qui est une $A$-alg\`ebre de type fini) dans son anneau total de
fractions n'est donc
rien d'autre que sa fermeture int\'egrale dans l'anneau total de
fractions de $A[T]$, et est un $\pcc(I)$-module de type fini d'apr\`es [2], 2e partie, 7.8. 

\th 1.15. Proposition|Soit $f$ un \'el\'ement entier sur $I$. Alors
$\bnu_I(f)\ge 1$.|

\dem Soit $f^k+a_1 f^{k-1}+\cdots + a_k=0$, $\nu _I(a_i)\ge i$, $i=1\ld k$, une
relation de d\'ependance int\'egrale de $f$ sur $I$. $\bnu_I$ \'etant une
fonction d'ordre~:
$$\bnu_I(f^k)=k\bnu_I(f)\ge\inf_{i=1\cdots
k}\bnu_I(a_i)+\bnu_I(f^{k-i})=\inf_{i=1\cdots s}\bnu_I(a_i)+(k-i)\bnu_I(f)~.$$
Or $\bnu_I(a_i)\ge\nu _I(a_i)\ge i$. On en d\'eduit que~:
$$k\bnu_I(f)\ge\inf_{i=1\cdots k}i+(k-i)\bnu_I(f)~.$$
Soit $i_0$ l'indice r\'ealisant cet $\inf$
$$k\bnu_I(f)\ge i_0 +(k-i_0)\bnu_I(f)~.$$
Ceci montre que $\bnu_I(f)\ge 1$.

\th 1.16. Corollaire|Si $J$ est un id\'eal de type fini et si $J\subset\oi$, alors
$\bnu_I(J)\ge 1$.|

\dem En effet si $x_1\ld x_n$ est un syst\`eme de g\'en\'erateurs de $J$, on
sait que $\bnu_I(J)=\inff_{1\le i\le n}\bnu_I(x_i)$.

\rque 1.17. Remarque|Nous montrerons la r\'eciproque de 1.15 si $A$ est une
alg\`ebre analytique locale ou un anneau local complet ou le localis\'e d'une
$\cb$-alg\`ebre de type fini.

\th 1.18. Proposition|Soit $A$ un anneau local n\oe th\'erien. Si $I_1$ et
$I_2$ sont deux id\'eaux primaires pour l'id\'eal maximal de $A$ ayant m\^eme
cl\^oture int\'egrale, ils ont m\^eme multiplicit\'e.|

\dem Supposons d'abord $I_2\subset\ol{I_1}$. Alors d'apr\`es 1.16,
$\bnu_{I_1}(I_2)\ge 1$. Choisissons $\varepsilon >0$. Il existe $N$ tel que si
$n\ge N$
$$\nu _{I_1}(I^n_2)\ge n(1-\varepsilon )~.$$
Soit $m$ le plus petit entier sup\'erieur ou \'egal \`a $n(1-\varepsilon )$. On a
$$\nu _{I_1}(I^n_2)\ge m\quad\hbox{et}\quad I^n_2\subset I^m_1~.$$
De la d\'efinition de la multiplicit\'e d'un id\'eal primaire, on d\'eduit
imm\'ediatement que, d\'esignant par $e(I_2)$ $\big(\resp e(I_1)\big)$ la
multiplicit\'e de $I_2$ $(\resp I_1)$ et $d$ la dimension de $A$
$$\eqalign{
e(I^n_2)=n^d e(I_2)&\ge m^d e(I_1)=e(I^m_1)\cr
{e(I_2)\over e(I_1)}&\ge ({m\over n})^d~.\cr}\leqno{\rm et~que}$$
Or $n(1-\varepsilon )\le m$. Par suite~:
$${e(I_2)\over e(I_1)}\ge (1-\varepsilon )^d~.$$
Ceci \'etant vrai pour tout $\varepsilon >0$, $e(I_2)\ge e(I_1)$. On obtient
l'in\'egalit\'e oppos\'ee en utilisant $I_1\subset \oi_2$.

\rque 1.19. Remarque|Si $A$ est une alg\`ebre analytique locale
\'equi-dimensionnelle, D.~Rees [3] a montr\'e que r\'eciproquement si $I_1$ et
$I_2$ sont des id\'eaux primaires pour l'id\'eal maximal ayant m\^eme
multiplicit\'e et tels que $I_1\subset I_2$, alors $I_1$ et $I_2$ ont m\^eme
cl\^oture int\'egrale.

\titre R\'ef\'erences|

\div 1|Bourbaki n|Alg\`ebre commutative. Chapitres 5 et 6|{\rm~Hermann}|

\div 2|Grothendieck a|\'El\'ements de g\'eom\'etrie alg\'ebrique, {\bf
IV}|{\rm~Publications de l'IHES, PUF}|

\livre 3|Rees d|$a$-transforms of local rings and a theorem on multiplicities
of ideals|Proceedings Camb. Philos., {\bf 57, 1}|{\rm 8--17}|

\div 4|Zariski o. {\rm et} Samuel p|Commutative algebra, Vol.I, Chap. V et Vol. II  {\rm~Appendice 4\ } Van Nostrand| 1960|
\par\medskip\noindent
\titre 2. Les avatars de la cl\^oture int\'egrale d'un id\'eal en g\'eom\'etrie
analytique complexe|

Les anneaux qui vont nous int\'eresser maintenant sont les $\cb$-alg\`ebres
analytiques locales \ie, celles obtenues comme quotient d'un anneau de s\'eries
convergentes $\cb\{x_1\ld x_n\}$. Nous nous pr\'eparons \`a d\'ecrire la
filtration associ\'ee \`a la fonction d'ordre $\bnu_I$. Nous emploierons
syst\'ematiquement le langage g\'eom\'etrique.

\th 2.1. Th\'eor\`eme|Soient $X$ un espace analytique complexe r\'eduit, $Y$ un
sous-espace analytique ferm\'e rare de $X$, $x$ un point de $Y$. Soit $\ic$
l'id\'eal coh\'erent de $\oc_X$ d\'efinissant $Y$. Soit $\jc$ un
$\oc_X$-id\'eal coh\'erent. Soit $I (\resp J)$ le germe de $\ic (\resp\jc)$ en
$x$. Soit $\oi$ la cl\^oture int\'egrale de $I$ dans $\oc_{X,x}$. Les
conditions suivantes sont \'equivalentes~:

\decale{i)} $J\subset\oi$.

\decale{ii)} $\bnu_I(J)\ge 1$.

\decale{iii)} Soit $\db$ le disque unit\'e du plan complexe $(\db=\{t\in\cb,
|t|<1\})$. Pour tout germe de morphisme $h:(\db,0)\to (X,x)$
$$h^*J\cdot\oc_{\db,0}\subset h^*I\cdot\oc_{\db,0}~.$$

\decale{iv)} Pour tout morphisme $\pi :X'\to X$ tel que 1) $\pi $ soit propre
et surjectif, 2) $X'$ soit un espace analytique normal, 3) $\ic\cdot\oc_{X'}$
soit un $\oc_{X'}$-module inversible, il existe un ouvert $U$ de $X$ contenant
$x$ tel que~:
$$\jc\cdot\oc_{X'|\pi ^{-1}(U)}\subset\ic\cdot\oc_{X'|\pi ^{-1}(U)}~.$$

\decale{v)} Il existe un $\oc_X$-id\'eal coh\'erent $\kc$ dont le support est
rare dans $X$ tel que si $\pi :\wtx\to X$ est l'\'eclatement de $\kc$, il
existe un ouvert $U$ de $X$ contenant $x$ tel que~:
$$\jc\cdot\oc_{\wtx |\pi ^{-1}(U)}\subset\ic\cdot\oc_{\wtx |\pi ^{-1}(U)}~.$$

\decale {vi)} Soit $V$ un voisinage de $x$ sur lequel $\jc$ et $\ic$ sont
engendr\'es par leurs sections globales. Pour tout syst\`eme de g\'en\'erateurs
$g_1\ld g_m$ de $\Gamma (V,\ic)$ et tout \'el\'ement $f$ de $\Gamma (V,\jc)$,
on peut trouver un voisinage $V'$ de $x$ et une constante $C$ tels que~:
$$|f(y)|\le C\sup_{i=1\ld m} |g_i(y)|,\hbox{\rm pour tout} \  y\in V'~.$$|

Avant de donner la d\'emonstration de ce th\'eor\`eme, nous allons \'enoncer et
d\'emontrer quelques lemmes que nous aurons \`a utiliser.

\th 2.1.1. Lemme|Soit $A$ un anneau local n\oe th\'erien, $I=(g_1\ld g_p)$ un
id\'eal $\ne 0$ de $A$ qui est principal. Alors $I$ est engendr\'e par l'un des
$g_i$.|

\dem Soit $g$ un g\'en\'erateur de $I$. Il existe $a_1\ld a_p$ et $b_1\ld b_p$
dans $A$ tels que~:
$$g_i=a_i g, i=1\ld p\quad\hbox{et}\quad g=\summ_{i=1\ld p}b_ig_i~.$$
Alors~:
$$g\big(1-\summ_{i=1\ld p}a_ib_i\big)=0~.$$
Il s'agit de voir que l'un des $a_i$ au moins est une unit\'e dans $A$. S'il
n'en \'etait pas ainsi, $1-\summ_{i=1\ld p}a_ib_i$ serait une unit\'e et $g$
serait nul.

\th 2.1.2. Lemme|Soient $A$ un anneau n\oe th\'erien de caract\'eristique $0$,
$I$ un id\'eal contenant au moins un \'el\'ement non diviseur de z\'ero. Alors
$I$ admet un syst\`eme de g\'en\'erateurs form\'es d'\'el\'ements non diviseurs
de z\'ero.|

\dem Soit $h_1$ l'\'el\'ement de $I$ non diviseur de z\'ero et soit $(h_1\ld
h_n)$ un syst\`eme de g\'en\'erateurs de $I$. Supposons que $h_1\ld h_k$
soient non diviseurs de z\'ero. Nous allons montrer que nous pouvons modifier
$h_{k+1}$.

Consid\'erons~:
$$g_s=h_{k+1}+sh_k,\quad s\in\nnf~.$$
On sait qu'un anneau n\oe th\'erien poss\`ede un nombre fini d'id\'eaux
premiers minimaux et qu'un \'el\'ement est diviseur de z\'ero si et seulement
s'il appartient \`a l'un d'entre eux. $A$ \'etant de caract\'eristique $0$, si
tous les $g_s$, $s\in\nnf$, \'etaient diviseurs de z\'ero, d'apr\`es le
principe des tiroirs, on pourrait d\'eterminer 2 entiers $s_1$ et $s_2$ et un
id\'eal premier minimal $P$ tels que $g_{s_1}$ et $g_{s_2}$ appartiennent \`a
$P$. On en d\'eduirait que $(s_1-s_2)h_k$ et donc $h_k$ lui-m\^eme est dans
$P$, ce qui est impossible puisque $h_k$ est non diviseur de z\'ero. Il existe
donc $s_0\in\nnf$ tel que $g_{s_0}=h_{k+1}+s_0h_k$ est non diviseur de z\'ero.
On remplace $h_{k+1}$ par $g_{s_0}$ et l'on construit ainsi par r\'ecurrence
le syst\`eme de g\'en\'erateurs d\'esir\'e.

\th 2.1.3. Lemme|Soient $X$ un espace analytique r\'eduit et $x$ un point de $X$.
On suppose que $X$ est normal au voisinage de $x$. Soit $f$ un germe de
fonction m\'eromorphe, non holomorphe, au voisinage de $x$. Alors il existe un
germe de morphisme $h:(\db,0)\to(X,x)$ o\`u $\db$ d\'esigne le disque unit\'e
du plan complexe tel que $f\circ h$ soit un germe de fonction m\'eromorphe,
non holomorphe, \`a l'origine de~$\db$.|

\dem Soit $\pi :\wtx\to X$ une r\'esolution des singularit\'es de $X$.
$f\circ\pi $ est un germe de fonction m\'eromorphe au voisinage de $\pi
^{-1}(x)$ et il existe un point $\tx\in\pi ^{-1}(x)$ tel que $f\circ\pi $ ne
soit pas holomorphe au voisinage de $\tx$. En effet s'il n'en \'etait pas
ainsi, $f\circ \pi $ serait born\'ee sur un voisinage de $\pi ^{-1}(x)$ et $\pi
$ \'etant propre, $f$ elle-m\^eme serait born\'ee au voisinage de $x$. $X$
\'etant normal en $x$, $f$ m\'eromorphe et born\'ee serait holomorphe. Nous
sommes donc ramen\'es \`a montrer 2.1.3 en supposant en plus que $X$ est non
singulier au voisinage de $x$. L'anneau local de $X$ en $x$ est alors
factoriel et le germe de fonction m\'eromorphe $f$ consid\'er\'e a un
repr\'esentant ${p\over q}$ o\`u $p$ et $q$ sont holomorphes au voisinage de
$x$ sans facteurs irr\'eductibles en commun. L'ensemble des z\'eros de $q$
contient donc une hypersurface $H$ au voisinage de $x$ non contenue enti\`erement
dans l'ensemble des z\'eros de $p$. On peut alors trouver un germe de courbe
irr\'eductible (en g\'en\'eral singulier en $x$) $\Gamma $ tel que $\Gamma $
soit contenu dans $H$ et non dans les z\'eros de $p$. La normalisation de
$\Gamma $ nous fournit un germe de morphisme $\tih :(\db,0)\to(X,x)$ tel que,
identifiant $\oc_{\db,0}$ \`a $\cb\{t\}$ anneau des s\'eries convergentes \`a
une variable et d\'esignant par $v$ la valuation naturelle sur $\cb\{t\}$,
$$v(p\circ\tih)=\alpha \quad\hbox{et}\quad v(q\circ\tih)=\infty~.$$
Nous allons montrer que modifiant $\tih$ par des termes en $t$ d'ordre assez
grand on peut trouver $h:(\db,0)\to(X,x)$ tel que~:
$$v(p\circ h)=\alpha \quad\hbox{et}\quad v(q\circ h)=\beta,\quad \beta
>\alpha ~.$$  
En effet soit $(x_1\ld x_n)$ un syst\`eme de coordonn\'ees sur $X$ au voisinage
de $x$ et posons $\tih_i(t)=x_i\circ\tih$. Alors,
$$\leqalignno{
\ \ \ \ p\big(\tih_1(t)+t^N,\tih_2(t)\ld\tih_n(t)\big)-p\big(\tih_1(t)\ld\tih_n(t)\big)&=
t^NR(t) \hbox{~o\`u~}R\in\cb\{t\},&(*)\cr
\ \ \ \ q\big(\tih_1(t)+t^N,\tih_2(t)\ld\tih_n(t)\big)-q\big(\tih_1(t)\ld\tih_n(t)\big)&=
t^NS(t) \hbox{~o\`u~}S\in\cb\{t\}.&(**)\cr}$$
Si donc $N>\alpha $, on d\'eduit de $(*)$ que 
$v\Big(p\big(\tih_1(t)+t^N,\tih_2(t)\ld\tih_n(t)\big)\Big)=\alpha $ et de $(**)$
que $v\Big(q\big(\tih_1(t)+t^N,\tih_2(t)\ld\tih_n(t)\big)\Big)\ge N>\alpha $.
On peut donc choisir $h_1(t)=\tih_1(t)+t^{\alpha +1}$, $h_i(t)=\tih_i(t)$, $i\ge
2$.

\rque D\'emonstration du th\'eor\`eme 2.1|Il suffit de montrer 2.1 si $J$ est
principal. Notons $f$ son g\'en\'erateur.

i) $\Rightarrow$ ii) Voir proposition 1.15.

ii) $\Rightarrow$ iii) Soit $\varphi :\oc_{X,x}\to\cb\{t\}$ le morphisme
associ\'e \`a $h$ et soit $v$ la valuation naturelle sur $\cb\{t\}$. Il s'agit
de voir que $\varphi (f)\in I\cdot\cb\{t\}$. Or, il existe un entier $m\ge 1$
tel que $I\cdot\cb\{t\}=t^m\cdot\cb\{t\}$ et il suffit en fait que~:
$$v\big(\varphi (f)\big)\ge m~.$$
Mais $\bnu_{I\cdot\cb\{t\}}\big(\varphi (f)\big)\ge\bnu_I(f)\ge 1$ d'apr\`es la
d\'efinition de $\bnu$ et on a vu que~:
$$\bnu_{I\cdot\cb\{t\}}\big(\varphi (f)\big)=
\bnu_{t^m\cdot\cb\{t\}}\big(\varphi (f)\big)=
{1\over m}\bnu_{t\cdot\cb\{t\}}\big(\varphi (f)\big)={1\over m}v\big(\varphi
(f)\big)~.$$

iii) $\Rightarrow$ iv) Nous utilisons ici le lemme 2.1.3. Le morphisme $\pi $ \'etant
propre, il s'agit en fait de montrer que pour tout $x'\in\pi ^{-1}(x)$
$$f\cdot\oc_{X',x'}\in I\cdot\oc_{X',x'}~.$$
Supposons qu'il n'en soit pas ainsi en $x'$. Soit $g$ un g\'en\'erateur de
$I\cdot\oc_{X',x'}$. D'apr\`es l'hypoth\`ese 3), $g$ est non diviseur de z\'ero
et ${f\over g}$ est un \'el\'ement de $\tot\oc_{X',x'}$ qui n'est pas dans
$\oc_{X',x'}$.

Il existe alors un germe de morphisme $h':(\db,0)\to(X',x')$ tel que, $v$
d\'esignant la valuation naturelle sur $\cb\{t\}$, si $\varphi
':\oc_{X',x'}\to\cb\{t\}$ est le morphisme associ\'e \`a $h'$, on ait $v\big(\varphi
'(f)\big)<v\big(\varphi '(g)\big)<\infty$.

Nous obtenons la contradiction cherch\'ee avec le morphisme $h=\pi \circ
h':(\db,0)\to(X,x)$.

iv) $\Rightarrow$ v) Soit $\pi :\oxp\to X$ l'\'eclatement normalis\'e
de $Y$. $Y$ \'etant rare dans $X$, $\pi $ satisfait les conditions 1), 2) et 3)
de iv). Or, \'etant donn\'e un espace analytique $\xc$ r\'eduit, il existe un
id\'eal coh\'erent sur cet espace dont l'\'eclatement est la normalisation de
$\xc$. De plus, le compos\'e de 2 \'eclatements est un \'eclatement (non
canoniquement).

v) $\Rightarrow$ i) L'id\'eal $\kc$ d\'efinissant un sous espace rare dans $X$, en tout point son germe
contient un \'el\'ement non diviseur de z\'ero. Utilisant le lemme 2.1.2 et la
coh\'erence de $\kc$ on peut supposer que $U$ est assez petit pour que $\kc |U$
soit engendr\'e par un nombre fini de sections globales $g_1\ld g_m$ sur $U$, 
chaque $g_i$ \'etant non diviseur de z\'ero. On peut supposer \'egalement que $U$ est assez petit
pour que $\ic |U$ soit engendr\'e par ses sections globales. Alors $\pi
^{-1}(U)$ est recouvert par un nombre fini d'ouverts $V_i$ tels que~:
$$\wtx |V_i\simeq\spe\oc_{X|U}\big[{g_1\over g_i}\ld {\hat{ g_i}\over g_i}\ld
{g_m\over g_i}\big]~.$$
(On utilise ici la convention habituelle, l'\'el\'ement sous le $\wedge$ est
omis). Sur chaque ouvert $V_i$, $f$ s'exprime donc polynomialement en fonction
des $g_k/g_i$, $k\ne i$. $K$ d\'esignant le germe en $x$ de $\kc$, on
d\'etermine un entier $n_i\ge 1$ tel que~:
$$f\cdot g^{n_i}_i\in I\cdot K^{n_i}~.$$
Posons $n=\summ_{i=1\ld m}n_i$ et soit $fg^{\alpha _1}_1\cdots g^{\alpha
_n}_n$. Si $\alpha _1+\cdots +\alpha _n\ge n$, il existe certainement $i_0$
tel que $\alpha _{i_0}\ge n_{i_0}$. On peut alors \'ecrire que~:
$$fg^{\alpha _1}_1\cdots g^{\alpha _n}_n=fg_{i_0}^{n_{i_0}}\cdot
g_{i_0}^{\alpha _{i_0}-n_{i_0}}\prod_{j\ne i_0}g_j^{\alpha _j}\in I\cdot
K^{n_{i_0}}K^{n-n_{i_0}}=IK^n~.$$
Ceci montre que~:
$$f\cdot K^n\subset I\cdot K^n~.$$
$K$ contenant un \'el\'ement de $\oc_{X,x}$ non diviseur de z\'ero, $K^n$ est
un $\oc_{X,x}$-module de type fini fid\`ele. D'apr\`es 1.10, ceci suffit \`a
assurer que $f\in\oi$.

iv) $\leftrightarrow$ vi) Soit $g_1\ld g_m$ un syst\`eme de g\'en\'erateurs de
$I$ form\'e d'\'el\'ements non diviseurs de z\'ero. Le morphisme $\pi $ \'etant
propre, on peut recouvrir $\pi ^{-1}(x)$ par un nombre fini d'ouverts
$V_\alpha $ relativement compacts sur chacun desquels un des $g_i$ engendre
$\ic\cdot \oc_{X'}$ (lemme 2.1.1). Le morphisme $\pi $ $\pi $ \'etant surjectif, la condition vi)
est satisfaite si et seulement si sur chaque $V_\alpha $ le quotient ${|f\circ\pi
|\over\sup |g_i\circ\pi |}$ est born\'e. Or si $g_1$ est le g\'en\'erateur de
$\ic\cdot\oc_{X'}$ sur $V_\alpha $, cette derni\`ere fonction est born\'ee, si
et seulement si ${|f\circ\pi |\over |g_1\circ\pi |}$ l'est. Mais en tout point
$y$ de $V_\alpha $, le germe de $g_1$ est non diviseur de z\'ero dans
$\oc_{X',y}$ et par cons\'equent ${f\circ\pi \over g_1\circ\pi }$ induit un \'el\'ement de
$\tot\oc_{X',y}$ et $X'$ \'etant un espace normal, on sait que le fait que cet \'el\'ement soit born\'e \'equivaut
au fait que ${f\circ\pi \over g_1\circ\pi }$ appartient \`a $\oc_{X',y}$ pour
tout $y$ de $V_\alpha $, \ie, $f\cdot\oc_{X'|V_\alpha
}\in\ic\cdot\oc_{X'|V_\alpha }$.

\th 2.2. Corollaire|M\^emes hypoth\`eses qu'au th\'eor\`eme 2.1. Soit $k$ un
entier $\ge 1$. Les conditions suivantes sont \'equivalentes~:

\decale{i)} $J\subset\overline{I^k}$.

\decale{ii)} $\bnu_I(J)\ge k$.| 

\dem C'est imm\'ediat puisque $\bnu_{I^k}(J)={1\over k}\bnu_I(J)$.

\th 2.3. Corollaire|M\^emes hypoth\`eses qu'au th\'eor\`eme 2.1. Les conditions
suivantes sont \'equivalentes~:

\decale{i)} $\bnu_I(f)=\infty$.

\decale{ii)} $f=0$.|

\dem Il suffit de montrer que i) $\Rightarrow$ ii). Or, dire que
$\bnu_I(f)=\infty$ signifie que $\bnu_I(f)\ge k$ pour tout entier $k$.
D'apr\`es 2.2, ceci entra\^{\i}ne que~:
$$f\in\cap_{k\in\nnf}\overline{I^k}~.$$
De plus, d'apr\`es 1.14, puisque une alg\`ebre analytique locale est un anneau
excellent il existe $N$ tel que si $k\ge N$, $\overline{I^k}=I^{k-N}\cdot\overline{I^N}$. Ceci
entra\^{\i}ne que~:
$$f\in\cap_{k\ge N}I^{k-N}=0~.$$

\th 2.4. Corollaire|M\^emes hypoth\`eses qu'au th\'eor\`eme 2.1. La topologie
de $\oc_{X,x}$ associ\'ee \`a la fonction d'ordre $\bnu_I$ est la topologie
$I$-adique.|

\dem D'apr\`es 2.2, la topologie de $\oc_{X,x}$ associ\'ee \`a la fonction
d'ordre $\bnu_I$ est la topologie associ\'ee \`a la filtration de $\oc_{X,x}$
par les id\'eaux $\overline{I^k}$. Cette filtration est $I$-bonne d'apr\`es 1.14. Elle
d\'efinit donc la topologie $I$-adique [1] \S3.

\th 2.5. Corollaire|M\^emes hypoth\`eses qu'au th\'eor\`eme 2.1. Soit
$k\in\nnf$. Si $\bnu_{\ic_x}(f)\ge k$, il existe un voisinage $U$ de $x$ dans
$X$ tel que $\bnu_{\ic_y}(f)\ge k$ pour tout $y\in U$.|

\dem D'apr\`es 2.2, $f$ v\'erifie une relation de d\'ependance int\'e\-grale~:
$$f^k+\Sigma a_i f^{k-i}=0, \nu _{\ic_x}(a_i)\ge ik~.$$
Il existe un voisinage $U$ de $x$ dans $X$ tel que $\nu _{\ic_y}(a_i)\ge ik$,
si $y\in U$. Ceci montre que $f_y\in\ol{\ic^k_y}$ et donc, toujours d'apr\`es 2.2
que $\bnu_{\ic_y}(f)\ge k$, si $y\in U$.

\th 2.6. Th\'eor\`eme|Soit $X$ un espace analytique complexe r\'eduit, $Y$ un
sous-espace analytique ferm\'e rare de $X$, $\ic$ l'id\'eal coh\'erent de
$\oc_X$ d\'efinissant $Y$. Il existe un $\oc_X$-id\'eal coh\'erent not\'e
$\oic$ tel que pour tout point $y$ de $X$~:
$$(\oic)_y=\ol{\ic_y}=\{f\in\oc_{X,y},\bnu_{\ic_y}(f)\ge 1\}~.$$
On l'appelle la cl\^oture int\'egrale de $\ic$.|

\dem Soit $\oic$ l'id\'eal de $\oc_X$ dont les sections sur un ouvert
quelconque $U$ de $X$ sont donn\'ees par~:
$$\Gamma (U,\oic)=\{f\in\Gamma (U,\oc_X):\bnu_{\ic_y}(f)\ge 1,\quad \forall
y\in U\}~.$$
Il est imm\'ediat que $\oic$ ainsi d\'efini est un faisceau d'id\'eaux de
$\oc_X$ et que, d'apr\`es 2.5, son germe en $y$ est $\ol{\ic_y}$. Montrons que
c'est un id\'eal coh\'erent. Soit $\pi :\oxp\to X$ l'\'eclatement normalis\'e
de $\ic$. On a le diagramme (dans la cat\'egorie des $\oc_X$-modules)
$$\matrix{
\ic &\kern-5pt \la \pi _*(\pi ^*\ic)\la &\kern-5pt\pi _*(\ic\oc_{\oxp})\cr
\noalign{\vskip3pt}
\big\downarrow &&\big\downarrow \cr
\noalign{\vskip3pt}
\oc_X &\kern-5pt \hfll &\kern-5pt \pi _*(\oc_{\oxp})\cr}$$
Les fl\`eches dispos\'ees aux 4 c\^ot\'es du carr\'e sont des injections. En
effet, $X$ est r\'eduit et $\pi $ surjectif, et $\pi _*$ est exact \`a gauche.
L'\'equivalence de i), iv) et v) dans 2.1 montre que~:
$$\oic=\pi _*(\ic\oc_{\oxp})\cap\oc_X~.$$
Mais d'apr\`es le th\'eor\`eme de Grauert, l'image directe d'un module
coh\'erent en est un et l'intersection de 2 sous-modules coh\'erents d'un
module coh\'erent est lui-m\^eme un module coh\'erent.

\th 2.7. Corollaire|Les hypoth\`eses sont celles du th\'eor\`eme 2.6. Soit
$k\in\nnf$. Il existe un $\oc_X$-id\'eal coh\'erent not\'e $\ol{\ic^k}$ tel que
pour tout $y\in X$
$$(\ol{\ic^k})_y=\ol{\ic^k_y}=\{f\in\oc_{X,y},\bnu_{\ic_y}(f)\ge k\}~.$$|

\dem C'est la cl\^oture int\'egrale de $\ic^k$.

\th 2.8. Proposition|Les hypoth\`eses sont celles du th\'eor\`eme 2.6.
$\opp_{n\in\nnf}\ol{\ic^n}$ est une $\oc_X$-alg\`ebre de pr\'esentation finie.
$\opp_{n\in\nnf}\ol{\ic^n}$ est un $\opp_{n\in\nnf}\ic^n$-module de type fini.|

Nous allons d'abord montrer le lemme suivant dont nous aurons besoin dans la
d\'emonstration et dans la suite~:

\th 2.8.1. Lemme|Soit $K$ un polycylindre. Il existe un entier $N$ (d\'ependant
uniquement de $K$ et de $\ic$) tel que si $n\ge N$, $n\in\nnf$, on a
$$\Gamma (K,\ol{\ic^n})=\Gamma (K,\ic)^{n-N}\cdot\Gamma (K,\ol{\ic^N})~.$$|

\dem 2.5 entra\^{\i}ne imm\'ediatement que pour tout polycylindre $K$,
$$\Gamma (K,\ol{\ic^n})=\{f\in\Gamma (K,\oc_X);\bnu_{{\cal{I}}_{y}}(f)\ge n,\forall y\in
K\}~.$$
Soit $\kc$ le $\oc_X$-id\'eal coh\'erent d\'efinissant un sous-espace rare dans $X$ et
dont l'\'eclatement est $\oc_X$-isomorphe \`a l'\'eclatement normalis\'e de $\cal{I}$
dans $X$. Soit $x$ un point de $K$ et soit $f_1\ld f_s$ un syst\`eme de
g\'en\'erateurs de $\ol{\ic^n_x}$. Puisque $\bnu_{\ic_x^n}(f_i)\ge 1$, d'apr\`es
2.1 (voir la partie de la d\'emonstration v) $\Rightarrow$ i), il existe
$p_i\in\nnf$ tel que~:
$$f_i\kc_x^{p_i}\subset\ic^n_x\cdot\kc_x^{p_i}~.$$
Soit $p_x=\sup p_i$
$$\ol{\ic^n_x}\cdot\kc_x^{p_x}\subset\ic^n_x\kc_x^{p_x}~.$$
$\ol{\ic^n}$, $\kc$ et $\ic^n$ \'etant coh\'erents, il existe $U_x$ un voisinage
de $x$ tel que~:
$$\ol{\ic^n_y}\cdot\kc_y^{p_x}\subset\ic^n_y\kc_y^{p_x}\quad\hbox{\rm pour tout}\  y\in U_x~.$$
Le polycylindre $K$ \'etant compact est recouvert par un nombre fini d'ouverts. D\'esignons les
$U_{x_1}\ld U_{x_t}$ et soit $p=\supp_{i=1\ld t}p_{x_i}$ et $U=\cupp_{i=1\ld
t}U_{x_i}$. L'ouvert $U$ est un voisinage de $K$ dans $X$.
$$\ol{\ic^n_y}\cdot\kc^p_y\subset\ic^n_y\cdot\kc^p_y,\quad\hbox{\rm pour tout}\  y\in U,$$
autrement dit
$$\ol{\ic^n}|U\cdot \kc^p|U\subset\ic^n|U\cdot\kc^p|U.$$
D'apr\`es [2] lemme 1.12 (c'est une cons\'equence
imm\'ediate du th\'eor\`eme B), on a 
$$\Gamma (K,\ol{\ic^n})\cdot\Gamma (K,\kc)^p\subset\Gamma (K,\ic)^n\cdot \Gamma
(K,\kc)^p~.$$
De plus, $\Gamma (K,\kc)^p$ est un module fid\`ele. En effet si $a\in\Gamma
(K,\kc)$ est tel que $a\Gamma (K,\kc)=0$, ceci entra\^{\i}ne que pour tout $x$
de $K$, on ait $a_x\kc_x=0$. Le support du sous espace d\'efini par $\kc$ \'etant rare, $\kc_x$ contient
certainement un \'el\'ement non diviseur de z\'ero de $\oc_{X,x}$. Ainsi
$a_x=0$ pour tout $x$ de $K$ et $a=0$.

D'apr\`es 1.9, $\Gamma (K,\ol{\ic^n})$ est donc contenu dans la cl\^oture
int\'egrale de $\Gamma (K,\ic)^n$. R\'eciproquement, il est facile de voir que
si $g\in\Gamma (K,\oc_X)$ est entier sur $\Gamma (K,\ic)^n$, alors $g\in\Gamma
(K,\ol{\ic^n})$, puisque la relation de d\'ependance int\'egrale dans $\Gamma
(K,\oc_X)$
$$g^k+\Sigma a_i g^{k-i}, \quad a_i\in\Gamma (K,\ic)^{ni}$$
induit pour tout $x$ de $K$, une relation de d\'ependance int\'egrale de $g_x$
sur $\ic^n_x$ donc que $\bnu_{\ic_x}(f)\ge n$. Nous avons donc montr\'e
finalement que $\Gamma (K,\ol{\ic^n})$ est la cl\^oture int\'egrale de $\Gamma
(K,\ic)^n$ dans $\Gamma (K,\oc_X)$.\par\noindent Appliquant maintenant le lemme 1.7, nous
obtenons que $\opp_{n\in\nnf}\Gamma (K,\ol{\ic^n})T^n$ est la fermeture
int\'egrale de ${\cal P}\big(\Gamma (K,\ic)\big)$ dans $\Gamma (K,\oc_X)[T]$. L'anneau $\Gamma
(K,\oc_X)$ \'etant excellent, $\opp_{n\in\nnf}\Gamma (K,\ol{\ic^n})$ est
 un ${\cal P}\big(\Gamma (K,\ic)\big)$-module de type fini et d'apr\`es 1.13, il
existe un entier $N$ tel que si $n\ge N$
$$\Gamma (K,\ic)\cdot\overline{\Gamma (K,\ic)^n}=\overline{\Gamma
(K,\ic)^{n+1}}$$
ou encore
$$\Gamma (K,\ic)\cdot\Gamma (K,\ol{\ic^n})=\Gamma (K,\ol{\ic^{n+1}})~.$$

\rque D\'emonstration de 2.8|Puisque $\opp_{n\in\nnf}\Gamma
(K,\ol{\ic^n})\big(\resp\opp_{n\in\nnf}\Gamma (K,\ic^n)\big)$ est canoniquement
isomorphe \`a $\opp_{n\in\nnf}\Gamma (K,\ol{\ic^n})T^n$ $\big(\resp {\cal P}(\Gamma
(K,\ic))\big)$, le lemme pr\'ec\'edent a montr\'e que $\opp_{n\in\nnf}\Gamma
(K,\ol{\ic^n})$ est un $\opp_{n\in\nnf}\Gamma (K,\ic^n)$-module de type fini. Ceci
signifie qu'il existe un entier $s$ et un morphisme surjectif, $\opp_{n\in\nnf}\Gamma
(K,\ic^n)$ lin\'eaire~:
$$\bigl(\opp_{n\in\nnf}\Gamma (K,\ic^n)\bigr)^s\longrightarrow\opp_{n\in\nnf}\Gamma
(K,\ol{\ic^n})~.\leqno(*)$$
Ceci permet de construire un morphisme de $\opp_{n\in\nnf}\ic^n|K$-modules
$$(\opp_{n\in\nnf}\ic^n|K)^s\la\opp_{n\in\nnf}\ol{\ic^n}|K$$
dont il s'agit de voir qu'il est surjectif. Pour qu'il en soit ainsi, il faut
et il suffit que pour tout $x\in K$, le morphisme
$$(\opp_{n\in\nnf}\ic^n_x)^s\la\opp_{n\in\nnf}\ol{\ic^n_x}$$
le soit.

Or, tensorisons $(*)$ au-dessus de $\Gamma (K,\oc_X)$ par $\oc_{X,x}$.
$$\opp_{n\in\nnf}\big(\Gamma (K,\ic^n)\otimes_{\Gamma (K,\oc_X)}\oc_{X,x}\big)^s\la
\opp_{n\in\nnf}\big(\Gamma (K,\ol{\ic^n})\big)\otimes_{\Gamma (K,\oc_X)}\oc_{X,x}$$
est aussi surjectif. Mais $\ic^n$ et $\ol{\ic^n}$ \'etant des $\oc_X$-id\'eaux
coh\'erents et $K$ \'etant un polycylindre, le th\'eor\`eme $A$ nous dit
justement que $\Gamma (K,\ic^n)\ott_{\Gamma
(K,\oc_X)}\oc_{X,x}=\ic^n_x$ et que $\Gamma (K,\ol{\ic^n})\ott_{\Gamma
(K,\oc_X)}\oc_{X,x}=\ol{\ic^n_x}$.

L'alg\`ebre $\opp_{n\in\nnf}\ol{\ic^n}$ est donc un $\opp_{n\in\nnf}\ic^n$-module de type fini. La $\oc_X$-alg\`ebre $\opp_{n\in\nnf}\ic^n$
\'etant  elle m\^eme de type fini (en tant qu'alg\`ebre), ceci
entra\^{\i}ne \`a son tour que $\opp_{n\in\nnf}\ol{\ic^n}$ est une $\oc_X$-alg\`ebre de
type fini. Puisque $\ol{\ic^n}$ est un $\oc_X$-module coh\'erent pour chaque $n$, $\opp_{n\in\nnf}\ol{\ic^n}$ est
une $\oc_X$-alg\`ebre de pr\'esentation finie.

\titre R\'ef\'erences|

\div 1|Bourbaki n|Alg\`ebre commutative, chapitres 3 et 4|{\rm~Hermann}|

\article 2|Frisch j|Points de platitude d'un morphisme d'espaces analytiques
complexes|Inventiones|4|1967|118--138|
\par\medskip\noindent

\titre 3. Cl\^oture int\'egrale d'un id\'eal et \'eclatement normalis\'e|

\rque 3.1. Remarque|Soit $A$ un anneau n\oe th\'erien r\'eduit, $\oa$ sa
normalisation, $I$ un id\'eal propre de $A$ et soit $J=I. \oa$
$$\ol{J^n}=\{f\in\oa\hbox{~v\'erifiant une relation~}f^k+\summ_{i=1}^k a_i
f^{k-i}=0,\quad a_i\in A,~\nu _I(a_i)\ge n i\}.$$

\dem Il suffit de remarquer que la fermeture int\'egrale de $\pcc(I)$ dans
$\oa[T]$ est aussi celle de $\pcc(J)$ dans $\oa[T]$. En effet  $JT=I. \oa T$
est form\'e d'\'el\'ements de $\oa[T]$ entiers sur $\pcc(I)$.

D'autre part, un calcul analogue \`a celui de 1.7 montre que la fermeture
int\'egrale de $\pcc(I)$ dans $\oa[T]$ est $\opp_{n\in\nnf}J_nT^n$ o\`u
$J_n=\{f\in\oa$ v\'erifiant une relation $f^k+\Sigma a_if^{k-i}=0$, $a_i\in A$,
$\nu _I(a_i)\ge ni\}$.

Mais on sait aussi (1.7) que la fermeture int\'egrale de $\pcc(J)$ dans $\oa[T]$
est $\oplus\ol{J^n}T^n$.

\th 3.2. Proposition|Soient $X$ un espace analytique complexe r\'eduit et $n:\ox\to X$ le morphisme de normalisation. Soit $Y$ un sous-espace
analytique ferm\'e rare de $X$, $W$ son image r\'eciproque par $n$. Soit
$\ic\ (\resp\jc)$ l'id\'eal coh\'erent de $\oc_X(\resp\oc_{\ox})$ d\'efinissant
$Y$ $(\resp W)$.

L'\'eclatement normalis\'e de $Y$ est $X$-isomorphe au morphisme compos\'e
$$\pro_{\ox}\oplus_{n\in\nnf}\ol{\jc^n}\build\la^\pi \fin\ox\build\la^n \fin
X$$ o\`u $\pi $ d\'esigne le morphisme structural.|

\dem Posons $Z=\pro\opp_{n\in\nnf}\ol{\jc^n}$ et soit
$p:X'\simeq\pro\opp_{n\in\nnf}\ic^n\to X$ l'\'eclatement de $Y$ dans $X$. Par
fonctorialit\'e de la formation du $\pro$, il existe $q:Z\to X'$ tel que $q\circ
p=\pi \circ n$. On sait d'apr\`es 2.8 que $\oplus\ol{\ic^n}$ est un
$\oplus\ic^n$-module de type fini. Apr\`es la remarque 3.1 le m\^eme
raisonnement montre que $\oplus\ol{\jc^n}$ est un $\oplus\ic^n$-module de type
fini. De [1] expos\'e 19, on d\'eduit que le morphisme canonique $\spe
\oplus\ol{\jc^n}\to\spe\oplus\ic^n$ est un morphisme fini et a fortiori $q$. Soit
$N(X)$ l'ouvert des points normaux de $X$ et soit $F=\complement N(X)\cup |Y|$.
C'est un ferm\'e analytique rare de $X$ dont l'image r\'eciproque $F'$ dans $X'$
est \'egalement un ferm\'e rare et $q$ induit un isomorphisme analytique de
$Z-(\pi \circ n)^{-1}(F)$ sur $X'-p^{-1}(F)$.

Montrons maintenant que $Z$ est un espace normal. Pour cela, il suffit que
$C_Z=\spe\oplus\ol{\jc^n}$ le soit, \ie, il suffit, que pour tout $x$ de $\ox$, le
germe en $x$ de $\oplus\ol{\jc^n}$ soit int\'egralement ferm\'e dans son anneau
total de fractions. Or posant $J=\jc_x$, $\oa=\oc_{\ox,x}$, nous avons d\'ej\`a
remarqu\'e (1.14) que $\oplus\ol{J^n}T^n$ (qui est canoniquement isomorphe \`a
$\oplus\ol{J^n}$) a m\^eme anneau total de fractions que $\oa[T]$. L'espace $\ox\times\cb$
\'etant normal et l'anneau $\oplus\ol{J^n}T^n$ \'etant la fermeture int\'egrale
de $\pcc(J)$ dans $\oa[T]$, il est int\'egralement clos. D'apr\`es [1] expos\'e 21,
cor.~3, $q:Z\to X'$ est le morphisme de normalisation de $X'$.

\th 3.3. Proposition|Soit $X$ un espace analytique complexe r\'eduit, $Y$ un
sous-espace analytique ferm\'e rare dont le support contient l'ensemble des
points non normaux de $X$. Soit $\ic$ l'id\'eal coh\'erent de $\oc_X$
d\'efinissant $Y$. L'\'eclatement normalis\'e de $Y$ est $X$-isomorphe au
morphisme canonique
$$\pro\opp_{n\in\nnf}\ol{\ic^n}\la X~.$$|

\dem Soit $\cc$ le conducteur de $\oc_X$ dans $\oc_{\ox}$. $Y$ contenant tous
les points non normaux de $X$, $\sqrt{\ic}$ est contenu dans $\sqrt{\cc}$.
Localement sur $X$, il existe donc un entier $k$ tel que $\ic^k\subset\cc$.
Soit, comme en 3.2, $\jc$ l'id\'eal coh\'erent de $\oc_{\ox}$ d\'efinissant
l'image r\'eciproque $W$ de $Y$ dans $\ox$. Localement sur $\ox$, d'apr\`es
2.8, il existe un entier $N$ tel que si $n\ge N$,
$\ol{\jc^n}=\jc^{n-N}\ol{\jc^N}$. Si donc $n\ge N+k$
$\ol{\jc^n}\subset\cc\cdot\jc^{n-N-k}\cdot\ol{\jc^N}\subset\oc_X\cap\ol{\jc^n}
=\ol{\ic^n}$ d'apr\`es 3.1.

Or (3.2), $\pro\oplus\ol{\jc^n}$ est isomorphe \`a l'\'eclatement normalis\'e de
$Y$ et on sait que pour tout entier $d$, $\pro\opp_{n\in\nnf}\ol{\jc^n}$
$(\resp\pro\opp_{n\in\nnf}\ol{\ic^n})$ est canoniquement isomorphe \`a
$\pro\opp_{n\in\nnf}\ol{\jc^{nd}}$ $(\resp\pro\opp_{n\in\nnf}\ol{\ic^{nd}})$. De
plus, le terme de degr\'e 0 dans les sommes directes est inessentiel.

On prendra garde que $\spe\opp_{n\in\nnf}\ol{\ic^n}$ peut ne pas \^etre un espace
normal comme le montre l'exemple suivant~: soit $X$ le cusp. Son anneau local
est $\cb\{t^2,t^3\}$. Soit $I$ l'id\'eal maximal. On v\'erifie que
$\ol{I^n}=\{\Sigma a_it^i\in\cb\{t\},a_i=0,i<2n\}$ si $n\ge 1$. Soit $\varphi
:\cb\{t^2,t^3\}[U,V]\to\opp_{n\in\nnf}\ol{I^n}$ le morphisme gradu\'e qui envoie
$U$ de degr\'e 1 sur $t^2\in\oi$ et $V$ de degr\'e 1 sur $t^3$. Le morphisme $\varphi $ est
surjectif. En effet, on remarque que $\ol{I^n}=I^n$. Il suffit donc que la
composante homog\`ene de degr\'e 1 de $\varphi $ soit surjective. Or $t^2$ est
l'image de $U$, $t^3$ celle de $V$, $t^{2n}$ est l'image de $t^{2(n-1)}U$,
$t^{2n+1}$ est l'image de $t^{2(n-1)}V$. D'autre part $\ker\varphi
=(t^3U-t^2V)$ et $\oplus\ol{I^n}$ n'est donc pas normal puisque $t={V\over U}$
est un \'el\'ement de son corps des fractions entier et ne lui appartenant pas.

\th 3.4. Proposition|Les notations et hypoth\`eses sont celles de 3.3. Pour que
l'\'eclatement $X'$ de $Y$ soit un espace analytique normal, il faut et il
suffit que localement sur $X$, il existe un entier $N$ tel que l'on ait $\ol{\ic^n}=\ic^n$ pour tout $n\ge N$.|

\dem Supposons $X'$ normal. D'apr\`es 2.8, pour tout $x\in X$, il existe un
voisinage $U$ de $x$ dans $X$ et un entier $N$ tel que~:

\ssection 3.4.1|$\ol{\ic^{n+N}}|U=\ic^n|U\cdot\ol{\ic^N}|U,\quad n\in\nnf$.

D'autre part, l'\'eclatement de $\ic^N$ \'etant canoniquement isomorphe \`a
celui de $\ic$ est \'egalement un espace analytique normal et d'apr\`es 2.1 iv)
et v) quitte \`a restreindre $U$, on d\'etermine un entier $k$ tel que~:

\ssection 3.4.2|$\ol{\ic^N}|U\cdot\ic^{Nk}|U=\ic^{N(k+1)}|U$.

En effet, on peut choisir pour $\kc$ et $\ic$, $\ic^N$, et pour $\jc$,
$\ol{\ic^N}$~; dans la d\'emonstration de v)$\to$i) on avait montr\'e que
localement il existe $k\in\nnf$ tel que $\ojc\kc^k\subset\ic\kc^k$. 
Rempla\c cons $n$ par $Nk$ dans 3.4.1~; nous obtenons~:
$$\ol{\ic^{N(k+1)}}|U=\ic^{Nk}|U\cdot\ol{\ic^N}|U~.$$
Comparant avec 3.4.2, il vient~:
$$\ol{\ic^{N(k+1)}}|U=\ic^{N(k+1)}|U~.$$
Soit $M=N(k+1)$, il suffit pour conclure de remarquer que si $s\ge M$
$$\ol{\ic^s}|U=\ic^{s-M}|U\cdot\ol{\ic^M}|U~.$$
R\'eciproquement, on sait que pour tout entier $d$, l'\'eclatement de $\ic$ est
isomorphe \`a $\pro\opp_n\ic^{nd}$ et d'apr\`es 3.3 que l'\'eclatement
normalis\'e de $\ic$ est isomorphe \`a $\pro\oplus\ol{\ic^n}$ donc aussi \`a
$\pro\oplus\ol{\ic^{nd}}$.

\rque 3.5. Remarque|Les notations et hypoth\`eses sont celles de 3.4. Si
l'\'eclate\-ment de $\ic$ est un espace analytique normal, pour tout $x$ de $X$,
il existe un entier $N_x$ tel que pour tout $f\in\oc_{X,x}$ tel que $\nu
_{\ic_x}(f)\ge N_x$, $\nu _{\ic_x}(f)$ est la partie enti\`ere de
$\bnu_{\ic_x}(f)$.

\rque 3.6. Avis|Nous recherchons un contre-exemple \`a la proposition 3.4 avec
$N=1$. ($X$ non r\'eduit s'abstenir).

\titre R\'ef\'erence|

\div 1|Cartan h|Familles d'espaces complexes et fondement de la g\'eom\'etrie
analytique, S\'eminaire Henri Cartan, {\bf 13, 2}|{\oldstyle~1960}{\rm --}{\oldstyle 1961}|
\par\medskip\noindent
\titre 4. $\bnu$ et \'eclatement normalis\'e|

\section 4.0|Dans ce paragraphe nous allons montrer comment l'on peut calculer
$\bnu$ apr\`es \'eclatement normalis\'e, et en d\'eduire d'importants
r\'esultats de finitude, notamment la rationalit\'e de $\bnu$, et le fait que
les alg\`ebres gradu\'ees associ\'ees \`a la ``filtration par le $\bnu$" sont
de pr\'esentation finie.

\section 4.1. Calcul de $\bnu$ dans l'\'eclatement normalis\'e|

\ssection 4.1.1|Soient $X$ un espace analytique complexe r\'eduit et $\ic$ un
$\oc_X$-id\'eal coh\'erent tel que le support de $\oc_X/\ic$ soit rare dans $X$.

Soient $\pi _0:X'_0\to X$ l'\'eclatement de $X$ d\'efini par $\ic$, et $\pi
:X'\to X$ l'\'eclatement normalis\'e de $X$ d\'efini par $\ic$, d\'efini comme
morphisme compos\'e $\pi :X'\build\to^n\fin X'_0\build\to^{\pi _0}\fin X$ o\`u
$n$ d\'esigne la normalisation de $X'_0$. On remarquera que, puisque $X$ est
r\'eduit et $\hbox{\rm supp}\oc_X/\ic$ rare dans $X$, l'espace $X'_0$ est r\'eduit et donc la
normalisation a un sens.

\ssection 4.1.2|Ainsi, $\ic\cdot\oc_{X'}$ est un id\'eal inversible dans
l'espace normal $X'$, et pour tout ouvert $U\subset X$ on peut d\'efinir un
ensemble de fonctions d'ordre sur l'anneau $\Gamma (U,\oc_X)$ comme suit~:

Consid\'erons les composantes irr\'eductibles $(D_\alpha )_{\alpha \in A(U)}$
de $D|U\big(=D\cap\pi ^{-1}(U)\big)$, o\`u $D$ est le diviseur exceptionnel de
$\pi $, \ie, le diviseur de $X'$ d\'efini par $\ic\cdot\oc_{X'}$.

Puisque $X'|U\big(=\pi ^{-1}(U)\big)$ est normal, pour chaque $\alpha \in A(U)$,
il existe un ouvert analytique dense $V_\alpha $ de $D_\alpha $ tel que, pour
tout point $x'\in V_\alpha $, $X'$ et $D_{\red}$ soient lisses en $x'$. On peut
alors choisir un syst\`eme de coordonn\'ees locales $(u,t_1\ld t_m)$ pour $X'$
en $x'$ tel que~:

\decale{i)} $\oc_{X',x'}\simeq\cb\{u,t_1\ld t_m\}$

\decale{ii)} $\ic\cdot{\cal O}_{X',x'}=(u^{e_\alpha })\cdot\cb\{u,t_1\ld t_m\}$, o\`u
$e_\alpha \in\nnf$.

\ssection 4.1.3|Consid\'erons maintenant $f\in\Gamma (U,\oc_X)$, et le
sous-espace $Z_f$ de $X'|U$ d\'efini par l'id\'eal $f\cdot\oc_{X'|U}(=f\circ\pi
|U)$. Puisque $D$ est un diviseur de $X'|U$, et que $X'|U$ \'etant normal,
d'apr\`es le ``hauptidealsatz", si $f\cdot\oc_{X'|U}$ ne s'annule pas
identiquement au voisinage d'un point $x'\in X'|U$, toutes les composantes
irr\'eductibles de $Z_f$ sont de codimension pure 1 en ce point, nous pouvons
trouver un ouvert analytique $U_\alpha \subset V_\alpha $ tel que, en tout
point $x'\in U_\alpha $ nous ayons~:

$ii)_f : f\cdot\oc_{X',x'}=(u^{m_\alpha })\cdot\cb\{u,t_1\ld t_m\}$

\noindent o\`u : 

$m_\alpha \in\nnf$, si $|D_\alpha |$ co\"{\i}ncide
ensemblistement avec une composante irr\'eductible de~$Z_f$,

$m_\alpha =0$, si $\dim(D_\alpha \cap Z_f)<\dim D_\alpha $,

\noindent et l'on pose~: $m_\alpha =+\infty$ si $f\circ \pi $ s'annule
identiquement au voisinage d'un point de $V_\alpha $, \cad en fait si $D_\alpha
$ est contenu dans une composante irr\'eductible de $X'|U$ qui est contenue
dans $Z_f$.

\ssection 4.1.4|Puisque $D_\alpha $ est irr\'eductible, $U_\alpha $ est
connexe, et puisque les entiers $m_\alpha $ et $e_\alpha $ que nous venons de
d\'efinir sont clairement localement constants sur $U_\alpha $, ils sont en
fait ind\'ependants du choix de $x'\in U_\alpha $.

\rque 4.1.5. Remarque|Pour tout $\alpha \in A(U)$, l'application $v_\alpha
:\Gamma (U,\oc_X)\to\nnf\cup\{+\infty\}$ d\'efinie par $v_\alpha (f)=m_\alpha $
satisfait~:

01) $v_\alpha (f+g)\ge\min\big(v_\alpha (f),v_\alpha (g)\big)$

02) $v_\alpha (f\cdot g)=v_\alpha (f)+v_\alpha (g)$

\noindent $v_\alpha $ est donc une fonction d'ordre.

On d\'efinit comme d'habitude $v_\alpha (I)$ par un id\'eal $I$ de $\Gamma
(U,\oc_X)$ par $v_\alpha (I)=\inff_{f\in I}v_\alpha (f)$, et l'on a alors~:
$$e_\alpha =v_\alpha \big(\Gamma (U,\ic)\big)~.$$

\th 4.1.6. Th\'eor\`eme|Soient $X$ un espace analytique complexe r\'eduit,
$\ic$ un $\oc_X$-id\'eal coh\'erent tel que $\hbox{\rm supp}\oc_X/\ic$ soit rare dans $X$,
et $K$ un sous-ensemble compact de $X$.

Il existe un voisinage ouvert $U$ de $K$ dans $X$ tel que~:

1) L'ensemble $A(U)$ des composantes irr\'eductibles du diviseur exceptionnel
de l'\'eclatement normalis\'e $\pi |U:X'|U\to X|U$ de $\ic$ est fini.

2) Pour tout $f\in\Gamma (K,\oc_X)$, il existe un voisinage ouvert $\wtu$ de
$K$ dans $X$ contenu dans $U$ et un prolongement $\tf$ de $f$ \`a $\wtu$ tel
que~:
$$\bnu^K_\ic(f)=\min_{\alpha \in A(\wtu)}{v_\alpha (\tf)\over v_\alpha
(\ic|\wtu)}~.$$
(Notations de 4.1.5 : $v_\alpha (\ic|\wtu)=v_\alpha \big(\Gamma
(\wtu,\ic)\big)$ o\`u l'on a pos\'e par d\'efinition :
$\bnu^K_\ic(f)=\inff_{x\in K}\bnu^x_\ic(f)$).|

Avant la d\'emonstration, donnons le

\th 4.1.7. Corollaire|Il existe un entier $q=q(K,y)$, que nous appellerons
``d\'enominateur universel" tel que pour tout ouvert $V$ contenant $K$, et tout
$f\in\Gamma (V,\oc_X)$, on ait
$$\bnu^K_\ic(f)\in{1\over q}\nnf\cup\{+\infty\}~.$$|

En particulier, $\bnu^K_\ic(f)$, s'il est fini, est un nombre rationnel. (On
peut prendre pour $q$ le p.p.c.m. des $v_\alpha (\ic|U)$).

La d\'emonstration du th\'eor\`eme 4.1.6 repose sur la proposition suivante.

\th 4.1.8. Proposition|Soient $X$ un espace analytique normal, d\'enom\-brable
\`a l'infini, $\ic$ un $\oc_X$-id\'eal inversible et $f\in\Gamma (X,\oc_X)$.
Notons $D$ le diviseur de Cartier d\'efini par $\ic$ et $D=\cupp_{\alpha \in
A}D_\alpha $ sa d\'ecomposition en composantes irr\'eductibles. Une condition
n\'ecessaire et suffisante pour que $f_x\in\ic_x$ pour tout $x\in X$ est que
pour tout $\alpha \in A$, il existe un point $x_\alpha \in D_\alpha $ tel que
$f_{x_\alpha }\in\ic_{x_\alpha }$.|

\dem La condition est \'evidemment n\'ecessaire, montrons qu'elle est
suffisante. Notons $\varphi \in\Gamma (X,\mc_X)$ la fonction m\'eromorphe
$\ic^{-1}\cdot f$. Nous appellerons sous-espace-polaire de $\varphi $, le
sous-espace analytique ferm\'e $P_\varphi $ de $X$ associ\'e \`a l'id\'eal
coh\'erent ${\cal P}_\varphi $ de $\oc_X$ d\'efini par
$$\Gamma (U,{\cal P}_\varphi )=\{h\in\Gamma (U,\oc_X):h\cdot(\varphi |U)\in\Gamma
(U,\oc_X)~.$$
Il est clair que $P_\varphi \subset D$ et que $f_x\in\ic_x$ si et seulement si
$x\notin P_\varphi $. Le germe en tout point $x\in X$ de $P_\varphi $ est soit
vide, soit de codimension 1. (Il contient une composante irr\'eductible de
codimension~1). En effet, \'ecrivons la d\'ecomposition primaire de ${\cal P}_{\varphi
,x}$ dans $\oc_{X,x}$
$${\cal P}_{\varphi ,x}={\bf q}_1\cap\cdots\cap {\bf q}_k~.$$
Supposons qu'aucun des $\sqrt{{\bf q}_i}$ ne soit de hauteur 1. Alors, pour tout
id\'eal premier ${\bf p}$ de hauteur 1, il existe $h\in {\cal P}_{\varphi ,x}$ tel que
$h\notin {\bf p}$. Sinon ${\cal P}_{\varphi ,x}\subset {\bf p}$ et puisque ${\bf p}$ est de hauteur 1,
il doit co\"{\i}ncider avec $\sqrt{{\bf q}_i}$ pour au moins un $i=1\cdots k$. Donc
$h\varphi _x\in\oc_{X,x}$ et $\varphi _x\in(\oc_{X,x})_{\bf p}$. L'anneau $\oc_{X,x}$ \'etant
normal, d'apr\`es le crit\`ere de Serre, ceci entra\^{\i}ne que $\varphi
_x\in\oc_{X,x}$, que $1\in {\cal P}_{\varphi ,x}$, et donc que $x\notin P_\phi$.

Puisque $P_\varphi \subset D$ et que tout $\alpha \in A$, $x_\alpha \notin
P_\varphi $, $P_\varphi =\emptyset$ (sinon il co\"{\i}nciderait avec une
composante de $D$) et $f_x\in\ic_x$ pour tout $x\in X$.

\rque 4.1.8.1. D\'emonstration du th\'eor\`eme (4.1.6)|Soit $\pi :X'\to X$
l'\'eclate\-ment normalis\'e de $\ic$ dans $X$, comme au 4.1.1. Le morphisme $\pi $ est propre, puisque compos\'e de deux morphismes propres, et donc tout
compact $K$ de $X$ poss\`ede un voisinage $U$ tel que~:

i) $D|U$ n'ait qu'un nombre fini de composantes irr\'eductibles (o\`u $D$ est
le diviseur exceptionnel de $\pi $, d\'efini par $\ic\cdot\oc_{X'}$, et $D|U$
signifie $D\cap\pi ^{-1}(U)$).

ii) Pour tout voisinage ouvert $\wtu\subset U$ de $K$, $\pi ^{-1}(\wtu)$
rencontre toutes les composantes irr\'eductibles de $D|U$. Si l'on pr\'ef\`ere,
toutes les composantes irr\'eductibles de $D|U$ rencontrent $\pi ^{-1}(K)$.

En effet, $\pi ^{-1}(K)$ est compact, et par la finitude locale du nombre des
composantes irr\'eductibles d'un espace analytique, poss\`ede un voisinage $V$
tel que $D\cap V$ n'ait qu'un nombre fini de composantes irr\'eductibles, qui
rencontrent toutes $\pi ^{-1}(K)$. $\pi $ \'etant propre, on peut choisir $V$
de la forme $\pi ^{-1}(U)$ avec $U\subset X$ voisinage de $K$.

Supposons maintenant $\inff_{x\in K}\bnu_{\ic_x}(f)\ge{p\over q}$, nombre
rationnel. Nous savons gr\^ace \`a (2.2 et 2.5) que $\bnu_{\ic_x}(f)\ge{p\over
q}$ \'equivaut \`a l'existence d'un voisinage ouvert $U_x$ de $x$ dans $X$ tel
que pour tout $x_1\in U_x$ on ait~:
$$f^q\cdot\oc_{X,x_1}\in\overline{\ic^p\cdot\oc_{X,x_1}}~.$$
Ainsi, si nous posons $\wtu=(\cupp_{x\in K}U_x)\cap U$, (qui est un voisinage
de $K$ dans $U$), nous avons, en tout point $x'\in\pi ^{-1}(\wtu)$
$$f^q\cdot\oc_{X',x'}\in\ic^p\cdot\oc_{X',x'}$$
et donc
$$q\cdot v_\alpha (f|\wtu)\ge p\cdot v_\alpha (\ic |\wtu)$$
et donc
$$\min_{\alpha \in A(U)}{v_\alpha (f|\wtu)\over v_\alpha (\ic |\wtu)}\ge {p\over
q}~.$$
Ceci nous montre que
$$\inff_{\alpha \in A(\wtu)}{v_\alpha (f|\wtu)\over v_\alpha (\ic
|\wtu)}\ge\inf_{x\in K}\bnu_{\ic_x}(f)~.$$
Or, supposons que cette in\'egalit\'e soit stricte~; nous pouvons choisir un
rationnel ${p'\over q'}$ et un point $x_0\in K$ tels que
$$\inf_{\alpha \in A(\wtu)}{v_\alpha (f|\wtu)\over v_\alpha (\ic
|\wtu)}>{p'\over q'}>\bnu_{\ic_{x_0}}(f)\ge\inf_{x\in K}\bnu_{\ic_x}(f)~.$$
Mais d'apr\`es la proposition 4.1.8, $\inff_{\alpha \in A(\wtu)}
{v_\alpha (f|\wtu)\over v_\alpha (\ic |\wtu)}\ge {p'\over q'}$ indique qu'en
tout point $x\in K$, $f^{q'}\oc_{X,x}\in \overline{\ic^{p'}\cdot\oc_{X,x}}$, et
donc, d'apr\`es (2.2), $\bnu_{\ic_x}(f)\ge {p'\over q'}$, d'o\`u la contradiction
cherch\'ee.

\section 4.2. D\'efinition des $\ol{\ic^{p/q}}$|

\th 4.2.1. D\'efinition|Soient $X$ un espace analytique complexe r\'eduit,
$\ic$ un $\oc_X$-id\'eal coh\'erent tel que $\hbox{\rm supp}\oc_X/\ic$ soit rare. Pour
tout ouvert $U\subset X$ on d\'efinit~:
$$\bnu^U_\ic : \Gamma (U,\oc_X)\la \rb$$
par
$$\bnu^U_\ic(f)=\inf_{x\in U}\big(\bnu^x_{\ic_x}(f)\big)~.$$|

\th 4.2.2. D\'efinition {\rm(hypoth\`eses et notations de 4.2.1)}|Soit $\nu
\in\rb_+$. Consid\'erons le faisceau~:
$$U\longmapsto \big\{f\in\Gamma (U,\oc_X)/\bnu^U_\ic(f)\ge\nu \big\}~.$$
Puisque clairement si $U'\subset U$, on a $\bnu^{U'}_\ic(f|U')\ge\bnu^U_\ic(f)$, ce
faisceau est un id\'eal de $\oc_X$, que nous noterons $\ol{\ic^\nu}$, et dont le
germe en $x\in X$ est $(\ol{\ic^\nu})_x=\{f\in\oc_{X,x}/\bnu_{\ic_x}(f)\ge\nu
\}$ et plus g\'en\'eralement, pour tout compact $K$, l'anneau des sections sur
$K$ est
$$\Gamma (K,\ol{\ic^\nu})=\{f\in\Gamma (K,\oc_X),\bnu^K_\ic(f)\ge\nu \}~.$$|

\th 4.2.3. Proposition {\rm (hypoth\`eses et notations de 4.2.1)}|Pour tout
nombre rationnel positif ${p\over q}$, l'id\'eal $\ol{\ic^{p/q}}$ de
$\oc_X$ est coh\'erent.|

\dem 1\`ere \'etape. V\'erifier le r\'esultat quand $X$ est normal et $\ic$
inversible.

Soit $X=\cupp_K V_k$ un recouvrement ouvert tel que pour tout $k$,
$$\ic|V_k=\varphi _k\cdot\oc_{X|V_k},\varphi _k\in\Gamma (V_k,\oc_X)~.$$

Consid\'erons, pour chaque $k$, le morphisme $\pi _k:\ol{V^q_k}\to V_k$ d\'efini
comme compos\'e $\ol{V^q_k}\build \to^n\fin V^q_k\build \la^{\omega _k}\fin V_k$
o\`u $n$ est la normalisation, et $\omega _k$ le morphisme structural
$\spe_{V_k}B_k\la V_k$, o\`u $B_k$ d\'esigne la $\oc_{V_k}$-alg\`ebre finie
$\oc_{V_k}[T]/(T^q-\varphi _k)$. Nous disposons de l'id\'eal inversible $J_k$
sur $\ol{V^q_k}$ engendr\'e par $(t\circ n)\oc_{\ol{V^q_k}}$, o\`u $t\in\Gamma
(V^q_k,\oc_{V^q_k})$ est l'image de $T$.

\th 4.2.3.1. Lemme|$\ol{\ic^{p/q}}|V_k{=}\pi_{k*}(J^p_k)\cap\oc_X|V_k$
(intersection dans $\pi_{k*}\oc_{\ol{V^q_k}}$).|

\dem Soient $U$ un ouvert de $V_k$, et $f\in\Gamma (U,\oc_X)$. Supposons
$\bnu^U_\ic(f)\ge{p\over q}$. Alors, puisque $\ic$ est suppos\'e inversible, et
que $X$ est normal, d'apr\`es 2.1, on a l'inclusion $f^q\cdot\oc_{X|U}\subset\ic^p\oc_{X|U}$ et a
fortiori $f^q\cdot\oc_{\ol{V^q_k}|U}\subset \ic^p\cdot\oc_{\ov^q_k|U}$.

Mais par d\'efinition de $J_k$,
$(\ic^p|U)\oc_{\ol{V^q_k}|U}=J^{pq}_k\cdot\oc_{\ol{V^q_k}|U}$. (La restriction
\`a $U$ dans $\ol{V^q_k}$ signifie bien s\^ur \`a $\pi ^{-1}_k(U)$). Et puisque
$\ol{V^q_k}$ est normal et $J_k$ inversible, on peut d\'eduire du fait que
$f^q\in J^{pq}_k$ que $f\in J^p_k$~; en effet $(J^{-p}_k\cdot
f)^q\subset\oc_{\ol{V^q_k}}$ ce qui, puisque $\ol{V^q_k}$ est normal, implique
$J^{-p}_k\cdot f\subset\oc_{V^q_k}$, donc $f\in J^p_k$.

Ainsi, nous venons de v\'erifier que~:
$$\bnu^U_\ic(f)\ge{p\over q}\Rightarrow f\in\Gamma (U,\pi_{k*}J^p_k\cap\oc_{V_k})$$
\cad encore~:
$$\ol{\ic^{p/q}}|U\subseteq (\pi_{k *}J^p_k\cap\oc_{V_k})|U~.$$
Montrons l'inclusion inverse. Soit $f\in\Gamma (U,\pi_{k*}J^p_k\cap\oc_{V_k})$ il vient, par d\'efinition de l'image directe~:
$$f\circ \pi _k=f\cdot\oc_{\ol{V^q_k}|U}\in J^p_k\cdot\oc_{\ol{V^q_k}|U}~,$$
d'o\`u
$$f^q\cdot\oc_{\ol{V^q_k}|U}\subset \varphi ^p_k\cdot\oc_{\ol{V^q_k}|U}$$
et donc, d'apr\`es 2.1 iv), puisque $\ol{V^q_k}$ est normal, que $\pi _k$ est
surjectif, pour tout $x\in U$, $f^q\cdot\oc_{X,x}\in\ol{\ic^p\cdot\oc_{X,x}}$,
donc $\bnu^x_\ic(f)\ge{p\over q}$ et $f\in\Gamma (U,\ol{\ic^{p/q}})$. QED pour le
lemme.

\th 4.2.3.2. Corollaire|Si $\ic$ est un id\'eal inversible d'un espace normal
$X$, tel que $\hbox{\rm supp}\oc_X/\ic$ soit rare dans $X$, $\ol{\ic^{p/q}}$ est un
$\oc_X$-id\'eal coh\'erent, pour tout rationnel positif ${p\over q}$.|

En effet, $\pi _k$ est propre, et $J^p_k$ un id\'eal coh\'erent pour tout $p$.
Donc d'apr\`es le th\'eor\`eme de Grauert, son image directe $\pi_{k*}J^p$
est coh\'erente, et son intersection dans le $\oc_{V_k}$-module coh\'erent $\pi_{k*}\oc_{\ol{V^q_k}}$ avec $\oc_{V_k}$ est encore un module coh\'erent, et donc
un id\'eal coh\'erent de $\oc_{V_k}$.

\rque 4.2.4. Remarque|$\ol{\ic^{p/q}}$ n'est pas en g\'en\'eral inversible, comme
le montre l'exemple suivant~: soit $X$ le c\^one quadratique d\'efini dans
$\cb^3$ par $\xi \eta =z^2$ et soit $\ic$ l'id\'eal engendr\'e par $\xi $.
L'id\'eal $\ol{\ic^{1/2}}$ est engendr\'e par $\xi $ et $z$.

2e \'etape de la d\'emonstration de 4.2.2~:

\th 4.2.5. Lemme|Soient $X$ un espace analytique complexe r\'eduit, $\ic$ un
$\oc_X$-id\'eal coh\'erent d\'efinissant un sous-espace rare dans $X$. Soit
$\pi :X'\to X$ l'\'eclatement normalis\'e de $\ic$ (4.1.1). On a~:
$$\ol{\ic^{p/q}}=\pi _*\ol{(\ic\cdot\oc_{X'})^{p/q}}\cap\oc_X~.$$|

\dem Soient $U$ un ouvert de $X$, et $f\in\Gamma (U,\ol{\ic^{p/q}})$. Pour tout
$x\in U$, $\bnu^x_\ic(f)\ge{p\over q}$, \ie,
$\bnu_{\ic_x}(f\cdot\oc_{X,x})\ge{p\over q}$, ce qui d'apr\`es (2.1)
entra\^{\i}ne que pour tout $x'\in\pi ^{-1}(U)$, on a $f^q\cdot\oc_{X',x'}\subset
\ic^p\cdot\oc_{X',x'}$, d'o\`u~:
$\bnu^{x'}_{\ic\cdot\oc_{X'}}(f\cdot\oc_{X'})\ge{p\over q}$  pour tout
$x'\in\pi ^{-1}(U)$. Donc
{\parindent=0pt
$$\leqalignno{
f\cdot\oc_{X'|U}&\subset \ol{(\ic\cdot\oc_{X'|U})^{p/q}},&\cr
\noalign{ce qui montre}
\ol{\ic^{p/q}}&\subset\pi _*\ol{(\ic\cdot\oc_{X'})^{p/q}}\cap\oc_X.&\cr}$$}
R\'eciproquement si $f\in\Gamma (U,\oc_X)$ est tel que $f\cdot\oc_{X'|U}\subset
\ol{(\ic\cdot\oc_{X'|U})^{p/q}}$, on a (2.1)
$$f^q\cdot\oc_{X'|U}\subset\ol{(\ic\cdot\oc_{X'|U})^p}=(\ic\cdot\oc_{X'|U})^p$$
puisque $\ic\cdot\oc_{X'}$ est un id\'eal inversible d'un espace normal. Donc
$$f^q\cdot\oc_{X'|U}\subset\ic^p\cdot\oc_{X'|U}$$
ce qui, toujours d'apr\`es (2.1), signifie que pour tout $x\in U$
$$f^q\cdot\oc_{X,x}\subset\ol{\ic^p}\cdot\oc_{X,x}~,$$
\ie, $\bnu^x_\ic(f)\ge{p\over q}$, $\forall x\in U$ et donc $f\in\Gamma
(U,\ol{\ic^{p/q}})$. QED pour le lemme.

\th 4.2.6. Corollaire|Soient $X$ un espace analytique r\'eduit et $\ic$ un
$\oc_X$-id\'eal coh\'erent d\'efinissant un sous-espace rare dans $X$.
Le $\oc_X$-id\'eal $\ol{\ic^{p/q}}$ est coh\'erent, et son germe
$(\ol{\ic^{p/q}})_x$ en un point $x\in X$ est
$\ol{\ic^{p/q}_x}=\{f\in\oc_{X,x}/\bnu_{\ic_x}(f)\ge p/q\}$.|

\dem D'apr\`es 4.2.3.2, $\ol{(\ic\cdot\oc_{X'})^{p/q}}$ est un $\oc_{X'}$-id\'eal
coh\'erent~; puisque $\pi :X'\to X$, \'eclatement normalis\'e de $\ic$, est
propre, $\pi _*\ol{(\ic\cdot\oc_{X'})^{p/q}}$ est un sous-$\oc_X$-module
coh\'erent du $\oc_X$-module coh\'erent $\pi _*\oc_{X'}$, donc $\ol{\ic^{p/q}}$,
intersection de deux sous-$\oc_X$-modules coh\'erents d'un $\oc_X$-module
coh\'erent, est coh\'erent.

\section 4.3. L'alg\`ebre gradu\'ee $\overline{\cal P}^{1/q}(\ic)$|

\th 4.3.0. D\'efinition|Soient $X$ un espace analytique complexe r\'eduit,
$\ic$ un $\oc_X$-id\'eal coh\'erent et $q$ un entier positif. On appelle alg\`ebre des
${1\over q}$-puissances de $\ic$ la $\oc_X$-alg\`ebre
$$\overline{\cal P}^{1/q}(\ic)=\oplus^\infty_{p=0}\ol{\ic^{p/q}}\cdot
T^{p/q}\subset\oc_X[T^{1/q}]$$
o\`u $\oc_X[T^{1/q}]$ d\'esigne la $\oc_X$-alg\`ebre $\oc_X[T,U]/(T-U^q)$.|

\th 4.3.1. Proposition|Pour tout entier positif $q$, la $\oc_X$-alg\`ebre gradu\'ee $\overline{\cal P}^{1/q}(\ic)$ est de pr\'esentation finie.| 

\dem D'apr\`es ([1] ch. I, 1.4) puisque $\ol{\ic^{p/q}}$ est un $\oc_X$-id\'eal
coh\'erent, il suffit de v\'erifier que $\overline{\cal P}^{1/q}(\ic)$ est une
$\oc_X$-alg\`ebre de type fini.

\th 4.3.2. Lemme| La $\oc_X$-alg\`ebre $\overline{\cal P}^{1/q}(\ic)$ est la fermeture int\'egrale dans
$\oc_X[T^{1/q}]$ de ${\cal P}(\ic)=\summ^\infty_{n=0}\ic^n
T^n\subset\oc_X[T]\subset\oc_X[T^{1/q}]$.|

\dem Nous savons gr\^ace \`a ([3] ch. VII \S2, ou Bourbaki, Alg. Comm., ch.
5-6, page 30) que la fermeture int\'egrale $\ol{{\cal P}(\ic)}$ de ${\cal P}(\ic)$ dans
$\oc_X[T^{1/q}]$ est une sous-$\oc_X$-alg\`ebre gradu\'ee de $\oc_X[T^{1/q}]$.
Nous pouvons donc nous restreindre au calcul de ses \'el\'ements homog\`enes.

Nous allons voir que le lemme 4.3.2 r\'esulte de~:

\th 4.3.3. Proposition|Soient $\oc$ une alg\`ebre analytique, $I$ un id\'eal de
$\oc$, $f\in\oc$ non inversible, et $d\in]0,+\infty]$. Les conditions suivantes
sont \'equivalentes~:

1) $\bnu_I(f)\ge d$

2) Il existe une relation de d\'ependance int\'egrale
$$f^k+a_1f^{k-1}+\cdots +a_k=0$$
avec $a_i\in\oc$ tels que $\nu _I(a_i)\ge id$.|

En fait, nous n'avons besoin ici que du cas particulier o\`u $d={p\over q}$,
que nous allons montrer directement, et qui d'ailleurs, apr\`es 4.1.7 suffit
pour montrer le cas g\'en\'eral. Une autre d\'emonstration, n'utilisant pas
le Th\'eor\`eme 2.1, est donn\'ee en appendice.

Supposons donc $\bnu_I(f)\ge{p\over q}$, \cad $\bnu_I(f^q)\ge p$. On a, gr\^ace
\`a (2.1)~: $f^q\in \ol{I^p}$, \cad que nous disposons d'une relation de
d\'ependance int\'egrale~:
$$f^{qk}+\summ^k_{i=1}a_i f^{q(k-i)}=0\quad\hbox{o\`u}\quad a_i\in I^{pi},$$
mais cette relation peut aussi se lire comme relation de d\'ependance
int\'egrale pour $f$
$$f^{qk}+\summ^{qk}_{j=q}a_j f^{qk-j}\quad\hbox{v\'erifiant}\quad \nu
_I(a_j)\ge j{p\over q}~.$$
R\'eciproquement, supposons que $f$ satisfasse
$$f^\ell+a_1 f^{\ell-1}+\cdots +a_\ell=0\quad\hbox{avec}\quad \nu _I(a_j)\ge
j{p\over q}~.\leqno (E)$$
Pour tout morphisme $\oc\build\to^\varphi \fin \cb\{t\}$, on a 
$$\varphi (f)^\ell+\varphi (a_1)\varphi (f)^{\ell-1}+\cdots +\varphi
(a_\ell)=0\leqno \big(\varphi (E)\big)$$
et
$$v\big(\varphi (a_j)\big)\ge j{p\over q}v\big(\varphi (I)\big)$$
o\`u $v$ d\'esigne la valuation $t$-adique. Nous allons en d\'eduire que
$$v\big(\varphi (f)\big)\ge{p\over q}v\big(\varphi (I)\big)~.$$
Pour cela remarquons d'abord que
$$v\big(\varphi (f)^{\ell-j}\cdot\varphi (a_j)\big)\ge (\ell-j)v\big(\varphi
(f)\big)+j{p\over q}v\big(\varphi (I)\big)$$
puisque $v\big(\varphi (a_j)\big)\ge j{p\over q}v\big(\varphi (I)\big)$. Donc
si nous avions
$$v\big(\varphi (f)\big)<{p\over q}v\big(\varphi (I)\big)~,$$
nous en d\'eduirions
$$v\big(\varphi (f)^{\ell-j}\cdot\varphi (a_j)\big)>\ell\cdot v\big(\varphi
(f)\big)\quad\hbox{pour tout}\quad 1\le j\le\ell$$
et d'apr\`es la relation de d\'ependance int\'egrale $\varphi (E)$~:
$$\ell\cdot v\big(\varphi (f)\big)\ge\min_{1\le j\le\ell}v\big(\varphi
(f)^{\ell-j}\cdot\varphi (a_j)\big)>\ell\cdot v\big(\varphi (f)\big)$$
et donc une contradiction.

Ceci montre que pour tout morphisme $\varphi :\oc\to\cb\{t\}$, nous avons
$v\big(\varphi (f)\big)\ge{p\over q}v\big(\varphi (I)\big)$, \ie, 
$v\big(\varphi (f^q)\big)\ge v\big(\varphi (I^p)\big)$ et d'apr\`es (2.1), ceci
implique $f^q\in\ol{I^p}$ dans $\oc$, \cad $\bnu_I(f)\ge{p\over q}$.

\ssection 4.3.5|D\'emontrons maintenant 4.3.2.

Il suffit de montrer que pour tout $x\in X$, le germe $\overline{\cal P}^{1/q}(\ic)_x$ en $x$
de $\overline{\cal P}^{1/q}(\ic)$ est la fermeture int\'egrale dans $\oc_{X,x}[T^{1/q}]$ de
${\cal P}(\ic)_x\subset\oc_{X,x}[T]$ et comme nous avons vu, il suffit de montrer
qu'un \'el\'ement homog\`ene de ${\cal O}_{X,x}[T^{1/q}]$ est entier sur ${\cal P}(\ic)_x$
si et seulement s'il appartient \`a $\overline{\cal P}^{1/q}(\ic)_x$.

Soit donc $f\cdot T^{p/q}\in\oc_{X,x}[T^{1/q}]$, entier sur ${\cal P}(\ic)_x$, \cad
satisfaisant une \'equation~:
$$(f\cdot T^{p/q})^k+A_1(f\cdot T^{p/q})^{k-1}+\cdots +A_k=0~;\quad A_i\in
{\cal P}(\ic)_x$$
dans $\oc_{X,x}[T^{1/q}]$. Par homog\'en\'eit\'e, nous pouvons supposer
$A_i=B_i\cdot T^{a(i)}$ avec $B_i\in\ic^{a(i)}_x$, et $a(i)={ip\over q}$, et
\'egaler \`a z\'ero le c\oe fficient de $T^{{kp\over q}}$ donne~:
$$f^k+B_1f^{k-1}+\cdots +B_k=0\quad\hbox{avec}\quad \nu _{\ic_x}(B_i)\ge
i{p\over q}$$
ce qui entra\^{\i}ne, d'apr\`es 4.3.3~:
$$f\in\ol{\ic^{p/q}_x},\quad\hbox{donc}\quad f\cdot T^{p/q}\in
\big(\overline{\cal P}^{1/q}(\ic)\big)_x=\overline{\cal P}^{1/q}(\ic_x)~.$$
R\'eciproquement, supposons $f\cdot T^{p/q}\in\overline{\cal P}^{1/q}(\ic_x)$, \cad
$f\in\ol{\ic^{p/q}_x}$, ou encore $f^q\in\ol{\ic^p_x}$. Ecrivons une relation
de d\'ependance int\'egrale~:
$$(f^q)^k+B_1(f^q)^{k-1}+\cdots + B_k=0\quad\hbox{avec}\quad B_i\in\ic^{p\cdot
i}_x~.$$
Apr\`es multiplication par $T^{kp}$, on peut r\'e\'ecrire ceci~:
$$(f\cdot T^{p/q})^{kq}+B_1T^p(f\cdot T^{p/q})^{(k-1)q}+\cdots + B_kT^{pk}=0$$
ce qui montre bien que $f\cdot T^{p/q}$ est entier sur ${\cal P}(\ic_x)$, et ach\`eve
la d\'emonstration de 4.3.2.

\rque Remarque|En fait 4.3.2 et 4.3.3 sont des \'enonc\'es \'equivalents. Le
m\^eme argument que celui d\'evelopp\'e en 2.8 montre que $\overline{\cal P}^{1/q}(\ic)$ est
un ${\cal P}(\ic)$-module de type fini donc une $\oc_X$-alg\`ebre de type fini, ce qui
ach\`eve la d\'emonstration de 4.3.1.

\th 4.3.6. Corollaire|Soient $X$ un espace analytique complexe r\'eduit, et
$\ic$ un $\oc_X$-id\'eal coh\'erent dont le support est rare dans $X$. Tout
point $x\in X$ poss\`ede un voisinage ouvert $U$ tel qu'il existe un entier
$N=N(U)$ tel que~:
$$(\ic|U)^k\cdot\ol{(\ic|U)^{p/q}}=\ol{(\ic|U)^{(p/q)+k}}\quad\hbox{d\`es
que}\quad {p\over q}\ge N~.$$|

\th 4.3.7. Corollaire|Dans la situation de 4.3, pour tout entier positif
$q$, la $\oc_X$-alg\`ebre gradu\'ee d\'efinie par
$$\ogr^{1/q}_\ic\oc_X=\oplus^\infty_{p=0}\ol{\ic^{p/q}}/\ol{\ic^{{p+1\over
q}}}$$ est une $\oc_X/\ic$-alg\`ebre gradu\'ee de pr\'esentation finie.|

\dem Tout d'abord, il r\'esulte imm\'ediatement du fait que
$\bnu_{\ic_x}(f\cdot g)\ge\bnu_{\ic_x}(f)+\bnu_{\ic_x}(g)$ pour tous $f,g\in {\cal O}_{X,x}$, que
$\ic\cdot\ol{\ic^{p/q}}\subset\ol{\ic^{{p+1\over q}}}$ et donc que
$\ogr^{1/q}_\ic\oc_X$ est en fait une $\oc_X/\ic$-alg\`ebre gradu\'ee. De plus,
ceci nous donne un homomorphisme surjectif de $\oc_X/\ic$-alg\`ebre gradu\'ees~:
$$\overline{\cal P}^{1/q}(\ic)\otimes_{\oc_X}\oc_X/\ic\la\ogr^{1/q}_\ic\oc_X\la 0$$
ce qui entra\^{\i}ne que $\ogr^{1/q}_\ic\oc_X$ est une $\oc_X/\ic$-alg\`ebre
gradu\'ee de type fini, d'apr\`es 4.3.1, et donc est de pr\'esentation finie
d'apr\`es ([1] ch.~I, 1.4) puisque chacune de ses composantes homog\`enes est
un $\oc_X/\ic$-module coh\'erent, comme quotient du $\oc_X/\ic$-module
coh\'erent $\ol{\ic^{p/q}}\ott_{\oc_X}\oc_X/\ic$.

\ssection 4.4.0|Nous allons maintenant utiliser l'existence du ``d\'enominateur
universel" de 4.1.7 pour montrer que localement, toutes les alg\`ebres
gradu\'ees que l'on a envie d'associer \`a la filtration par le $\bnu$ sont du
type \'etudi\'e ci-dessus.

\th 4.4.1. Proposition|Soient $X$ un espace analytique r\'eduit, et $\ic$ un
$\oc_X$-id\'eal coh\'erent. Consid\'erons, pour tout nombre r\'eel $\nu
\in\rb_+$ le faisceau $\ol{\ic^\nu} (\resp \ol{\ic^{\nu +}})$ associ\'e au
pr\'efaisceau $$\eqalign{
U &\longmapsto \{f\in\Gamma (U,\oc_X)/\bnu^U_\ic(f)\ge\nu \}\cr
\big(\resp U&\longmapsto \{f\in\Gamma (U,\oc_X)/\bnu^U_\ic(f)>\nu \}\big)\cr}$$
et la $\oc_X/\ic$-alg\`ebre gradu\'ee (par $\rb_0=\{\nu \in\rb,\nu \ge 0\}$)
$$\ogr_\ic\oc_X=\oplus_{\nu \in\rb_0}\ol{\ic^\nu} /\ol{\ic^{\nu +}}~.$$
Pour tout $x\in X$, il existe un voisinage ouvert $U$ de $x$ dans $X$ et un
entier $q$ tels que l'injection canonique
$\ogr^{1/q}_\ic\oc_X|U\hookrightarrow\ogr_\ic\oc_X|U$ soit un isomorphisme.|

\dem D'apr\`es 4.1.7, il existe un voisinage $U$ de $x\in X$ et un entier $q$
tel que $\forall x'\in U$, $\bnu_{\ic_{x'}}(f)\in{1\over q}\nnf$ pour tout
$f\in\oc_{X,x'}$. En effet, il suffit de choisir $U$ assez petit pour que
toutes les composantes irr\'eductibles du diviseur exceptionnel $D$ de
l'\'eclatement normalis\'e $\pi :X'\to X$ de $\ic$ qui rencontrent $\pi
^{-1}(U)$ rencontrent $\pi ^{-1}(x)$, et de prendre pour $q$ le p.p.c.m. de
$\{e_\alpha ,\alpha \in A(U)\}$.

\th 4.4.2. Corollaire|$\ogr_\ic\oc_X$ est une $\oc_X/\ic$-alg\`ebre gradu\'ee
de pr\'esen\-tation finie.|

En effet, \^etre de pr\'esentation finie est une condition locale, par
d\'efinition, et $\ogr^{1/q}_\ic\oc_X$ est de pr\'esentation finie d'apr\`es
4.3.5.

\rque 4.4.3. Remarque|$\ogr_\ic\oc_X$ \'etant en fait r\'eduite, elle est de
fa\c con naturelle une $\oc_X/\sqrt{\ic}$-alg\`ebre, et de pr\'esentation finie
en tant que telle, d'apr\`es 4.3.8.

\par\medskip\noindent
 
\titre Appendice au \S4|

D\'emonstration de 4.3.3 en g\'en\'eral, et une autre d\'emonstration de la
rationalit\'e de $\bnu_I(f)$ au moyen de la th\'eorie des installations.

\section A 1|Rappelons bri\`evement qu'on appelle installation la donn\'ee d'un
espace analytique complexe $Z$, de sous-espaces $X$ et $W$ de $Z$ et d'une
r\'etraction $r:Z\to W$, et que l'on note un tel objet $\trdot =(X,Z,W,r)$. La
cat\'egorie des installations est le cadre naturel de construction de polygones
de Newton. Nous nous pla\c cons dans le cas o\`u $r$ est une r\'etraction lisse \`a fibre isomorphe
\`a $\cb^t$. Etant donn\'e est un sous-espace analytique ferm\'e $Y$ de $X\cap
W$, pour tout $d\in\rb_+\cup\{\infty\}$, on associe \`a $\trdot $ son c\^one
normal anisotrope de long de $Y$, de tropisme $d$, not\'e $C^d_{\trdott
,Y}$. C'est le
sous-c\^one de $C_{W,Y}\build\times_Y\fin [C_{Z,W}\build\times_W\fin Y]$
d\'efini par l'id\'eal $\inn_Y(\trdot ,\ul{\trdot\kern-3pt}\kern3pt ,d)$ de
$\gr_Y\oc_W[Z_1\ld Z_t]$ construit comme suit: pour$$ g=\summ_{a\in\nnf^t}g_a z^a\ \in \oc_{Z,x}\simeq\oc\{z\},$$ on pose $\nu _Y(g,d)=\sup\{\mu /g\in I(d,\mu
)\}$, o\`u $I(d,\mu )$ est l'id\'eal de $\oc\{z\}$ engendr\'e par les $hz^a$
tels que $a+\nu _Y(h)/d\ge\mu $, et $$\inn_Y(g,d)=\summ_{a+\nu _Y(g_A)/d=\nu _Y(g,d)}\inn_Yg_aZ^a .$$
L'id\'eal $\inn_Y(\trdot ,\ul{\trdot\kern-3pt}\kern3pt ,d)$ est l'id\'eal de $\gr_Y\oc_W[Z_1\ld Z_t]$
 engendr\'e par les \'el\'ements $\inn_Y(g,d)$ lorsque
$g$ parcourt l'id\'eal $J$ d\'efinissant $X$ dans $Z$.

On a alors le

\th Th\'eor\`eme [1], [2]|\'Etant donn\'e $x\in Y\subset W\cap X$, il existe une
suite finie de nombres rationnels $0<d_1<\cdots<d_s<\infty$ telle que l'on ait~:

1) Pour tout $i$, $1\le i\le s+1$, le germe de $C^d_{\trdott,Y}$ en $x$ est
ind\'ependant de $d\in]d_{i-1},d_i[$. Notons $C_{\trdott,Y}(i)$ ce germe.

2) Les $C_{\trdott,Y}(i)$ sont tous distincts.

Les $d_i$ sont appel\'es tropismes critiques de l'installation $\trdot$ le
long de $Y$ en~$x$.

De plus, si nous consid\'erons l'installation~:
$$\big(C_{\trdott,Y}^d, C_{W,Y}\build\times_Y\fin [C_{Z,W}\build\times_W\fin
Y],C_{W,Y},P_1\big)$$
il existe $\varepsilon >0$ tel que le germe en $x$ du c\^one normal
anisotrope le long de $Y$, de tropisme 0, de cette installation s'identifie
canoniquement au germe en $x$ de $C^\delta _{\trdott,Y}$ pour $\delta
\in[d-\varepsilon ,d[$.|

\section A 2|Reprenons maintenant notre alg\`ebre analytique $\oc$,
correspondant \`a un germe d'espace analytique $(W,x)$, notre id\'eal $I$
d\'efinissant un sous-espace ferm\'e $(Y,x)\subset (W,x)$, et soient $Z=W\times \C$ et $y=x\times \{0\}$.
On consid\`ere $(W,x)$ comme un sous-espace de $(Z,y)$ en identifiant $(W,x)$ \`a $(W\times \{0\},y)$. 
Etant donn\'e $f\in\oc$ tel que $\bnu_I (f)\neq 0$, l'id\'eal $J=(z-f)\subset {\cal O}_{Z,y}\cong {\cal O}\{z\}$ d\'efinit un sous-espace ferm\'e $(X,y)$ de $(Z,y)$, et quitte \`a \'elever $f$ \`a une puissance assez grande, nous pouvons supposer que $f\in I$, d'o\`u $(Y,y)\subset (X\cap W,y)$. Nous pouvons donc consid\'erer l'installation 
$\trdot (f)=\trdot =(X,Z,W,\hbox{\rm pr}_1)$
 o\`u $\hbox{\rm pr}_1$ est la projection  $W\times \C\to W$, et le germe en $y$ des c\^ones normaux anisotropes $C^d_{{\trdott},Y}$ le long de $Y$.

\section A 3|Nous sommes maintenant en position pour d\'emontrer 4.3.3,\cad
l'\'equivalence, pour un $d\in]0,+\infty]$ des assertions~:

1) $\bnu_I(f)\ge d$.

2) Il existe une relation de d\'ependance int\'egrale
$$f^\ell+a_1f^{\ell-1}+\cdots + a_\ell=0\quad\hbox{avec}\quad\nu _I(a_i)\ge
d\cdot i\cdot$$

\ssection A.3.1|

1) $\Rightarrow$ 2). Le point crucial est que par d\'efinition de $\bnu$,
pour tout $\varepsilon >0$ (resp. pour tout $A$ si $d=+\infty$) il existe un
$k_0$ tel que si $k\ge k_0$, ${\nu _I(f^k)\over k}>d-\varepsilon $ $(\resp
{\nu _I(f^k)\over k}>A)$ ce qui signifie que dans notre installation
$\trdot=\trdot(f)$, puisque $z^k-f^k=(z-f)(z^{k-1}+z^{k-2}f+\cdots
+f^{k-1})\in J$ id\'eal d\'efinissant $X$ dans $Z$, par d\'efinition des
c\^ones normaux anisotropes le long de $Y$, (en \'ecrivant $z$ pour $\inn_Y\big((z),\trdot,d-\varepsilon )\big)$ nous avons
$$z^k\in\inn_Y(\trdot,\ul{\trdot\kern-3pt}\kern3pt ,d-\varepsilon
)\subset\gr_I\oc_{W,x}[z]$$
$\big(\inn_Y(\trdot,\ul{\trdot\kern-3pt}\kern3pt ,d-\varepsilon )\big)$ est
l'id\'eal d\'efinissant $C^{d-\varepsilon }_{\trdot,Y}$ dans
$C_{X,Y}\times\cb=C_{X,Y}\times [C_{Z,W}\build\times_W\fin Y]$. Ainsi, il doit
exister $g\in(z-f)\oc_{W,x}\{z\}=J$ tel que
$$\inn_N\big(\inn_Y(g,\trdot,d)\big)=z^k\quad (\resp
\hbox{~avec~}d=+\infty)$$
o\`u $N$ est l'id\'eal de $\gr_I{\cal O}[z]$ engendr\'e par
$\opp^\infty_{i=1}\gr^i_I{\cal O}$ (o\`u ${\cal O}_{W ,x}={\cal O})$. Nous pouvons
\'ecrire un tel \'el\'ement $g$ sous la forme~:
$$\eqalign{
g&=(b_0+b_1 z+\cdots + b_{k-1}z^{k-1}+\cdots)(z-f)\cr
&=-b_0 f+(b_0-b_1 f)z+\cdots + (b_{k-1}-b_k f)z^k+\cdots \cr}$$
et le fait que $\inn_N\big(\inn_Y(g,\ul{\trdot\kern-3pt}\kern3pt
,d)\big)=z^k$ impose que~:
$$\left\{
\eqalign{
b_{k-1}-b_k\cdot f & = 1 \mod I\cr
\nu _I(b_{i-1}-b_i\cdot f)  &\ge (k-i)d\qquad i=1\ld k-1\cr
\nu _I(b_0\cdot f)  &\ge k\cdot d\cr}\right.$$
La premi\`ere relation implique que $b_{k-1}$ est inversible dans ${\cal O}$, car $f\in I$. Nous pouvons
donc remplacer $g$ par $b^{-1}_{k-1}\cdot g$,
et donc supposer $b_{k-1}=1$.

D\'efinissons maintenant $a_i\in{\cal O},\ 1\leq i\leq k,$ par~:
$$\eqalign{
b_{k-2} &=f+a_1\cr
\vdots &\cr
b_{i-1} &=b_i\cdot f+a_{k-i}\cr
\vdots &\cr
b_0 &=b_1 f+a_{k-1}\cr
b_0 f&=-a_k\cr}$$
avec $\nu _I(a_1)\ge d$, $\nu _I(a_{k-i})\ge (k-i)d$, $0\le i\le k-1$.

Nous avons donc~:
$$b_0=f^{k-1}+a_1f^{k-2}+\cdots + a_{k-1}$$
et donc
$$f^k+a_1f^{k-1}+\cdots + a_k=0\quad\hbox{avec}\quad \nu _I(a_i)\ge id$$
ce qui ach\`eve de prouver 1) $\Rightarrow$ 2).

\ssection A.3.2|

Pour montrer que 2) $\Rightarrow$ 1), souvenons-nous qu'en 4.3.3, nous
l'avons montr\'e que tout $d$ rationnel. Si donc nous avons une relation
$$f^\ell+a_1 f^{\ell-1}+\cdots + a_\ell=0\quad\hbox{avec}\quad \nu
_I(a_i)\ge id~,$$
pour tout rationnel ${p\over q}<d$, nous avons $\bnu_I(f)\ge{p\over q}$, ce
qui montre bien $\bnu_I(f)\ge d$.

\ssection A.3.3|

\th Corollaire|$\bnu_I(f)=+\infty\Longleftrightarrow \exists k\  \hbox{\rm tel que}\  f^k=0$.|

{\bf A 4.}  ---   \th Th\'eor\`eme|$\bnu_I(f)$ est le plus grand tropisme critique $d_i$
de l'installation $\trdot(f)$ tel que le c\^one normal anisotrope
$C^{d_i}_{\trdott,Y}\subset C_{W,Y}\times\cb$ ne contienne pas $Y\times\cb$ au voisinage de $y$
($Y$ \'etant vu comme section nulle du c\^one relatif $C_{W,Y}\to Y$).|

Posons $d=\bnu_I(f)$.

Si $d$ n'\'etait pas un tropisme critique, il existerait $\varepsilon >0$ tel
que $\inn_Y(\trdot ,\ul{\trdot\kern-3pt}\kern3pt ,\delta)$, germe en $y$ de l'id\'eal d\'efinissant
$C^\delta _{\trdott,Y}$ dans $C_{W,Y}\times\cb$ ne d\'epende pas de $\delta
\in[d-\varepsilon ,d+\varepsilon ]$. Le m\^eme raisonnement qu'en A.3.1
permettrait de construire $g\in(z-f)\oc\{z\}$ tel que
$$z^k=\inn(g,d+\varepsilon )$$
et une relation de d\'ependance int\'egrale
$$f^k+a_1f^{k-1}+\cdots a_k=0$$
avec $\nu _I(a_i)\ge i(d+\varepsilon )$.

D'apr\`es A.3.2, on aurait $\bnu_I(f)\ge d+\varepsilon $. $\bnu_I(f)$ est
donc un tropisme critique.

Montrons maintenant que $Y\times\cb\not\hookrightarrow C^d_{\trdott,Y}$ et
que $Y\times\cb\hookrightarrow C^\delta _{\trdott,Y}$ si $\delta >d$, au voisinage de $y$. En
effet, on a une relation de d\'ependance int\'egrale~:
$$f^k+a_1f^{k-1}+\cdots + a_k=0\quad \nu _I(a_i)\ge id$$
soit
$$g=(z^k-f^k)+a_1(z^{k-1}-f^{k-1})+\cdots a_{k-1}(z-f)~.$$
On a~:
$$g\in(z-f)\oc\{z\}\quad\hbox{et}\quad\nu _Y(g,d)=k$$
$$\inn_Y(g,d)=z^k+\inn a_1z^{k-1}+\cdots +\inn a_k\notin \oplus_{i\ge
1}\gr^i_Y W[z]~.$$
Si par contre pour un $\delta >d$, on avait $Y\times\cb\not\hookrightarrow
C^\delta _{\trdott,Y}$, c'est \`a dire
$$\inn_Y(\trdot,\ul{\trdot\kern-3pt}\kern3pt ,\delta
)\not\hookrightarrow\oplus_{i\ge 1}\gr^i_Y W[z]$$
il existerait $g\in(z-f)\oc\{z\}$ tel que
$$\inn_Y(g,\delta )\notin\bigoplus_{i\ge i}\gr^i_Y W[z]$$
et on aurait avec $N=\opp_{i\ge 1}\gr^i_Y W[z]$
$$z^k=\inn_N\big(\inn(g,\delta )\big)$$
donc comme en A.3.1, $\bnu_I(f)\geq \delta>d$.

\ssection A.4.1|

\th Corollaire|$\bnu_I(f)\in\qb\cup\{\infty\}$.|

\titre R\'ef\'erences|

\livre 1|Lejeune M., Teissier B|Contribution \`a l'\'etude des
singularit\'es du point de vue du polyn\^ome de Newton|Th\`ese, Paris
VII|1973|

\livre 2|Lejeune M., Teissier B|Transversalit\'e, polygone de Newton et
installations|Ast\'erisque 7.8|1973|

\div 3|Zariski, Samuel|Commutative algebra, Van Nostrand| 1960|
\par\medskip\noindent
\titre 5. $\bnu$ et arcs analytiques|

\section 5.0|Les r\'esultats de ce paragraphe ont pour but de justifier le
calcul de $\bnu$, dans certains cas, par restriction \`a des arcs
``suffisamment g\'en\'eraux", et de pr\'eciser les limites de cette
m\'ethode.

\th 5.1. D\'efinition|Soit $X$ un espace analytique complexe. On appelle arc
analytique sur $X$ centr\'e en un point $x\in X$ un germe de morphisme
$h:(\db,0)\to(X,x)$ o\`u $\db=\{t\in\cb, |t|<1\}$. On notera $\ac_{X,x}$
l'ensemble des arcs non triviaux centr\'es en $x$, \cad des arcs tels que
$\im(h)\ne\{x\}$, ou encore tels que le morphisme
$h^*:\oc_{X,x}\to\oc_{\db,0}$ ne soit pas nul. [Un arc analytique centr\'e
dans un sous-espace $Y\subset X$, \ie, tel que $h(0)\in Y$, sera not\'e
$h:(\db,0)\to (X,Y)$]. Si $h\in\ac_{X,x}$, puisque $\oc_{\db,0}$ est un
anneau de valuation discr\`ete, on peut associer \`a $f\in\Gamma (X,\oc_X)$
l'entier $v(f\circ h)\in\zb_+$ o\`u $f\circ h=h^*(f\cdot\oc_{\db,0})$. De
m\^eme on peut d\'efinir $v(\ic\circ h)$ pour un id\'eal $\ic$ de $\oc_X$
par $v(\ic\circ h)=v\big(h^*(\ic\cdot\oc_{X,x})\big)$, $v$ d\'esignant
toujours la valuation naturelle de $\oc_{\db,0}$ (\ie, l'ordre en $t$, si
$\oc_{\db,0}\cong\cb\{t\}$).|

\th 5.2. Th\'eor\`eme|Soient $X$ un espace analytique complexe, $\ic$ un
$\oc_X$-id\'eal coh\'erent, et $x\in X$.

Pour tout $f\in\Gamma (X,\oc_X)$, on a~:
$$\bnu^x_\ic(f)=\inf_{h\in\ac_{X,x}}\big\{{v(f\circ h)\over v(\ic\circ
h)}\big\}$$|

\dem Apr\`es 4.1.7 ou A.4.1, nous pouvons \'ecrire~: $\bnu^x_\ic(f)={p\over
q}$ ($p,q$, entiers) et, apr\`es 0.2.9 et 2.1~:
$$f^q\cdot\oc_{X,x}\in\ol{\ic^p\cdot\oc_{X,x}}$$
et le crit\`ere valuatif de d\'ependance int\'egrale implique~:
{\parindent=0pt
$$\leqalignno{
v(f^q\circ h) &\ge v(\ic^p\circ h),\quad \forall h\in\ac_{X,x}\cr
\noalign{d'o\`u~:}
{v(f\circ h)\over v(\ic\circ h)} &\ge {p\over q}=\bnu^x_\ic(f),\quad\forall
h\in\ac_{X,x}\cr
\noalign{et~:}
\bnu^x_\ic(f) &\le \inf_{h\in\ac_{X,x}}\big\{{v(f\circ h)\over v(\ic\circ
h)}\big\}.\cr}$$}
Pour achever la preuve de 5.2, raisonnons par l'absurde et supposons
l'in\'egalit\'e ci-dessus stricte~: soit ${p'\over q'}$ un rationnel tel que
$$\bnu^x_\ic(f) <{p'\over q'}\le \inf_{h\in\ac_{X,x}}\big\{{v(f\circ h)\over
v(\ic\circ h)}\big\},$$
et donc tel que $v(f^{q'}\circ h)\ge v(\ic^{p'}\circ h)$, $\forall
h\in\ac_{X,x}$.

Le crit\`ere valuatif de d\'ependance int\'egrale fournit~:
$$f^{q'}\cdot\oc_{X,x}\in\ol{\ic^{p'}\cdot\oc_{X,x}}$$
et [2.1] nous donne alors~:
$$\bnu^x_\ic(f)\ge {p'\over q'}$$
et la contradiction cherch\'ee.

\section 5.3|Commentaires g\'eom\'etriques sur 5.2~: apr\`es la construction
faite au paragraphe pr\'ec\'edent, il n'est pas difficile d'imaginer un cas
o\`u nous saurons construire un arc donnant exactement le minimum. Reprenons
les notations de 4.1.2, et supposons que le minimum du th\'eor\`eme 4.1.6,
soit atteint sur une composante $D_\alpha $ (\it cf. \rm 4.1.3 et 4.1.4) du diviseur exceptionnel $D$ de
l'\'eclatement normalis\'e $\pi :X'\to X$ de $\ic$, telle que $\pi ^{-1}(x)$
rencontre l'ouvert $U_\alpha $ des points ``assez g\'en\'eraux" de $D_\alpha
$. Choisissons un germe de courbe analytique lisse centr\'e en un point
$x'\in\pi ^{-1}(x)\cap U_\alpha $, et transversal \`a $D_\alpha $ en $x'$,
\ie, un arc $h'\in\ac_{X',x'}$ tel que $h'^*(u)\in\oc_{\db,0}$ soit de
valuation 1, o\`u $(u)^{e_\alpha }\in\cb\{u,t_1\ld t_n\}$ est l'id\'eal de
$D_\alpha $ dans $X'$ au voisinage de $x'$. On aura alors clairement,
toujours avec les notations de 4.1~:
$$e_\alpha =v\big((\ic\cdot\oc_{X',x'})\circ h'\big) ;\ \  m_\alpha
=v\big((f\cdot\oc_{X',x'})\circ h'\big)$$
et donc, pour $h=\pi \circ h'\in\ac_{X,x}$,
$${v(f\circ h)\over v(\ic\circ h)}={m_\alpha \over e_\alpha
}=\bnu^x_\ic(f).\qquad (\hbox{d'apr\`es}\  4.1.6)$$

Soit maintenant, pour chaque composante irr\'eductible $|Y|_i$ de
$|Y|=\hbox{\rm supp}\oc_X/\ic$, $A(i)$ l'ensemble des $\alpha \in A$ ($A$ indiciant les
composantes irr\'eductibles du diviseur exceptionnel $D\subset X'$ \cf 4)
tels que $\pi (D_\alpha )=|Y|_i$ ($\pi (D_\alpha )$ est un sous-espace
ferm\'e de $X$ par la propret\'e de $\pi $). Et pour chaque $i$, soit
$\bnu_i=\min_{\alpha \in A(i)}({e_\alpha \over m_\alpha })$. Au moins au
voisinage d'un compact donn\'e de $X$, les $A(i)$ sont finis, et
$\bnu_i\in\qb\cup\{+\infty\}$. Pour $\alpha \in A(i)$, $\pi (U_\alpha )$
contient un ouvert analytique dense de $|Y|_i$ et donc en particulier pour
les $\alpha \in A(i)$ tels que ${e_\alpha \over m_\alpha }=\bnu_i$.

Ceci suffit pour montrer la

\th 5.4. Proposition|\'Etant donn\'e un id\'eal coh\'erent $\ic$ sur un espace
analytique complexe $X$, et $f\in\Gamma (X,\oc_X)$, l'ensemble des points
$x\in |Y|=\hbox{\rm supp}\oc_X/\ic$ tels qu'il existe $h_0\in\ac_{X,x}$ tel que~:
$$\bnu^x_\ic(f)={v(f\circ h_0)\over v(\ic\circ h_0)}=\inf_{h\in\ac_{X,x}}   
\big\{{v(f\circ h)\over v(\ic\circ h)}\big\}$$
contient un ouvert analytique partout dense de $|Y|$.|

En particulier, si $\hbox{\rm supp}\oc_X/\ic=\{x\}$, $|D|=|\pi ^{-1}(x)|$ et l'on peut
toujours calculer $\bnu^x_\ic(f)$ en prenant pour $h$ un disque $\pi \circ h'$,
o\`u $h'$ est construit comme en 5.3.

\th 5.5. Proposition|Soient $X$, $\ic$ et $f\in\Gamma (X,\oc_X)$ comme dans le
th\'eor\`eme 1.

Soit $p:X_1\to X$ un morphisme d'espaces analytiques complexes propre et
surjectif. Alors, posant $\ic_1=\ic\cdot\oc_{X_1}$, $f_1=f\circ p\in\Gamma
(X_1,\oc_{X_1})$, on a, pour tout $x\in X$
$$\bnu^x_\ic(f)=\bnu^{p^{-1}(x)}_{\ic_1}(f_1)\build
=^{\hbox{d\'ef}}\fin\min_{x_1\in p^{-1}(x)}\bnu^{x_1}_{\ic_1}(f_1)~.$$|

\dem Ceci r\'esulte imm\'ediatement de 5.2 et du crit\`ere valuatif de
propret\'e.

\rque 5.6. Remarque|La comparaison de 5.5 et du th\'eor\`eme 4.1.6 (\S4) peut
surprendre, puisque 5.5 implique que tous les arcs analytiques
$h':(\db,0)\to\big(X',\pi ^{-1}(x)\big)$ ($\pi :X'\to X$ est toujours
l'\'eclatement normalis\'e de $\ic$) nous donnent
$${v(f'\circ h')\over v(\ic\circ h')}\ge\min\big\{{v_\alpha (f)\over e_\alpha
}\big\}$$
ce qui est relativement peu ais\'e \`a v\'erifier directement.

\section 5.7|Un cas o\`u l'on peut calculer $\bnu$ \`a l'aide d'un arc
analytique.

Soient $X$ un espace analytique r\'eduit de dimension $d\ge 1$, $x\in X$ tel
que $\oc_{X,x}$ soit un anneau de Cohen-Macaulay, et $\ic$ un $\oc_X$-id\'eal
tel que $\hbox{\rm supp}\oc_X/\ic=\{x\}$, \ie, $I=\ic\cdot\oc_{X,x}$ est primaire pour
l'id\'eal maximal $\mc$ de $\oc_{X,x}=\oc$. Supposons d'abord que $I$ puisse
\^etre engendr\'e par une suite r\'eguli\`ere $(\varphi _1\ld\varphi _d)$
$(d=\dim\oc)$. Nous pouvons alors d\'efinir une famille de germes de courbes
dans $(X,x)$, param\'etr\'ee par $\pb^{d-1}$, comme suit~: $(\varphi _1\ld
\varphi _d)$ d\'efinissent un germe de morphisme $\Phi :(X,x)\to(\cb^d,0)$,
fini d'apr\`es le th\'eor\`eme de pr\'eparation de Weierstrass, puisque $\Phi
^{-1}(0)$ est d\'efini par $I$ et que $\sqrt{I}=\mc$. \`A chaque point
$\ell\in\pb^{d-1}$?? correspond une droite $\ell$ dans $(\cb^d,0)$ et $\big(\Phi
^{-1}(\ell),x\big)$ est un germe de courbe contenu dans $(X,x)$~; que nous
noterons $(C_\ell,x)$.

Il y a une bien meilleure fa\c con de d\'ecrire cette famille de courbes.
L'\'eclatement de $\ic$ dans $X$ peut \^etre d\'ecrit comme l'adh\'erence
$X'_0$ dans $X\times\pb^{d-1}$ du graphe de morphisme $X-\{x\}\to\pb^{d-1}$
d\'efini par $x'\mapsto \big(\varphi _1(x'):\cdots :\varphi
_d(x')\big)\in\pb^{d-1}$.
Le morphisme $\pi _0\colon X'_0\to X$ d\'eduit de la premi\`ere projection est l'\'eclatement de $\ic$ dans $X$. Nommons
$G_\ic:X'_0\to\pb^{d-1}$ le morphisme d\'eduit par la seconde projection.
Alors, on v\'erifie sans mal que $C_\ell=\pi _0\big(G^{-1}_\ic(\ell)\big)$.

On peut utiliser par exemple le fait que si l'on choisit des coordonn\'ees
homog\`enes $(T_1:\cdots :T_d)$ sur $\pb^{d-1}$, $X'_0$ est d\'efini dans
$X\times\pb^{d-1}$ par l'id\'eal engendr\'e par les $\{(T_i\varphi
_j-T_j\varphi _i),i\ne j\}$. De plus, $\pi ^{-1}_0(x)$ est isomorphe \`a
$\pb^{d-1}$ par $G_\ic |\pi_0 ^{-1}(x)$. Je noterai $\sigma $
l'isomorphisme inverse.

Ainsi nous pouvons consid\'erer $G_\ic:X'_0\build\leftrightarrows^\sigma
\fin\pb^{d-1}$ comme une famille de germes de courbes, et d'apr\`es le
th\'eor\`eme de Bertini-Sard, puisque $X'_0$ est r\'eduit, il existe un ouvert
de Zariski dense $U\subset\pb^{d-1}$ tel que si $\ell\in U$,
$\big(G^{-1}_\ic(\ell),\sigma (\ell)\big)$ est un germe de courbe r\'eduite.
Par ailleurs, quitte \`a restreindre $U$ on peut supposer que la normalisation
$X'\build\la^n\fin X'_0$ v\'erifie~:

\ssection 5.7.1|

1) Le morphisme compos\'e $G_\ic\circ n=\ol{G_\ic}:X'\to\pb^{d-1}$ est lisse
(=plat et \`a fibre lisse) en tout point de $\ol G^{-1}_\ic(\ell)$, $\ell\in
U$.

2) Le morphisme induit $\ol G^{-1}_\ic(\ell)\to G^{-1}_\ic(\ell)$ est la
normalisation, pour tout $\ell\in U$.

3) Le nombre des composantes irr\'eductibles de $\big(G^{-1}_\ic(\ell),\sigma
(\ell)\big)$ est ind\'ependant de $\ell\in U$.

4) Chaque composante irr\'eductible de $G^{-1}_\ic(\ell)$ est un arc sur $X'$
passant par un point non singulier de $X'$, o\`u $D_{\red}$ est aussi non
singuli\`ere, et transverse \`a $D_{\red}$ en ce point. [$D$ est comme
d'habitude le diviseur exceptionnel de l'\'eclatement normalis\'e $X'\to X$].

Puisque $\pi _0$ est un isomorphisme hors de $\pi ^{-1}(x)$, $\pi _0$ induit
un isomorphisme $G^{-1}_\ic(\ell)-\{\sigma (\ell)\}\simto C_{\ell}-\{x\}$.

Soit $C_\ell=\cupp^r_1 \Gamma _q$ la d\'ecomposition de $C_\ell$ en
composantes irr\'eductibles. Si $\ell\in U$, chaque $\Gamma _q$ est r\'eduite
et irr\'eductible, et fournit un arc $h_q:(\db,0)\to(X,x)$.

\th Proposition|Pour tout $f\in\oc_{X,x}$, il existe un ouvert de Zariski non
vide $V$ de $\pb^{d-1}$ tel que si $\ell\in V$, il existe une composante
irr\'eductible $\Gamma _q$ de $C_\ell$ telle que|

\ssection 5.7.2|$\bnu^x_\ic(f)={v(f\circ h_q)\over v(\ic\circ h_q)}~.$

\dem Il suffit d'appliquer 5.3 \`a la situation cr\'e\'ee ci-dessus.

Ainsi, nous pouvons affirmer dans ce cas-ci qu'une composante irr\'eductible
d'une courbe d\'efinie par $d-1$ combinaisons lin\'eaires "g\'en\'eriques" de
g\'en\'erateurs de $\ic_x$, calcule $\bnu$ pour nous.
\par\medskip\noindent
\titre 6. $\bnu$ et exposants de \L ojasiewicz|

\section 6.0|Nous allons montrer ici qu'un calcul de $\bnu$ est en fait un
calcul d'exposant de \L ojasiewicz, et en d\'eduire la rationalit\'e de ces
derniers en g\'eom\'etrie analytique complexe. Ce paragraphe-ci est clairement
le seul que l'on ne puisse pas transcrire en g\'eom\'etrie alg\'ebrique~!

\defi 6.1. D\'efinition|Soient $X$ un espace analytique complexe r\'eduit,
$\ic$ un $\oc_X$-id\'eal coh\'erent, $f\in\Gamma (X,\oc_X)$ et $K$ un
sous-ensemble compact de $X$. L'exposant de \L ojasiewicz $\theta _K(f,\ic)$ de
$f$ par rapport \`a $\ic$ sur $K$ est la borne inf\'erieure de l'ensemble des
$\theta \in\rb_+$ tels qu'il existe un voisinage ouvert $U$ de $K$ dans $X$ et
une constante $C\in\rb_+$ tels que
$$|f(x)|^\theta \le C\cdot\sup_{g\in\Gamma (U,\ic)}|g(x)|\quad\ \hbox{\rm pour tout}\ x\in U.$$
Si l'ensemble de ces $\theta $ est vide, on convient de poser $\theta
_K(f,\ic)=+\infty$. Si par ailleurs $\Gamma (U,\ic)$ est de type fini quand $U$
est un voisinage assez petit de $K$ (ce sera le cas si $K=\{x\}$), $\theta
_K(f,\ic)$ est aussi la borne inf\'erieure dans $\rb_+\cup\{\infty\}$ de
l'ensemble des $\theta $ tels qu'il existe $U$ et $C$ tels que
$$|f(x)|^\theta \le C\cdot\sup^m_{i=1} |g_i(x)|\ \hbox{\rm pour tout}\ x\in U,$$o\`u
$(g_1\ld g_m)$
engendrent $\Gamma (U,\ic)$.|

\rque 6.2 Remarque|On peut aussi d\'efinir l'exposant de \L ojasiewicz
$\theta _K(\ic',\ic)$ d'un id\'eal $\ic'$ par rapport \`a $\ic$ sur $K$~:
$$\theta _K(\ic',\ic)=\sup_{f\in\Gamma (X,\ic')}\theta _K(f,\ic)~.$$

\th 6.3. Th\'eor\`eme|

1) $\theta _K(f,\ic)={1\over \bnu^K_\ic(f)}$, o\`u $\bnu^K_\ic(f)=\inff_{x\in
K}\bnu^x_\ic(f)$.

\noindent [On convient bien s\^ur que $\bnu^K_\ic(f)=0\Rightarrow\theta
_K(f,\ic)=+\infty$].

2) Et de plus, il existe un voisinage $U$ de $K$ dans $X$ et une constante
$C\in\rb_+$ tels que
$$|f(x)|^{\theta _K(f,\ic)}\le C\cdot\sup_{g\in\Gamma (U,\ic)}
|g(x)|\quad\hbox{pour tout}\quad x\in U,$$
c'est \`a dire que la borne inf\'erieure de 6.1 est atteinte.|

\th 6.4. Corollaire|Pour tout id\'eal coh\'erent $\ic$ sur un espace analytique
complexe r\'eduit $X$, tel que $\hbox{\rm supp}\oc_X/\ic$ soit rare, pour tout
$f\in\Gamma (X,\oc_X)$ et tout compact $K\subset X$,
$$\theta _K(f,\ic)\in\qb_0\cup\{+\infty\}~.$$|

\rque D\'emonstration de 6.3|Posons $\bnu^K_\ic(f)={p\over q}$ (4.1.6). Apr\`es
(4.2.3), il existe un voisinage $U_0$ de $K$ dans $X$ tel que pour tout $x\in
U_0$, $f^q\cdot\oc_{X,x}\in\ol{\ic^p\cdot\oc_{X,x}}$, et le th\'eor\`eme de
majoration (2.1.vi) nous fournit une constante $C$ telle que $|f(x)|^q\le
C\cdot\supp_{g\in\Gamma (U_0,\ic^p)}|g(x)|$ pour tout $x\in U_0$. Mais
d'apr\`es (2.1.iv), l'id\'eal (contenu dans $\ic^p$) engendr\'e par les
puissances $p$-i\`emes d'\'el\'ements de $\ic$ a m\^eme cl\^oture int\'egrale
que $\ic^p$, et donc en appliquant \`a nouveau le th\'eor\`eme de majoration
nous pouvons \'ecrire au prix d'un changement de la constante $C$
{\parindent=0pt
$$\leqalignno{
|f(x)|^q &\le C\cdot\sup_{g\in\Gamma (U_0,\ic)}|g(x)|^p,\hbox{~i.e.}\cr
|f(x)|^{q/p} &\le C^{1/p}\cdot\sup_{g\in\Gamma (U_0,\ic)}|g(x)|.\cr
\noalign{D'o\`u~:}
\theta _K(f,\ic) &\le{q\over p}={1\over \bnu^K_\ic(f)}.\cr}$$}
Mais, si nous supposons l'in\'egalit\'e stricte, il existe ${q'\over
p'}<{q\over p}$, un voisinage $U'$ de $K$ dans $X$ et une constante
$D\in\rb_+$ tels que~:
$$|f(x)|^{q'}\le D\cdot\sup_{g\in\Gamma (U',\ic)}|g(x)|^{p'}\quad\hbox{\rm pour tout}\  x\in U',$$
et le th\'eor\`eme de majoration, avec l'argument pr\'ec\'edent, nous donne~:
$$f^{q'}\cdot\oc_{X,x}\in\ol{\ic^{p'}\cdot\oc_{X,x}}\quad\hbox{\rm pour tout}\  x\in U',$$
donc~:
$$\bnu^K_\ic(f)\ge{p'\over q'}~,$$
(apr\`es (2.1)) et la contradiction cherch\'ee. Ceci d\'emontre 1) et 2) de
6.3, avec $U=U_0$.

\section 6.5|Comme les morphismes propres conservent les in\'egalit\'es du
type consid\'er\'e ici, au prix \'eventuellement d'une modification des
constantes, on peut tr\`es bien utiliser 6.3 pour d\'emontrer 5.5 (\S5) sans
utiliser le crit\`ere valuatif de propret\'e.

\titre 7. Th\'eor\`eme r\'ecapitulatif|

\section 7.1. Notations|Soient $X$ un espace analytique complexe r\'eduit,
$\ic$ un $\oc_X$-id\'eal coh\'erent tel que $|Y|=\sup\oc_X/\ic$ soit rare dans
$X$, et $K$ un sous-ensemble compact de $X$. Soient $\pi :X'\to X$
l'\'eclatement normalis\'e de $\ic$, $D$ le diviseur exceptionnel, sous-espace
de $X'$ d\'efini par l'id\'eal inversible $\ic\cdot\oc_{X'}$. Soit $A(K)$
l'ensemble fini tel que les composantes irr\'eductibles $D_\alpha $ de $D$,
avec $\alpha \in A(K)$ soient exactement celles qui rencontrent $\pi ^{-1}(U)$
pour tout voisinage ouvert $U$ de $K$ dans $X$. Soit enfin $e_\alpha $ la
multiplicit\'e de $\ic$ en un point $x'\in V_\alpha $, o\`u $V_\alpha \subset
D_\alpha $ est un ouvert analytique dense dans $D_\alpha $ en chaque point
duquel $X'$ et $D_{\alpha ,\red}$ sont non singuliers, et
$\ic\cdot\oc_{X',x'}=u^{e_\alpha }\cdot\oc_{X',x'}$, $D_{\alpha ,\red}$
\'etant d\'efini par $(u)\cdot\oc_{X',x'}$.

\th 7.2. Th\'eor\`eme|\'Etant donn\'es un nombre rationnel ${p\over q}>0$ et une fonction $f\in\Gamma (X,\oc_X)$, les
conditions suivantes sont \'equivalentes~:

1) $f^q\cdot\oc_{X,x}\in\ol{\ic^p\cdot\oc_{X,x}}\quad\forall x\in K$.

2) $\bnu^K_\ic(f)=\inff_{x\in K}\bnu^x_\ic(f)\ge{p\over q}$.

3) Pour tout $x\in K$, il existe des $a_i\in\oc_{X,x}$, tels que $\nu _{\ic_
x}(a_i)\ge{p\over q}\cdot i$ (o\`u $\ic_x=\ic\cdot\oc_{X,x})$ et que
$f^k_x+a_1f^{k-1}_x+\cdots +a_k=0$ dans $\oc_{X,x}$ (o\`u
$f_x=f\cdot\oc_{X,x}$).

3') Si $K$ est un polycylindre, il existe $a_i\in\Gamma (K,\oc_X)$, $i=1\cdots
k$, tels que $\nu _{\Gamma (K,\ic)}(a_i)\ge{ip\over q}$ et que $f^k+a_1
f^{k-1}+\cdots +a_k=0$ dans $\Gamma (K,\oc_X)$.

4) Pour tout arc $h:(\db,0)\to(X,K)$ (\ie, $h(0)\in K$) on a~: ${v(f\circ
h)\over v(\ic\circ h)}\ge {p\over q}$, o\`u $v$ d\'esigne la valuation naturelle de ${\cal O}_{\db, 0}\simeq\C\{t\}$, c'est \`a dire l'ordre en $t$.

5) Pour tout morphisme $\pi :X'\to X$ propre, dont l'image contient $K$, et
tel que $\ic\cdot\oc_{X'}$ soit inversible, et que $X'$ soit un espace
analytique normal, il existe un voisinage ouvert $U'$ de $\pi ^{-1}(K)$ dans
$X'$ tel que~: $f^q\cdot\oc_{U'}\in\ic^p\cdot\oc_{U'}$.

6) Il existe un voisinage ouvert $U$ de $K$ dans $X$, et une constante
$C\in\rb_+$ tels que
$$|f(x)|^{q/p}\le C\cdot\sup_{g\in\Gamma (U,\ic)}|g(x)|\quad\hbox{pour
tout}\quad x\in U~.$$| 

De plus,

A) Il existe $\alpha _0\in A(U)$ tel que $\bnu^K_\ic(f)={\mc_{\alpha _0}\over
e_{\alpha _0}}$ o\`u $\mc_{\alpha _0}$ est la multiplicit\'e de $f$ le long de
$D_{\alpha ,\red}$ en tout point $x'$ d'un ouvert analytique dense $V_{\alpha
_0}\subset U_{\alpha _0}\subset D_{\alpha _0}$, et donc, pour tout arc
$h:(\db,0)\to(X,U)$ de la forme $\pi \circ h'$, o\`u $h':(\db,0)\to\big(X',\pi
^{-1}(U)\big)$ est tel que $h'(0)\in V_{\alpha _0}$ et que $h'^*(U)$ soit de
valuation 1 dans $\oc_{\db,0}$ (6.1.1), on a~:
$${v(f\circ h)\over v(\ic\circ h)}={\mc_{\alpha _0}\over e_{\alpha
_0}}=\bnu^K_\ic(f)=\inf_{h\in\ac_{X,K}}{v(f\circ h)\over v(\ic\circ h)}$$
et tout voisinage $U$ de $K$ dans $X$ contient de tels arcs, \cad que l'on peut
trouver de tels arcs avec $h(0)\in U$.

B) Le faisceau d'id\'eaux de $\oc_X$ d\'efini par
$$\ol{\ic^{p/q}}(U)=\{f\in\Gamma (U,\oc_X)/\bnu^U_\ic(f)\ge{p\over q}\}$$
(o\`u $\bnu^U_\ic=\inff_{x\in U}\bnu^x_\ic)$ est coh\'erent pour tout rationnel
${p\over q}$, et la $\oc_X$-alg\`ebre gradu\'ee
$\opp^\infty_{p'=0}\ol{\ic^{p/q}}T^{p/q}$ est de pr\'esentation finie pour tout
entier $q$. Enfin, la $\oc_X$-alg\`ebre gradu\'ee $\ogr_\ic\oc_X=\opp_{\nu
\in\rb_+}\ol{\ic^\nu }/\ol{\ic^{\nu +}}$ co\"{\i}ncide localement sur $X$ avec
une alg\`ebre du type
$$\oplus^{+\infty}_{p=0}\ol{\ic^{p/q}}/\ol{\ic^{{p+1\over q}}}$$
et est donc de pr\'esentation finie, puisque cette derni\`ere l'est comme
quotient de
$$\big(\oplus^{+\infty}_{p=0}\ol{\ic^{p/q}}T^{p/q}\big)
\otimes_{\oc_X}\oc_X/\ol{\ic^{1/q}}~.$$
 
\titre Appendice par J.J.~Risler|

\titre Les exposants de \L ojasiewicz dans le cas analytique r\'eel|

Dans le cas r\'eel, l'exposant de \L ojasiewicz n'a pas d'interpr\'etation
alg\'ebrique simple analogue au $\bnu$ ou \`a la notion de cl\^oture
int\'egrale~; je vais cependant montrer que comme dans le cas complexe on peut
le calculer \`a l'aide d'arcs analytiques, ou de morphismes analytiques r\'eels
qui jouent un r\^ole analogue \`a celui de l'\'eclatement normalis\'e~; il en
r\'esultera que dans le cas r\'eel aussi les exposants de \L ojasiewicz sont
toujours rationnels.

Les r\'ef\'erences au s\'eminaire seront pr\'ec\'ed\'ees de la lettre $S$.

\titre 1. Pr\'eliminaires|

\th 1.1. D\'efinition {\rm~(\cf [R])}|Soit $A$ une $\rb$-alg\`ebre analytique~;
on dit qu'un id\'eal $I\subset A$ est r\'eel s'il satisfait \`a la condition
suivante~:
$$f_i\in A (1\le i\le p)\hbox{~et~}f^2_1+\cdots  + f^2_p\in I\Rightarrow f_i\in I
(1\le i\le p)~.$$|

On a alors la proposition suivante (\cf [R])~:

\th 1.2. Proposition|Soient $A$ une $\rb$-alg\`ebre analytique, $I$ un id\'eal
premier de $A$ tel que $\dim(A/I)=h$, $(X,x)$ un repr\'esentant du germe
analytique d\'efini par $A/I$~; les conditions suivantes sont \'equivalentes~:

a) $I$ est r\'eel.

b) $I$ est l'id\'eal de tous les \'el\'ements de $A$ nuls sur le germe de $X$
au point $x$.

c) $X$ poss\`ede un point lisse de dimension $h$ dans tout voisinage du point
$x$.|

\section 1.3|Soit $(X;\oc_{X,x})$ un germe analytique dans $\rb^n$~; on a alors
$\oc_{X,x}\simto\oc_{n}/I$ (avec $\oc_n=\rb\{x_1\ld x_n\}$). Notons $I(X)$
l'id\'eal de $\oc_n$ form\'e des s\'eries nulles sur $X$ ($I(X)$ est la racine
r\'eelle de $I$ ([R]))~; on dit que $X$ est normal en $x$ si~:

a) $I=I(X)$

b) l'anneau $\oc_{X,x}$ est int\'egralement clos.

Dans ce cas l'anneau $\oc_{\wtx,x}=\oc_{X,x}\ott_{\rb}\cb$ est aussi
int\'egralement clos, autrement dit $X$ poss\`ede un complexifi\'e $\wtx$ qui
est aussi un espace normal~: l'anneau $\oc_{X,x}\ott_{\rb}\cb$ est en effet
int\`egre ([R], proposition 6.1), et il r\'esulte d'un th\'eor\`eme
d'alg\`ebre classique qu'il est alors int\'egralement clos (\cf par exemple
Bourbaki, Alg. Comm., chap.~V). On dit qu'un espace analytique r\'eel
$(X,\oc_X)$ est normal s'il est normal en chaque point (rappelons que $\oc_X$
d\'esigne un faisceau coh\'erent de $\rb$-alg\`ebres analytiques~; \cf [H]
pour la notion d'espace analytique r\'eel).

\section 1.4|Si $\ib=]-1,+1[\subset\rb$, nous noterons comme dans S.5.1, $v$ la
valuation naturelle de l'anneau $\oc_{\ib,0}\simto\rb\{t\}$.

\titre 2. Arcs analytiques r\'eels et r\'esolution des singularit\'es|

On a dans le cas r\'eel le th\'eor\`eme suivant, analogue \`a une partie du
th\'eor\`eme S.7.2~:

\th 2.1. Th\'eor\`eme|Soient $(X,\oc_X)$ un espace analytique r\'eel, $K$ un
compact de $X$, $\ic$ un $\oc_X$-id\'eal coh\'erent, $f\in\Gamma (X,\oc_X)$ et
${p\over q}$ un nombre rationnel~; les conditions suivantes sont
\'equivalentes~:

1) Il existe un voisinage ouvert $U$ de $K$ dans $X$ et une constante $C>0$
tels que~:
$$|f(x)|^{q/p}\le C\sup_{g\in\Gamma (U,\ic)}|g(x)|\quad\hbox{\rm pour tout}\  x\in U~.$$

2) Pour tout arc analytique r\'eel $h:(]-1,1[,0)\to (X,K)$ (\ie, tel que
$h(0)\in K$), on a
$${v(f\circ h)\over v(I\circ h)}\ge {p\over q}~.$$

3) Pour tout morphisme analytique r\'eel $\pi =X'\to X$ dont l'image contient
$K$ et tel que~:

a) $\pi $ soit propre et $X'$ soit normal (\cf 1.3)

b) $\ic\oc_{X'}$ soit localement principal

c) $\forall x'\in \pi ^{-1}(K)$, l'id\'eal $\sqrt{\ic\oc_{X',x'}}$ soit un
id\'eal r\'eel (\cf 1.1~; $\sqrt{\ic\oc_{X',x'}}$ d\'esigne la racine de
l'id\'eal $\ic\oc_{X',x'}$),

\noindent il existe un voisinage ouvert $U'$ de $\pi ^{-1}(K)$ dans $X'$ tel
que~ $f^q\oc_{U'}\in\ic^p\oc_{U'}$.

4) Il existe un morphisme analytique r\'eel $\pi :X'\to X$ dont l'image
contient $K$ et v\'erifiant les propri\'et\'es a), b), c) ci-dessus et un
voisinage ouvert $U'$ de $\pi ^{-1}(K)$ dans $X'$ tels que
$f^q\oc_{U'}\in\ic^p\oc_{U'}$.|

\dem

1) $\Rightarrow$ 2) : Soit $h:(]-1,1[,0)\to (X,K)$ un arc analytique r\'eel~;
on a par hypoth\`ese $|f(x)|^q\le C\supp_{g\in\Gamma (U,\ic)}|g(x)|^p$
$\hbox{\rm pour tout}\  x\in U$, d'o\`u $|f\circ h(t)|^q\le C\supp_{g\in\Gamma
(U,\ic)}|g\circ h(t)|^p$ pour $t$ voisin de 0 dans $]-1,1[$, d'o\`u
imm\'ediatement $v\big((f\circ h)^q\big)\ge v\big((\ic\circ h)^p\big)$~; soit
$qv(f\circ h)\ge pv(\ic\circ h)$. 

2) $\Rightarrow$ 3) : Soit $\pi =X'\to X$ un morphisme analytique r\'eel
propre v\'erifiant les conditions \'enonc\'ees dans 3)~; si l'on suppose qu'il
existe $x'\in\pi ^{-1}(K)$ tel que $f^q\oc_{X',x'}\notin\ic^p\oc_{X',x'}$, il
faut montrer qu'il existe un arc analytique r\'eel $h:(]-1,1[,0)\to (X,K)$ tel
que $v(f\circ h)/v(\ic\circ h)<p/q$, ce qui va r\'esulter du lemme suivant~:

\th 2.2. Lemme {\rm (\cf le lemme S.~2.1.3 pour le cas complexe)}|Soient $X$
un espace analytique r\'eel normal, $x$ un point de $X$, $f$ et $g$ deux
\'el\'ements de $\oc_{X,x}$ tels que l'id\'eal $\sqrt{(g)}$ soit r\'eel et que
$f\notin (g)$~; il existe alors un arc analytique r\'eel $h :(]-1,1[,0)\to
(X,x)$ tel que $v(f\circ h)< v(g\circ h)$.|

\dem $X_{\reg}$ d\'esignera l'ouvert form\'e des points $y$ de $X$ o\`u
l'anneau local $\oc_{X,y}$ est r\'egulier (dans un voisinage de $x$, cet
ouvert co\"{\i}ncide avec l'ouvert des points lisses de dimension
$d=\dim\oc_{X,x}$).

Quitte \`a restreindre $X$, on peut supposer qu'il existe un morphisme de
r\'esolution des singularit\'es (\cf [H]) \ie, un morphisme analytique r\'eel
$\pi :X'\to X$ propre et surjectif tel que $X'$ soit lisse et que $\pi |\pi
^{-1}(X_{\reg}):\pi ^{-1}(X_{\reg})\to X_{\reg}$ soit un isomorphisme.

On raisonne maintenant comme dans [B-R], Section 2, lemme 3.

La fonction ``m\'eromorphe" $f/g$ a un lieu polaire $P$ non vide dans $X$ (car
$f/g\notin\oc_{X,x}$ par hypoth\`ese) dont le germe en $x$ est r\'eunion de
certaines composantes irr\'eductibles du germe $Z(g)$ d\'efini par $(g)$ (car
si $\wtx$ d\'esigne un complexifi\'e d'un voisinage de $x$ dans $X$ qui soit
un espace normal, et $\tf$ et $\tig$ des extensions de $f$ et $g$ \`a $\wtx$,
la fonction m\'eromorphe $\tf / \tig$ a un lieu polaire de codimension 1 dans
$\wtx$). $P$ est donc de codimension r\'eelle 1 dans $X$ au voisinage de $x$,
car $\sqrt{(g)}$ est par hypoth\`ese un id\'eal r\'eel ce qui implique que
tous les facteurs irr\'eductibles sont r\'eels (\cf 1.2)~; il en r\'esulte que
$P \cap X_{\reg}\ne \emptyset$, car $X$ \'etant normal est lisse en codimension 1.

Soit $x'\in\pi ^{-1}(x)\cap \ol{\pi ^{-1}(P\cap X_{\reg})}$~; comme $X'$ est
lisse, l'anneau $\oc_{X',x'}$ est factoriel et l'on peut \'ecrire~: $f\circ
\pi /g\circ \pi =\alpha /\beta $ dans le corps des fractions de $\oc_{X',x'}$,
avec $\alpha $ et $\beta $ premiers entre eux. $\beta $ s'annule alors sur
$\pi ^{-1}(P\cap X_{\reg})$ au voisinage de $x'$, car $\pi ^{-1}(P\cap
X_{\reg})$ fait partie du lieu polaire de la fonction $\alpha /\beta $ puisque
$\pi |\pi ^{-1}(X_{\reg})$ est un isomorphisme. D'autre part, si $P_1$ est une
composante irr\'eductible analytique locale de $\pi ^{-1}(P)$ en $x'$ telle
que $P_1\cap\pi ^{-1}(X_{\reg})\ne \emptyset$, $\alpha $ ne peut s'annuler
identiquement sur $P_1$, car $P_1$ \'etant de codimension r\'eelle 1 dans
$X'$, cela serait contradictoire avec le fait que $\alpha $ et $\beta $ sont
premiers entre eux (\cf proposition 1.2~: $\alpha $ et $\beta $ seraient tous
deux divisibles par un g\'en\'erateur de l'id\'eal $I(P_1)$).

On peut alors choisir par le ``curve selection lemma" (\cf [M]) un arc
analytique $h'=(]-1,1[,0)\to (X',x')$ tel que $\beta \circ h'\equiv 0$~; on a
alors $v(\beta \circ h')=+\infty$ et $v(\alpha \circ h')<+\infty$~; si
maintenant $h''$ est un arc analytique ayant un contact suffisamment grand avec
$h'$ on aura (exactement comme dans la d\'emonstration du lemme S.~2.1.3)~:
$v(\alpha \circ h'')<v(\beta \circ h'')<+\infty$, soit $v(f\circ\pi \circ
h'')< v(g\circ \pi \circ h'')$, d'o\`u le r\'esultat cherch\'e en posant
$h=\pi \circ h''$. C.Q.F.D.

3) $\Rightarrow$ 4)~: \'Etant donn\'e un espace analytique r\'eel $X$ et un
compact $K\subset X$, il faut montrer qu'il existe un morphisme analytique
r\'eel propre $\pi :X'\to X$ v\'erifiant les propri\'et\'es a), b) et c) de la
proposition 3)~; le probl\`eme est local en $X$, car si pour tout $x\in K$ on
trouve un voisinage $U_x$ de $x$ et un morphisme~: $X'_{U_x}\to U_x$
satisfaisant aux conditions demand\'ees, on prendra pour $X'$ la somme
disjointe des $X'_{U_i}$, o\`u $(U_i)$ est un recouvrement fini de $K$ extrait
du recouvrement $(U_x)$.

Soit donc $x\in K$~: on utilise ``d\'esingularisation I" ([H], 5.10) pour
trouver un voisinage $U$ de $x$ et un morphisme propre et surjectif $\pi
_1:X''\to U$ avec $X''$ lisse, et ``d\'esingularisation II" ([H], 5.11) qui
permet pour chaque point $x''\in\pi ^{-1}_1(x)$ de trouver un voisinage $V$
de $x''$ et un morphisme $\pi _2:X'_V\to V$ propre et surjectif tel que $X_V$
soit lisse et $\ic\oc_{X'}$ un diviseur \`a croisements normaux, ce qui
entra\^{\i}ne \'evidemment que $\forall x'\in X'_V$, l'id\'eal
$\sqrt{\ic\oc_{X'_V,x'}}$ est r\'eel. Il suffit alors de prendre pour $X'$ la
somme disjointe des $X'_{V_i}$ correspondant \`a un recouvrement fini $(V_i)$
de $\pi ^{-1}_1(x)$.

4) $\Rightarrow$ 1)~: Soit $\pi :X'\to X$ un morphisme v\'erifiant les
propri\'et\'es de la condition 3)~; comme par hypoth\`ese il existe un
voisinage ouvert de $\pi ^{-1}(K)$ dans $X'$ tel que
$f^q\oc_{U'}\in\ic^p\oc_{U'}$, il existe un voisinage $U''$ de $\pi ^{-1}(K)$
et une constante $C$ tels que~: $|f\circ\pi (x')|^q\le C\supp_{g\in\Gamma
(U,\ic)}|g\circ\pi (x')|^p$ $\hbox{\rm pour tout}\  x'\in U''$, d'o\`u le r\'esultat puisque
$\pi $ \'etant propre, $\pi (U'')$ contient un voisinage de $K$. C.Q.F.D.

\titre 3. Applications aux exposants de \L ojasiewicz et compl\'ements|

Je vais d'abord montrer un th\'eor\`eme analogue au corollaire S.~6.4
montrant que les exposants de \L ojasiewicz sont rationnels.

\section 3.1|Soient $(X,\oc_X)$ un espace analytique r\'eel, $\ic$ un
$\oc_X$-id\'eal coh\'erent, $f\in\Gamma (X,\oc_X)$ et $K$ un sous-ensemble
compact de $X$~; on d\'efinit de la m\^eme mani\`ere qu'en S.~6.1 l'exposant $\theta
_K(f,\ic)$ comme la borne inf\'erieure de l'ensemble des $\theta \in\rb_+$
tels qu'il existe un voisinage ouvert $U$ de $K$ dans $X$ et une constante
$C\in\rb_+$ avec~:
$$|f(x)|^\theta \le C\sup_{g\in\Gamma (U,\ic)}|g(x)|\quad\hbox{\rm pour tout}\ x\in U.$$
(Dans le cas o\`u $K=\{x\}$, $\theta _{(x)}(f,\ic)$ \'etait not\'e $\alpha
(f,\ic)$ dans [B-R]).

Nous poserons d'autre part $\tit_K(f,\ic)=\theta _K(\tf,\tic)$, $\tf$ et
$\tic$ \'etant des extensions de $f$ et $\ic$ \`a un complexifi\'e $\wtx$ de
$X$~; on a toujours $\theta _K(f,\ic)\le\tit_K(f,\ic)$, et $\tit_K(f,\ic)$
est toujours un nombre rationnel (S.~6.4).

Dans le cas r\'eel, on a le th\'eor\`eme suivant~:

\th 3.2. Th\'eor\`eme|

a) $\theta _K(f,\ic)\in\qb^+\cup\{+\infty\}$.

b) Il existe un voisinage $U$ de $K$ dans $X$ et une constante $C\in\rb_+$
tels que~:
$$|f(x)|^{\theta _K(f,\ic)}\le C\sup_{g\in\Gamma (U,\ic)}|g(x)|\quad\hbox{\rm pour tout}\ x\in U.$$|

Je n'\'ecrirai pas la d\'emonstration de ce th\'eor\`eme~: il suffit en effet
pour la partie a) de recopier la d\'emonstration du th\'eor\`eme S.~4.1.6,
$\bnu^K_I(f)$ \'etant remplac\'e par $1/\theta _K(f,\ic)$ et l'\'eclatement
normalis\'e par un morphisme analytique r\'eel $\pi :X'\to X$ satisfaisant aux
conditions du th\'eor\`eme 2.1.3~; et pour la partie b) de recopier la
d\'emonstration du th\'eor\`eme S.~6.3.

\rque 3.3. Remarque|Dans [B-R], nous avions pos\'e~:
$$\alpha _K(f,\ic)=\inf\big\{\alpha \in\rb_+:\exists
C>0\hbox{~avec~}|f(x)|^\alpha \le C\sup_{g\in \Gamma (K,x)}|g(x)|,\forall x\in
K\big\}$$
et montr\'e que si $K$ est sous-analytique dans $X$, $\alpha _K(f,\ic)$ est un
nombre rationnel~; ce r\'esultat n'a pas de rapport avec le th\'eor\`eme 3.2,
et sa d\'emonstration est tr\`es diff\'erente~: on d\'emontre que pour calculer
$\alpha _K(f,\ic)$ (dans le cas o\`u $K$ est sous-analytique), on peut toujours
se restreindre \`a un arc analytique, alors que c'est faux en g\'en\'eral pour
$\theta _K(f,\ic)$.

\section 3.4|\'Etudions maintenant la question suivante (d\'ej\`a envisag\'ee
dans [B-R]) qui se pose de mani\`ere naturelle~: sous quelles conditions peut-on
affirmer que $\theta _K(f,\ic)=\tit_K(f,\ic)$ (et donc que $\theta
_K(f,\ic)=1/\bnu^K_\ic(f)$)~?

Pour simplifier, nous supposerons dor\'enavant $K=\{x\}$, et poserons $\theta
_{\{x\}}(f,\ic)=\theta (f,I)$ o\`u $I=\ic\cdot\oc_{X,x}$.

\th 3.5. D\'efinition|Soient $(X,\oc_X)$ un espace analytique r\'eel, $x\in
X$. On dit qu'un id\'eal $I\subset\oc_{X,x}$ est r\'eellement r\'eel si
$\hbox{\rm pour tout}\  f\in\oc_{X,x}$, on a l'\'egalit\'e~:
$$\theta (f,I)=\tit(f,I)~.$$|

On a montr\'e dans [B-R] (proposition II.3) que si $X$ est normal et $I$
principal, $I$ est r\'eellement r\'eel si et seulement si $\sqrt{I}$ est un
id\'eal r\'eel, et donn\'e des conditions suffisantes dans le cas o\`u $I$
est engendr\'e par une suite r\'eguli\`ere.

\section 3.6|Une d\'esingularisation \`a la Hironaka d'un id\'eal
$I\subset\oc_{X,x}$ est par d\'efinition un morphisme $\pi :X'\to U$ (o\`u
$U$ est un voisinage convenable de $x$ dans $X$) ayant les propri\'et\'es
suivantes~:

a) $\pi $ est propre et surjectif.

b) $X'$ est lisse et $I\oc_{X'}$ est un diviseur \`a croisements normaux (on
consid\`ere par abus de langage $I$ comme un id\'eal de $\oc_U$).

c) $\pi $ est compos\'e d'une suite finie d'\'eclatements de sous-vari\'et\'es
lisses.

Soit $\wtx$ un complexifi\'e de $X$~; nous noterons $\tilde\pi :\wtx'\to\wtu$
le morphisme analytique complexe propre obtenu en faisant \'eclater les
sous-vari\'et\'es lisses complexifi\'ees des sous-vari\'et\'es lisses que
l'on fait \'eclater pour obtenir le morphisme $\pi $~; si $I\subset\oc_{X,x}$
est un id\'eal, nous poserons $\wt I=I\oc_{\wtx,x}$ (avec
$\oc_{\wtx,x}\simto\oc_{X,x}\otimes_{\rb}\cb$).

On a alors le th\'eor\`eme suivant~:

\th 3.7. Th\'eor\`eme|Soit $I\subset\oc_{X,x}$ un id\'eal~; supposons qu'il
existe une d\'esingularisation \`a la Hironaka de $I$~: $X'\build\la^\pi \fin
U$ telle que l'on ait~:
$$\wt I\oc_{\wtx'}\simto I\oc_{X'}\otimes_{\oc_{X'}}\oc_{\wtx'}~.$$
Alors l'id\'eal $I$ est r\'eellement r\'eel.|

\rque 3.8. Exemples|

a) Soient $X=\rb^2$, 
$I=(x^2+y^2)\subset\rb\{x,y\}$, et $\pi :X'\to X$ l'\'eclatement de l'origine $(0,0)$ de $\rb^2$.

Si $V$ d\'esigne la carte de l'\'eclatement $\pi $ avec coordonn\'ees $x'$ et
$y'$ d\'efinies par $\cases{
x=x'\cr
y=x'y'\cr}$, on a
{\parindent=0pt
$$\leqalignno{
I\Gamma (X',V)&=(x'^2)\cr
\noalign{d'o\`u}
I\Gamma (X',V)\otimes_{\Gamma (X',V)}\Gamma (\wtx',\wt V)&=(x'^2)\cr
\noalign{alors que}
\wt I\Gamma (\wtx',\wt V) &=\big(x'^2(y'+i)(y'-i)\big),\cr}$$}
ce qui montre que la d\'esingularisation $\pi $ ne satisfait pas \`a
l'hypoth\`ese du th\'eor\`eme 3.7 (il est d'ailleurs imm\'ediat de voir que
$I$ n'est pas r\'eellement r\'eel).

b) L'id\'eal $I=(x^2+y^2,y^5)\subset\rb\{x,y\}$ n'est pas non plus
r\'eellement r\'eel (car $\theta (y,I)=2$ et $\tit(y,I)=5$) bien que
contrairement \`a l'exemple pr\'ec\'edent $\sqrt{I}$ soit r\'eel (on a en
effet $\sqrt{I}=(x,y)$).

c) En revanche, l'id\'eal $I=(x^4,y^4)\subset\rb\{x,y\}$ est r\'eellement
r\'eel~: une d\'esingularisation satisfaisant aux hypoth\`eses du
th\'eor\`eme 3.7 est fournie par l'\'eclatement de l'origine~; on peut
remarquer que pourtant $I$ ne satisfait pas au crit\`ere de la proposition II.5
de [B-R].

\rque D\'emonstration du th\'eor\`eme 3.7|Soit $\pi :X'\to X$ une
d\'esingularisation de $I$ satisfaisant aux hypoth\`eses de 3.7~; supposons
que $I$ soit engendr\'e par $(g_1\ld g_p)$.

Il est clair qu'il suffit de montrer que si $f\in\oc_{X,x}$ est telle qu'il
existe un voisinage $V$ de $x$ et une constante $C$ avec $|f(y)|\le
C\supp_{1\le i\le p}|g_i(y)|$ $\forall y\in V$, $f$ est entier sur $I$ (\cf
S.~1.1).

Supposons donc que $|f(y)|\le C\supp_{1\le i\le p}|g_i(y)|$~; on en d\'eduit
$f\oc_{X'}\subset I\oc_{X'}$ par le th\'eor\`eme 2.1~; on a donc
{\parindent=0pt
$$\leqalignno{
f\oc_{\wtx'}&\in I\oc_{X'}\otimes_{\oc_{X'}}\oc_{\wtx'}\cr
\noalign{d'o\`u}
f\oc_{\wtx'}&\in\tilde I\oc_{\wtx'}\cr}$$}
\`a cause de l'hypoth\`ese~; ceci implique que $f$ est entier sur $\wt I$
dans $\oc_{\wtx,x}$ (th\'eor\`eme S.~2.1) d'o\`u imm\'ediatement que $f$ est
entier sur $I$ dans $\oc_{X,x}$. C.Q.F.D.

\titre R\'ef\'erences|
\lettre B-R|

\livre B-R| Bochnack J. et Risler J-J.|Sur les exposants de \L ojasiewicz|Comment. Mat.
Helvetici {\bf 50}|1975|

\livre H|Hironaka h|Introduction to real-analytic sets and real-analytic
maps| Quaderni dei Gruppi di Ricerca Matematica del Consiglio Nazionale delle Ricerche. Istituto Matematico "L. Tonelli" dell'Universit\`a di Pisa, Pisa| 1973|

\livre M|Milnor j|Singular points of complex hypersurfaces|Annals of Math. Studies no. 61 Princeton University Press|1968|

\livre R|Risler j.j|Le th\'eor\`eme des z\'eros en g\'eom\'etrie
alg\'ebrique et analytique r\'eelles|\bsmf|1976|  
\par\medskip\noindent
\titre Sept compl\'ements au s\'eminaire| Cette partie pr\'esente quelques travaux en rapport direct avec le contenu du s\'eminaire qui sont venus \`a notre
connaissance depuis, ce qui n'exclut pas
que certains aient \'et\'e \'ecrits bien avant! En particulier une partie du s\'eminaire appara\^\i t \it a posteriori \rm comme une traduction en g\'eom\'etrie analytique de r\'esultats de Rees (voir [R1], [R2], [R3]),   Northcott-Rees (voir [N-R])  et Nagata ([N]) cit\'es dans la bibliographie compl\'ementaire, dont nous ignorions l'existence.\par\noindent Nous ne pr\'etendons
nullement \^etre exhaustifs. En particulier nous ne mentionnerons ici que quelques travaux post\'erieurs au s\'eminaire concernant la d\'ependance int\'egrale sur les id\'eaux, et aucune des applications \`a l'\'equisingularit\'e, pour lesquels nous renvoyons au travaux de B. Teissier cit\'es dans la m\^eme bibliographie, ni les travaux sur la d\'ependance int\'egrale sur les modules  pour lesquels nous renvoyons \`a ceux de  T. Gaffney et S. Kleiman, \'egalement cit\'es, ainsi qu'\`a ceux de Kleiman-Thorup. Disons seulement que pour les familles d'hypersurfaces (voir [T1]) ou d'intersections compl\`etes \`a singularit\'es isol\'ees (voir [G-K 1]), les conditions de Whitney, la condition $a_f$ de Thom et la condition $w_f$ s'expriment toutes par le fait que le $\bnu$ de certains id\'eaux ou modules dans la d\'efinition desquels interviennent des mineurs jacobiens des \'equations de la fermeture des strates, par rapport \`a d'autres du m\^eme type, est $\geq 1$ ou $>1$. Signalons les travaux de Rees et B\"oger cit\'es dans la bibliographie compl\'ementaire, ainsi que le livre [H-I-O] et en particulier l'appendice de B. Moonen. La d\'ependance int\'egrale sur les id\'eaux et modules est \'etudi\'ee d'un point de vue alg\'ebrique dans le livre [V] de W. Vasconcelos ainsi que dans le livre r\'ecent [Hu-S] de C. Huneke et I. Swanson. Pour les rapports de la d\'ependance int\'egrale avec la "Tight closure" de Hochster, nous renvoyons \`a [Hu].\par\medskip\noindent
{\bf Compl\'ement 1: L'in\'egalit\'e de \L ojasiewicz pour le gradient}
\par\noindent
Nous avons d\'ej\`a discut\'e dans le \S 6 du s\'eminaire du rapport entre le $\overline\nu$ et l'in\'egalit\'e de
\L ojasiewicz. Nous allons montrer que ce dictionnaire est efficace en l'utilisant pour d\'emontrer
un des r\'esultats c\'el\`ebres de la th\'eorie des in\'egalit\'es de \L ojasiewicz, d\^u \`a celui-ci:\par\noindent
{\bf Th\'eor\`eme} (\L ojasiewicz): Soit $f(x_1,\ldots ,x_n)$ une fonction analytique r\'eelle d\'efinie dans un
voisinage de l'origine dans ${\bf R}^n$ et telle que $f(0)=0$. Il existe un voisinage $U$ de $0$ dans ${\bf R}^n$, un nombre r\'eel $\theta\
,\ 0<\theta\ <1$ et une constante $C>0$ tels que l'on ait
pour tout $x\in U$ l'in\'egalit\'e:
$$\vert \hbox{grad}f(x)\vert\geq C\vert f(x)\vert^\theta .$$

Remarquons d'abord qu'il suffit de prouver la m\^eme in\'egalit\'e pour la fonction
complexifi\'ee de $f$, que nous noterons encore $f$. Remarquons aussi que le point est
l'in\'egalit\'e $\theta <1$; l'existence d'une in\'egalit\'e se d\'eduit du th\'eor\`eme des z\'eros de Hilbert et du fait qu'au voisinage de $0$ la fonction analytique $f(x)$ s'annule sur le lieu des z\'eros de son gradient. Nous allons utiliser le mode de calcul du $\overline
\nu$ donn\'e au \S 5. Soit $W$ un voisinage de $0$ dans $\C^n$ o\`u $f$ converge, choisissons y
des coordonn\'ees holomorphes $(z_1,\ldots ,z_n)$ et consid\'erons l'\'eclatement 
$$\pi\colon Z\to W$$
de l'id\'eal jacobien $$j(f){\cal O}_W=({\partial f \over {\partial z_1}},\ldots ,{\partial f \over
{\partial z_n}}){\cal O}_W$$ et le morphisme
$$\overline \pi\colon \overline Z\to Z\to W$$
compos\'e de $\pi$ et de la normalisation $n\colon \overline Z\to Z$.\par
Puisque l'espace $\overline Z$ est normal, son lieu singulier est de codimension $\geq 2$, et
puisque l'image inverse par $\overline \pi$ de l'id\'eal $j(f){\cal O}_W$, c'est \`a dire l'id\'eal de
${\cal O}_{\overline Z}$ engendr\'e par les $({\partial f \over {\partial z_1}}\circ \overline
{\pi},\ldots ,{\partial f \over {\partial z_n}}\circ \overline {\pi})$ est par construction
localement principal et engendr\'e en tout point de $\overline Z$ par l'un des ${\partial f \over
{\partial z_j}}\circ \overline {\pi}$, cet id\'eal d\'efinit un sous-espace $D$ de codimension 1
dans $\overline Z$, qui n'a qu'un nombre fini de composantes irr\'eductibles si nous avons
pris la pr\'ecaution de choisir le voisinage $W$ relativement compact, puisque le morphisme
$\overline \pi$ est propre.\par
Soit $D_{red}=\bigcup_iD_i$ la d\'ecomposition ensembliste de $D$ en composantes irr\'eductibles.
Celles-ci sont de codimension 1 et donc chacune contient un ouvert analytique dense $V  _i$ en
chaque point $z$ duquel on a ceci:\par\noindent
1) L'espace $\overline Z$ est non singulier en $z$ ainsi que l'ensemble analytique $D_{red}$
sous-jacent \`a l'espace analytique $D$, et $D_{red}$ est donc d\'efini au voisinage de $D$ par
une \'equation $v=0$, o\`u $v$ est une coordonn\'ee locale sur $\overline Z$ en $z$.\par\noindent
2) Choisissons des coordonn\'ees locales $b_1,\ldots ,b_{n-1}$ sur $D$ en $z$; alors $v,b_1,\ldots
,b_{n-1}$ est un syst\`eme de coordonn\'ees locales sur $\overline Z$ en $z$ et l'on a une \'ecriture \par\noindent
$(f\circ\overline \pi)_z=Av^{\mu_i}$, o\`u $A\in \C\lbrace v,b_1,\ldots ,b_{n-1}\rbrace$ avec
$A(0,\ldots ,0)\neq 0$ et $\mu_i>0$ puisque $f$ s'annule l\`a o\`u toutes ses d\'eriv\'ees
s'annulent. \par\noindent
3) L'id\'eal $(j(f){\cal O}_{\overline Z})_z$ est engendr\'e par un des germes $ ({\partial f \over
{\partial z_j}}\circ \overline {\pi})_z$, disons $({\partial f \over {\partial
z_1}}\circ \overline {\pi})_z$, et l'on a \par\noindent
  $({\partial f \over {\partial
z_1}}\circ \overline {\pi})_z=Bv^{\nu_i}\ \ \ \ \hbox{o\`u}\ B\in\C\lbrace v,b_1,\ldots ,b_{n-1}\rbrace\ \ 
\hbox{\rm avec }\  B(0,\ldots ,0)\neq 0$ et $\nu_i >0$.\par\noindent
Remarquons que puisque $D_i$ est irr\'eductible, les entiers $\mu_i, \nu_i$ ne d\'ependent pas du
point $z\in V_i$ choisi. La r\`egle de Leibniz donne:
$${\partial{(f\circ \overline \pi)_z}\over\partial v}=\sum_{j=1}^n ({\partial f \over
{\partial z_j}}\circ \overline {\pi})_z{\partial{(z_j\circ \overline \pi)_z}\over
\partial v}=v^{\mu_i-1}(\mu_iA+v{\partial A\over \partial v}).$$
On en d\'eduit aussit\^ot l'in\'egalit\'e $$\mu_i -1\geq \nu_i$$
Cela donne, posant $\theta = \hbox{sup}_i{\nu_i\over
\mu_i}={1\over\overline{\nu}_{j(f)}(f)}$, l'in\'egalit\'e $\theta<1$ et d'autre part par les
m\'ethodes du \S 2 du s\'eminaire, il vient en posant $\theta={p\over q}$ l'inclusion
$f^p\in \overline{j(f)^q}$, ce qui \'equivaut d'apr\`es le th\'eor\`eme 7.2 du s\'eminaire \`a $\vert f(z)\vert^p\leq C'\vert
\hbox{grad}f(z)\vert^q$ pour $z$ assez proche de $0$ et donc au r\'esultat
cherch\'e.\par\noindent
{\bf Remarque :} Nous venons de prouver l'in\'egalit\'e $\bnu_{j(f)}(f)>1$. Dans le m\^eme ordre d'id\'ee, notant $(z)*j(f)$ l'id\'eal $(z_1\partial f/\partial z_1,\ldots, z_n\partial f/\partial z_n){\cal O}_{\C^n,0}$, l'in\'egalit\'e $\bnu_{(z)*j(f)}(f)\leq 1$ est facile \`a v\'erifier par restriction \`a une droite de l'espace ambient, en utilisant le fait que $\bnu$ ne peut qu'augmenter par passage au quotient. Il est prouv\'e dans la section 0.5 de [T 1] que l'on a pour tout $f\in m\subset \C\{z_1,\ldots ,z_n\}$ l'in\'egalit\'e $\bnu_{(z)*j(f)}(f)\geq 1$, et on a donc $\bnu_{(z)*j(f)}(f)= 1$, ce qui est une version en alg\`ebre asymptotique de l'identit\'e d'Euler pour les polynomes homog\`enes.\par\noindent
Si l'on veut un r\'esultat n'impliquant que des id\'eaux ind\'ependants des coordonn\'ees, on peut en utilisant le th\'eor\`eme de Bertini id\'ealiste de [T1] prouver par restriction \`a une droite assez g\'en\'erale l'in\'egalit\'e $\bnu_{m.j(f)}(f)\leq 1$ et d\'eduire du r\'esultat pr\'ec\'edent l'\'egalit\'e $\bnu_{m.j(f)}(f)= 1$.\par 
Il faut ajouter que la d\'efinition alg\'ebrique de l'exposant de \L ojasiewicz donn\'ee dans le s\'eminaire a \'et\'e \'etendue au cas analytique r\'eel, ou semi-alg\'ebrique, par Fekak ([F1], [F2]) en utilisant une d\'efinition due \`a Brumfiel [Br] des relations de d\'ependance semi-int\'egrale dans le cas r\'eel. Cela r\'epond au voeu exprim\'e par Risler au d\'ebut de son appendice.\par\medskip\noindent
{\bf Compl\'ement 2: L'ordre $\bnu$ et le polyg\^one de Newton}\par 
Dans l'appendice au \S 4 du s\'eminaire, nous donnons une seconde d\'emonstration de la rationalit\'e de $\bnu$ qui repose sur le fait que $\bnu$ peut \^etre vu comme tropisme critique d'une certaine installation, c'est-\`a-dire comme pente d'un c\^ot\'e d'un polyg\^one de Newton g\'en\'eralis\'e. Dans le cas d'un id\'eal $I$ primaire pour l'id\'eal maximal d'une alg\`ebre analytique r\'eduite ${\cal O}$ de dimension pure, le \S 4 de [T2] associe \`a tout \'el\'ement $f\in {\cal O}$ ou \`a tout id\'eal $J\subset {\cal O}$ un polygone de Newton tel que $(\bnu_I(f))^{-1}$ ou $(\bnu_I(J))^{-1}$ apparaisse parmi les oppos\'es des pentes de ses c\^ot\'es. Ce polygone ne d\'epend que de la cl\^oture int\'egrale de $I$.\par
Etant donn\'es $h, \ell\in \R_{\geq 0}$, notons $\bigl\{ {{\ell}\over{\overline{\ h \ }}} \bigr\}$ le polygone de Newton \'el\'ementaire (poss\'edant au plus un c\^ot\'e compact)  ayant pour sommets les points $(0,h)$ et $(\ell, 0)$. C'est le bord de l'enveloppe convexe de $((0,h)+\R^2_{\geq 0})\bigcup ((\ell, 0)+\R^2_{\geq 0}) $. Le mono\"\i de pour l'addition de Minkowski (point par point) de tous les polygones de Newton rencontrant les deux axes de coordonn\'ees est engendr\'e par les polygones \'el\'ementaires. Si l'on autorise $h$ et $\ell$ \`a prendre la valeur $+\infty$, en convenant que $\bigl\{{{\ell}\over{\overline{\ \infty\ }}}\bigr\}$ est constitu\'e de deux demi-droites parall\`eles aux axes et se rencontrant au point $(\ell,0)$ et de m\^eme pour $\bigl\{{{\infty}\over{\overline{\ h\ }}}\bigr\}$ et le point $(0,h)$, on engendre le semigroupe de tous les polygones de Newton, avec l'\'el\'ement neutre consistant en la r\'eunion des deux demi-axes positifs.\par
 Soient $X$ un espace analytique complexe r\'eduit et \'equidimensionel, $x\in X$ et $I$ un id\'eal primaire pour l'id\'eal maximal de ${\cal O}_{X,x}$. Consid\'erons le diviseur exceptionnel $D$ de l'\'eclatement normalis\'e $\overline E_I\colon Z\to X$ de $I$ dans $X$. Chacune des composantes irr\'eductibles $D_k$ du diviseur compact $D$ d\'etermine une fonction d'ordre $v_k$ sur ${\cal O}_{X,x}$; l'ordre de $f\in {\cal O}_{X,x}$ est l'ordre d'annulation de $f\circ \overline E_I$ le long de $D_k$. Puisque l'\'eclatement normalis\'e se factorise \`a travers la normalisation $n\colon \overline X\to X$ qui s\'epare les composantes analytiques de $X$ en $x$, cette fonction d'ordre est en fait compos\'ee d'une valuation divisorielle sur l'une des composantes analytiquement irr\'eductibles $X_{j(i)}$ du germe $(X,x)$ et de la surjection ${\cal O}_{X,x}\to {\cal O}_{X_{j(i)},x}$. On d\'efinit l'ordre d'un id\'eal $J$ de ${\cal O}_{X,x}$ comme l'infimum des ordres de ses \'el\'ements. Par ailleurs les composantes $D_k$ du diviseur exceptionnel sont des vari\'et\'es projectives plong\'ees par le faisceau tr\`es ample $I{\cal O}_Z/(I{\cal O}_Z)^2$, et les vari\'et\'es r\'eduites sous jacentes $\vert D_k\vert$ aussi; on peut donc parler de leur degr\'e $\hbox{\rm deg}\vert D_k\vert$.\par\noindent

On peut alors d\'efinir pour $g\in{\cal O}_{X,x}$ le \it Polygone de Newton de $g$ par rapport \`a $I$ \rm par 
$$N_I(g)=\sum_k\hbox{\rm deg}\vert D_k\vert \bigl\{{{v_k(I)}\over{\overline{\ v_k(g)\ }}}\bigr\},$$
et de m\^eme le \it Polygone de Newton de $J$ par rapport \`a $I$ \rm par 
$$N_I(J)=\sum_k\hbox{\rm deg}\vert D_k\vert \bigl\{{{v_k(I)}\over{\overline{\ v_k(J)\ }}}\bigr\}.$$\par\noindent
Il r\'esulte du \S 4 du s\'eminaire que $\bnu_I(g)$ est la valeur absolue $h/\ell$ de la pente du c\^ot\'e le plus horizontal (ou dernier c\^ot\'e) du polygone de Newton $N_I(g)$.\par
On peut montrer (voir [R3] et [T2], [T3]) que l'application qui \`a $g\in {\cal O}_{X,x}$ associe la longueur $\sum_k\hbox{\rm deg}\vert D_k\vert v_k(g)$ de la projection orthogonale de $N_I(g)$ sur l'axe vertical n'est autre que la \it degree function \rm de Pierre Samuel et David Rees, qui est d\'efinie comme la multiplicit\'e de l'image $(I +g{\cal O}_{X,x})/g{\cal O}_{X,x}$ de l'id\'eal $I$ dans ${\cal O}_{X,x}/g{\cal O}_{X,x}$, tandis que la longueur $\sum_k\hbox{\rm deg}\vert D_k\vert v_k (I)$ de sa projection sur l'axe horizontal est \'egale \`a la multiplicit\'e au sens de Samuel de l'id\'eal primaire $I\subset {\cal O}_{X,x}$ (voir [T2] et [R-S]). On peut observer que lorsque $g$ est dans $I$ le quotient des longueurs des deux projections de $N_I(g)$, ou mieux encore $N_I(g)$ lui-m\^eme \`a homoth\'etie pr\`es, est une mesure du d\'efaut de superficialit\'e (au sens de Samuel) de $g$ par rapport \`a $I$. Lorsque $g\in I$ est superficiel, le polyg\^one $N_I(g)$ n'a qu'un seul c\^ot\'e compact, de pente \'egale \`a $-1$. Cela r\'esulte aussit\^ot de l'\'egalit\'e des multiplicit\'es de $I$ dans ${\cal O}_{X,x}$ et ${\cal O}_{X,x}/g{\cal O}_{X,x}$ et des in\'egalit\'es $v_k(g)\geq v_k(I)$. On peut aussi  le v\'erifier en interpr\'etant [Bon] dans l'\'eclatement normalis\'e de $I$. La  "degree function" est \'etendue  \`a des filtrations noetheriennes et \`a des ${\cal O}_{X,x}$-modules au chapitre 9 de [R11].\par\noindent
Un fait int\'eressant d\'emontr\'e dans [R-S] est que la relation
$$e((I +g{\cal O}_{X,x})/g{\cal O}_{X,x})=\sum_k\hbox{\rm deg}\vert D_k\vert v_k(g)\ \ \hbox{\rm pour tout}\  g\in {\cal O}_{X,x}$$
d\'etermine de mani\`ere unique les coefficients $\hbox{\rm deg}\vert D_k\vert$.\par
Lorsque $g$ est un \'el\'ement g\'en\'eral de l'id\'eal $J$, on a l'\'egalit\'e $N_I(g)=N_I(J)$. Cela sert entre autres \`a montrer (voir [T2], [T3]) que la longueur de la projection verticale du polygone $N_I(J)$ est \'egale \`a la \it multiplicit\'e mixte \rm $e(I^{[d-1]},J^{[1]})$. Ce fait a \'et\'e utilis\'e en dimension 2 par Rees et Sharp dans [R-S].\par
D'autre part on peut \'etendre \`a ce cadre les r\'esultats de la section 5.7  du s\'eminaire. Si $X$ est de Cohen-Macaulay et si $(g_1,\ldots, g_d)$ sont des \'el\'ements de $I$ engendrant un id\'eal qui a m\^eme cl\^oture int\'egrale, en suivant exactement la preuve de 5.7.1, on prouve que si $C_{\hbox{\rm gen}} =\bigcup_{q=1}^r\Gamma_q$ est la d\'ecomposition en composantes irr\'eductibles de la courbe d\'efinie par $d-1$ combinaisons lin\'eaires assez g\'en\'erales\footnote{(*)}{Qui sont un peu abusivement dites "g\'en\'eriques" \`a la fin du \S 5. Le lecteur est encourag\'e \`a consulter [N-R] en se souvenant qu'un id\'eal $I$ est une {\it r\'eduction} d'un id\'eal $J$ si $I\subseteq J$ et $\overline I=\overline J$.} de $(g_1,\ldots, g_d)$ et si l'on note $h_q\colon (\D,0)\to (X,x)$ les arcs analytiques param\'etrant les $\Gamma_q$, on a (voir [T2]) les \'egalit\'es:
$$N_I(g)=\sum_{q=1}^r\bigl\{{{\nu_0(I \circ h_q)}\over{\overline{\ \nu_0(g\circ h_q)\ }}}\bigr\}\ \ \hbox{\rm et}\ \ N_I(J)=\sum_{q=1}^r\bigl\{{{\nu_0(I \circ h_q)}\over{\overline{\ \nu_0(J\circ h_q)\ }}}\bigr\},$$
o\`u $\nu_0$ d\'esigne comme plus haut l'ordre \`a l'origine de $\C$, et l'on a gard\'e les notations de 5.1 du s\'eminaire.
\par
Lorsque $f\colon (\C^n,0)\to (\C,0)$ est un germe de fonction holomorphe \`a singularit\'e isol\'ee, si l'on prend pour $I$ l'id\'eal $j(f)=({{\partial f}\over{\partial z_1}},\ldots ,{{\partial f}\over{\partial z_n}})$ et pour $J$ l'id\'eal maximal $m$ de $\C\{z_1,\ldots ,z_n\}$, le polygone de Newton $N_{j(f)}(m)$ prend le nom de \it polygone de Newton jacobien \rm. Il est d\'emontr\'e dans [T3] que le polygone de Newton jacobien est un invariant d'\'equisingularit\'e \`a la Whitney des hypersurfaces \`a singularit\'e isol\'ee. En particulier c'est un invariant du type topologique des singularit\'es de courbes planes r\'eduites. Dans le cas des branches planes, c'est m\^eme un invariant total du type topologique, et il y a un r\'esultat analogue pour les courbes planes r\'eduites (voir [GB]). Lorsque $n>2$, l'invariance par d\'eformation Whitney-\'equisinguli\`ere est le seul moyen dont on dispose pour prouver que l'exposant de \L ojasiewicz optimal des in\'egalit\'es du gradient $\vert \hbox{\rm grad}(f(z))\vert \geq C_1\vert f(z)\vert^{\theta_1}$ (resp. $\vert \hbox{\rm grad}(f(z))\vert \geq C_2\vert z\vert^{\theta_2}$) pour $\vert z\vert$ assez petit est invariant par de telles d\'eformations.\par

Des travaux r\'ecents de Evelia Garc\'\i a Barroso, Janusz Gwo\'zdziewicz, Tadeusz Krasi\'nski, Andrzej Lenarcik et Arkadiusz P\l oski donnent des exemples de nombres rationnels qui ne peuvent \^etre exposant de \L ojasiewicz pour la seconde in\'egalit\'e du gradient d'une singularit\'e de courbe plane (voir [GB-P], [GB-K-P 1], [GB-K-P 2]), ce qui implique un r\'esultat analogue pour la premi\`ere puisque $\theta_1={{\theta_2}\over{\theta_2+1}}$ d'apr\`es [T3].\par Tout r\'ecemment (voir [GB-G]) une caract\'erisation combinatoire des polygones de Newton jacobiens des branches planes a \'et\'e obtenue. Ce dernier travail montre en particulier que le polygone de Newton jacobien d'un germe de courbe analytique complexe plane permet de d\'ecider si elle est irr\'eductible, en fort contraste avec le polygone de Newton usuel.
\par Enfin, des r\'esultats pr\'ecis sur les relations analogues \`a celle qui vient d'\^etre cit\'ee entre les exposants des deux in\'egalit\'es du gradient pour des fonctions\break analytiques r\'eelles $f\colon (\R^n,0)\to (\R,0)$ telles que $f^{-1}(0)=\{0\}$ se trouvent dans [Gw].\par
Dans la m\^eme veine, consid\'erons une branche plane $(X,0)\subset (\C^2,0)$ donn\'ee param\'etriquement par $x(t)=t^n,\ y(t)=t^m+\cdots$ avec $m\geq n$ et ayant pour caract\'eristique de Puiseux $(\beta_0=n,\beta_1,\ldots, \beta_g)$ (les $\beta_j/n$ sont les exposants caract\'eristiques).\par\noindent Dans l'alg\`ebre $\C\{t,t'\}$ de $(\overline X\times \overline X,\{0\}\times \{0\})$, l'id\'eal $(x(t)-x(t'),y(t)-y(t'))$ qui d\'efinit le sous espace produit fibr\'e $\overline X\times_X \overline X$ s'\'ecrit $(t-t'){\cal N}$, o\`u ${\cal N}$ est un id\'eal primaire correspondant au sous-espace des points doubles du morphisme fini $(\C,0)\to (\C^2,0)$ param\'etrisant $X$. On a d'ailleurs $e({\cal N})=2\delta$ o\`u $\delta=\hbox{\rm dim}_\C \overline{{\cal O}_{X,x}}/{\cal O}_{X,x}$. Si l'on pose $e_0=n,\ e_i=\hbox{\rm pgcd}(\beta_0,\ldots, \beta_i)$ et $m=(t,t')\C\{t,t'\}$, on peut v\'erifier que la courbe ${{t^n-t'^n}\over {t-t'}}=0$ est assez g\'en\'erale pour permettre le calcul du polygone de Newton $N_{\cal N}(m)$ par la m\'ethode expliqu\'ee plus haut, ce qui donne l'\'egalit\'e
$$N_{\cal N}(m)=\sum_{j=1}^g(e_{j-1}-e_j)\bigl\{{{\beta_j-1}\over{\overline{\ \ \ \ 1\ \ \ }}}\bigr\}.$$
On retrouve ainsi en particulier l'\'egalit\'e $2\delta=1-\beta_0+\sum_{j=1}^g(e_{j-1}-e_j)\beta_j$ (voir [Mi], Remark 10.10 et [Z], 3.14) et une interpr\'etation du plus grand exposant de Puiseux comme exposant de \L ojasiewicz puisque $\bnu_{\cal N}(m)=(\beta_g-1)^{-1}$. Pour tout ceci, voir [P-T], [T4], Chap. II, \S 6 et [T7].
 
\par\medskip\noindent
{\bf Compl\'ement 3: Une remarque et un r\'esultat d'Izumi}\par Comme l'a remarqu\'e Izumi dans [I 1], une
d\'emonstration dans le cadre analytique de l'existence d'un nombre r\'eel $b(I)$ tel que pour tout $x$ on ait 
$$\n_I(x)-\nu_I(x)\leq b(I)$$ est implicitement donn\'ee dans le s\'eminaire. En fait, ce r\'esultat
est dans le cadre analytique un corollaire facile du fait prouv\'e dans le s\'eminaire (lemme 4.3.2) que pour
tout entier $q$ l'alg\`ebre gradu\'ee $\overline{{\cal P}^{1\over q}}(I) =\bigoplus_{p\in
\N}\overline{I^{p\over q}}$ est la fermeture int\'egrale dans $A[T^{1\over q}]$ de l'alg\`ebre de
Rees ${\cal P}(I)$. Cette \'egalit\'e implique en effet, puisque en G\'eom\'etrie
analytique les anneaux sont de Nagata, que $\overline{{\cal P}^{1\over q}}(I)$ est un ${\cal
P}(I)$-module gradu\'e de type fini (\it cf. \rm la remarque pr\'ec\'edant 4.3.6), et donc qu'il existe un entier $p_0$ tel que pour $p>p_0$ et tout entier positif $\ell$ on ait
$$\overline{I^{{p+\ell q}\over q}}=I^\ell\overline{I^{p\over q}}$$ d'o\`u r\'esulte l'in\'egalit\'e annonc\'ee. Dans le cas alg\'ebrique cette in\'egalit\'e est due \`a Rees ([R3]) et Nagata ([N]). Rees a donn\'e dans [R12] une interpr\'etation valuative des r\'esultats de ce type dans le cadre g\'en\'eral des anneaux locaux qui sont de Nagata.\par
Izumi a donn\'e dans [I 2] un crit\`ere pour qu'un morphisme injectif $\phi\colon A\to B$ d'alg\`ebres analytiques complexes satisfasse la condition du rang de Gabrielov, qui implique l'injectivit\'e du morphisme $\hat\phi\colon\hat A\to \hat B$ des compl\'et\'es : il faut et il suffit qu'il existe une constante $C> 0$ telle que pour tout $f\in A$ on ait $C\n_{m_A} (f)\geq \n_{m_B} (\phi(f))$. Notons que puisque clairement $\n_{m_B} (\phi(f))\geq\n_{m_A} (f)$, on doit avoir $C\geq 1$. Pour une bonne pr\'esentation de ses r\'esultats nous renvoyons \`a [I 3].\par\noindent

\par\medskip\noindent
 {\bf Compl\'ement 4: G\'en\'eralisation de la d\'efinition du $\n_I(J)$ et de sa rationalit\'e}\par
Dans [C-E-S] les auteurs prouvent le r\'esultat suivant: soient $J_1,\ldots, J_k, I$ des id\'eaux d'un anneau $A$ localement analytiquement non ramifi\'e tels que $J_i\subset \sqrt I$ pour tout $i$, que l'id\'eal $I$ ne soit pas nilpotent et v\'erifie $\bigcap_k I^k=(0)$.\par\noindent Soit $C=C(J_1,\ldots, J_k, I)$ le c\^one  de $\R^{k+1}$ engendr\'e par les $(m_1,\ldots ,m_k,n)\in \N^{k+1}$ tels que $J_1^{m_1}\ldots J_k^{m_k}\subset I^n$. Alors l'adh\'erence du c\^one $C$ est un c\^one polyh\'edral rationnel. \par
Le cas o\`u $k=1$ correspond \`a la rationalit\'e de $\bnu$. On pourrait d\'emontrer ce r\'esultat en appliquant les m\'ethodes du s\'eminaire \`a la compl\'etion des localis\'es de $A$.\par
Il serait int\'eressant d'\'etendre les r\'esultats du \S 4 du s\'eminaire aux alg\`ebres gradu\'ees
$$\bigoplus_{\ell\in \N^k} I_1^{\ell_1}\ldots  I_k^{\ell_k} T_1^{\ell_1}\ldots  T_k^{\ell_k}\subset A[T_1,\ldots  T_k],$$ et en particulier au vu de la section 4.3 du s\'eminaire, \`a leur fermeture int\'egrale dans des alg\`ebres du type $A[T_1^{1/q_1},\ldots , T_k^{1/q_k}]$.
\par\medskip\noindent
 {\bf Compl\'ement 5: {Le $\bnu_{I}(\sqrt I)$ et le type des id\'eaux }\par
 Le \it type \rm d'un id\'eal $I$ de fonctions analytiques en un point de $\C^n$ a \'et\'e introduit par D'Angelo dans [D] comme moyen de mesurer, \'etant donn\'e un domaine $\Omega$ de $\C^n$ dont le bord $\partial \Omega$ est suppos\'e lisse, le contact avec $\partial \Omega$ de germes $h\colon (\C,0)\to (\C^n,p)$ de courbes holomorphes en un point $p\in \partial\Omega$. Cela est li\'e \`a des estim\'ees hypoelliptiques, dont on trouvera un r\'esum\'e au d\'ebut de [Mc-N] et de [H], pour l'\'equation $\overline\partial u=\alpha$ sur $\Omega$.
 \par
 La d\'efinition du type est la suivante: \'etant donn\'es un point $x\in \C^n$ et un id\'eal $I\subset m_{\C^n,x}$, on pose\footnote{(*)}{La d\'efinition originelle du type \'etait en fait $\hbox{\rm sup}_h\bigl\{{{\nu_0(I\circ h)}\over{\nu_0(m\circ h)}}\bigr\}$, c'est \`a dire $(\bnu_I(m))^{-1}$ pour des id\'eaux primaires pour l'id\'eal maximal. Cette d\'efinition-ci est due \`a Heier [H].} (toujours avec les notations de 5.1) $$T_x(I)=\hbox{\rm sup}_h\bigl\{{{\nu_0(I\circ h)}\over{\nu_0(\sqrt I\circ h)}}\bigr\},$$ o\`u $h$ parcourt l'ensemble des germes d'arcs analytiques $h\colon (\db,0)\to (\C^n,x)$ tels que $h(\db)$ ne rencontre le sous-ensemble analytique d\'efini par $I$ qu'en $h(0)=x$ et $\nu_0$ d\'esigne l'ordre \`a l'origine de $(\db,0)$.\par
On peut alors d\'efinir le type d'un faisceau coh\'erent d'id\'eaux ${\cal I}$ sur un espace analytique r\'eduit $X$ en chaque point $x$ du sous-espace d\'efini par ${\cal I}$; c'est le type du germe ${\cal I}_x$. Au vu des r\'esultats du \S 5 du s\'eminaire, \bf le type en $x$ du faisceau d'id\'eaux ${\cal I}$ n'est autre que $(\bnu_{{\cal I}_x}(\sqrt{\cal I}_x))^{-1}$, \rm et on peut donc lui appliquer le reste des r\'esultats du s\'eminaire. En particulier l'article [H-L] retrouve essentiellement le contenu de la section 5.7 du s\'eminaire dans ce cas particulier. Lorsque $I$ est primaire pour l'id\'eal maximal $m=(z_1,\ldots ,z_n)$ de ${\cal O}_{X,x}$, le type $T_x(I)$ est le plus petit exposant possible pour une in\'egalit\'e de \L ojasiewicz
$$\vert g(z)\vert \geq C\vert z\vert^\theta  $$ au voisinage de $x$,  o\`u $g=(g_1,\ldots g_s)$ est un syst\`eme de g\'en\'erateurs de $I$ et $C$ une constante positive.\par D\`es le milieu des ann\'ees 1980, A. P\l oski (voir [P\l 1]) puis J. Chadzy\'nski et T. Krasi\'nski (voir [C-K 1]) s'\'etaient rendu compte que dans le plan les exposants de \L ojasiewicz du type $(\bnu_I(m))^{-1}$ pouvaient s'exprimer en termes de contact de germes de courbes planes. Rappelons que le contact en $0$ du germe $g=0$ avec $f=0$ est d\'efini comme le quotient de la multiplicit\'e d'intersection en $0$ de $f=0$ et $g=0$ par la multiplicit\'e de $f$ en $0$. Lorsque $f$ est irr\'eductible comme l'est l'image d'un arc analytique $h\colon (\D,0)\to (\C^2,0)$, le contact de $g=0$ avec $f=0$ est donc \'egal \`a $\nu_0(g\circ h)/\nu_0(h)$, et si l'on prend le sup sur tous les arcs $h$ de l'infimum lorsque $g$ parcourt les g\'en\'erateurs d'un id\'eal primaire $I$ on retrouve le type de $I$. Si l'on remplace le contact de deux courbes par le contact d'une hypersurface avec une courbe irr\'eductible, cette approche s'\'etend en toute dimension, mais dans le plan on dispose des outils issus du d\'eveloppement de Puiseux pour d\'ecrire le contact. Une des cons\'equences de ces travaux est que l'exposant de \L ojasiewicz par rapport \`a l'id\'eal maximal d'un id\'eal primaire $I$ de ${\cal O}_{\C^2,0}$ est calcul\'e en prenant le supremum sur l'ensemble fini des arcs analytiques $h$ correspondant aux composantes irr\'eductibles des \'el\'ements d'un syst\`eme de g\'en\'erateurs $(g_i)$ de $I$ des $\nu_0(g_i\circ h)/\nu_0(h)$ qui sont finis (voir [C-K 1], [C-K 2], [P\l 3]).
Ce r\'esultat a \'et\'e red\'ecouvert dans le langage du type des id\'eaux par  J.D. McNeal et A. N\'emethi dans [Mc-N], avec une preuve diff\'erente.\par Par ailleurs, apr\`es un premier r\'esultat de P. Philippon dans cette direction (voir [P]), Ein-Lazarsfeld (voir [E-L]) et M. Hickel (voir [Hi]) puis G. Heier (voir [H]) ont prouv\'e diverses formes locales et globales du th\'eor\`eme des z\'eros de Hilbert effectif en s'appuyant sur la d\'ependance int\'egrale sur les id\'eaux et le calcul de $\bnu$ par \'eclatement normalis\'e comme dans le \S 4 du s\'eminaire.\par Pour Hickel, il s'agit de prouver \it sur un corps $k$ quelconque \rm une conjecture de Berenstein et Yger concernant l'appartenance effective \`a un id\'eal $I$ de $k[X_1,\ldots ,X_n]$ d'une puissance d'un \'el\'ement de $\overline I$. Etant donn\'es un syst\`eme de g\'en\'erateurs $p_1,\ldots ,p_m$ de $I$ et un polyn\^ome $p\in \overline I$, on cherche \`a \'ecrire une relation $p^s=\sum_{i=1}^mq_ip_i$ avec $s\leq \hbox{\rm min}(m, n+1)$ et des bornes sur le degr\'e des $q_ip_i$ en fonction du degr\'e des $p_i$ et de la dimension $n$.  Pour Heier il s'agit, \'etant donn\'e un faisceau coh\'erent d'id\'eaux ${\cal I}$ sur une vari\'et\'e projective complexe non singuli\`ere $X$, de trouver le plus petit exposant $N$ possible pour une inclusion du type $(\sqrt {\cal I})^N\subseteq {\cal I}$ en fonction d'invariants globaux comme le degr\'e de g\'en\'erateurs de ${\cal I}$ et la dimension de $X$. Comme le soulignent Hickel puis \`a nouveau Heier, on peut distinguer trois \'etapes pour trouver l'exposant : on cherche d'abord une expression  pour un exposant $t$ apparaissant dans des inclusions du type $(\sqrt {\cal I})^{mt}\subseteq \overline{{\cal I}^m}$ pour tout $m$, et l'on constate qu'il suffit de borner inf\'erieurement $\bnu_{{\cal I}}(\sqrt{\cal I}))$ puis, \`a l'aide de l'interpr\'etation de $\bnu_{{\cal I}}(\sqrt{\cal I}))$ par \'eclatement normalis\'e et de r\'esultats comme ceux du \S 5 du s\'eminaire, on borne $\bnu_{{\cal I}}(\sqrt{\cal I}))$ en fonction des donn\'ees num\'eriques du probl\`eme au moyen de la th\'eorie des intersections, et enfin on applique le th\'eor\`eme de Brian\c con-Skoda (voir [L], 9.6) qui affirme que sur un espace r\'egulier de dimension $n$ on a $\overline {{\cal I}^n}\subseteq {\cal I}$. \par
Dans [H-L] et [H], on utilise le fait que $T_x({\cal I})$ est le plus petit nombre rationnel $t$ tel que l'on ait pour tout entier $m$ l'inclusion
$$(\sqrt {\cal I}_x)^{\lceil mt\rceil}\subseteq \overline{{\cal I}_x^m},$$ ce qui r\'esulte aussit\^ot de l'\'egalit\'e $T_x({\cal I})=(\bnu_{{\cal I}_x}(\sqrt {{\cal I}_x}))^{-1}$ que nous venons de voir.\par 
Soit maintenant $X$ une vari\'et\'e projective complexe non singuli\`ere et notons $n$ sa dimension. Soient ${\cal I}$ un faisceau coh\'erent d'id\'eaux de ${\cal O}_X$ et $L$ un faisceau inversible ample sur $X$ tels que $L\otimes_{{\cal O}_X}{\cal I}$ soit engendr\'e par ses sections globales. On peut d\'efinir $T({\cal I})$ comme le supremum des $T_x({\cal I})$ et Heier montre que l'on a  l'in\'egalit\'e $$T({\cal I}) \leq (L^n).$$
Appliquant le th\'eor\`eme de Brian\c con-Skoda, il en d\'eduit un r\'esultat qui se trouve aussi dans [E-L]:\par\noindent
Avec les notations ci-dessus on a l'inclusion
$$(\sqrt{\cal I})^{n(L^n)}\subseteq {\cal I}.$$\par\medskip\noindent
{\bf Compl\'ement 6: Les quotients $A/\overline{I^n}$ et l'alg\`ebre $\bigoplus_{n\geq 0}\overline{I^n}/\overline{I^{n+1}}$}\par
Dans [M1], Morales calcule le polyn\^ome de Hilbert Samuel avec lequel coincide pour $n$ grand la longueur du quotient $A/\overline{I^n}$ dans le cas o\`u $A$ est l'alg\`ebre d'un germe de courbe analytique plane r\'eduite et $I$ son id\'eal maximal. Il prouve qu'il est \'egal \`a $e_m(A)n-\delta$ o\`u $e_m(A)$ est la multiplicit\'e et $\delta =\hbox{\rm dim}\overline A/A$ est la diminution de genre. Dans [M2] il prouve qu'\'etant donn\'ee une vari\'et\'e projective $X$ sur un corps alg\'ebriquement clos, $x$ un point ferm\'e de $X$ et $I$ un id\'eal primaire de $A={\cal O}_{X,x}$, on peut interpr\'eter g\'eom\'etriquement les coefficients du polyn\^ome avec lequel lequel coincide pour $n$ grand la longueur du quotient $A/\overline{I^n}$. Dans [M3] Morales montre que si $A$ est un anneau local normal excellent de dimension 2 dont le corps r\'esiduel est alg\'ebriquement clos, l'alg\`ebre gradu\'ee $\bigoplus_{n\geq 0}\overline{I^n}/\overline{I^{n+1}}$ est de Cohen-Macaulay si la longueur du quotient $A/\overline{I^n}$  coincide avec un polyn\^ome d\`es que $n\geq 1$. Enfin, dans le cas o\`u $A$ est local excellent et $I$ est engendr\'e par un syst\`eme de param\`etres, on trouve dans [M-T-V] des conditions de r\'egularit\'e pour $A$ en fonction du comportement de la suite des longueurs des quotients $\overline {I^n}/I^n$. Ceci est reli\'e \`a l'avis de recherche 3.6 du s\'eminaire.\par\noindent
Un algorithme pour le calcul de l'alg\`ebre $\bigoplus_{n\geq 0}\overline{I^n}$ est pr\'esent\'e dans [P-U-V].
\par\medskip\noindent
{\bf Compl\'ement 7: Sp\'ecialisation sur le gradu\'e}\par
L'alg\`ebre gradu\'ee $$\overline{\hbox{\rm gr}}_{\cal I}{\cal O}_X=\bigoplus_{\nu\in \R_0}\overline{{\cal I}^\nu}/\overline{{\cal I}^{\nu +}}$$ du \S 4 est l'objet central d'\'etude du s\'eminaire, et le fait qu'elle soit une ${\cal O}_X/{\cal I}$-alg\`ebre de type fini un des r\'esultats principaux. Par analytisation, cette alg\`ebre gradu\'ee correspond donc \`a un germe d'espace analytique complexe muni d'une action de $\C^*$, qui est par construction \it r\'eduit \rm. Comme on sait  qu'une alg\`ebre filtr\'ee peut \^etre vue comme d\'eformation de l'alg\`ebre gradu\'ee associ\'ee\footnote{(*)}{Il semble que ce r\'esultat ait \'et\'e red\'ecouvert de nombreuses fois depuis depuis Gerstenhaber ([Ge]); \it cf. \rm [T5], [Po], \S 5, [Fu], Chap. 5.}, il existe une d\'eformation \`a un param\`etre du germe d'espace analytique associ\'e \`a $\hbox{\rm Spec}\overline{\hbox{\rm gr}}_{\cal I}{\cal O}_{X,x}$ dont toutes les fibres sauf la fibre sp\'eciale sont isomorphes \`a $(X,x)$ et il est int\'eressant d'explorer  g\'eom\'etriquement cette sp\'ecialisation d'un germe d'espace analytique (ou de sch\'ema excellent) $(X,x)$ sur un c\^one r\'eduit qui peut jouer le r\^ole d'un c\^one normal de $Y$ dans $X$.\par\noindent Le cas le plus simple est celui o\`u $(X,x)$ est une singularit\'e de branche plane et ${\cal I}$ est l'id\'eal maximal $m$ de l'alg\`ebre locale ${\cal O}_{X,x}$. On peut consid\'erer ${\cal O}_{X,x}$ comme une sous-alg\`ebre de sa normalisation $\C\{t\}$ et les valeurs que prend sur ${\cal O}_{X,x}$ la valuation $t$-adique forment un semigroupe $\Gamma$ d'entiers. Le premier auteur a remarqu\'e que $\overline{\hbox{\rm gr}}_m{\cal O}_{X,0}$ est dans ce cas l'alg\`ebre du semigroupe $\Gamma$ \`a coefficients dans $\C$ et que le spectre de cette alg\`ebre, c'est \`a dire la courbe monomiale correspondante, est une intersection compl\`ete. Il en r\'esulte que toutes les branches planes appartenant \`a la m\^eme classe d'\'equisingularit\'e, qui ont le m\^eme semigroupe associ\'e, sont des d\'eformations de la courbe monomiale; elles apparaissent donc dans la d\'eformation miniverselle de cette courbe. Le second auteur a dans [T5] appliqu\'e cela \`a l'\'etude des modules de branches planes, en particulier pour en donner une compactification naturelle. Pour d'autres applications de ce point de vue aux espaces de modules on renvoie \`a [C] et, en ce qui concerne la r\'esolution des singularit\'es, \`a [G-T]. \par
Dans [Kn] Allen Knutson \'etudie le remplacement  dans la th\'eorie des intersections du c\^one normal $C_{X,Y}$ d'une sous vari\'ete $Y$ de $X$ d\'efinie par l'id\'eal ${\cal I}$ par ce qu'il appelle le "balanced normal cone" $\overline C_{X,Y}$, qui est le spectre de $\overline{\hbox{\rm gr}}_{\cal I}{\cal O}_X$. Il montre que si $X, Y$ et $V$ sont des vari\'et\'es quasi-projectives r\'eduites, \'etant donn\'ee une immersion r\'eguli\`ere $X\hookrightarrow Y$ et un morphisme $V\to Y$, posant $W=V\times_YX$, le produit d'intersection raffin\'e de Fulton-MacPherson qui est un \'el\'ement de l'anneau de Chow $A^\bullet(W)$ d\'efini \`a l'aide de la sp\'ecialisation sur le c\^one normal $C_{Y,X}$ peut \^etre d\'efini aussi bien \`a l'aide de la sp\'ecialisation sur $\overline C_{X,Y}$. Cet article contient aussi des exemples. Knutson sugg\`ere que le remplacement de $C_{X,Y}$ par $\overline C_{X,Y}$ conduit \`a une th\'eorie dynamique des intersections qui garde la trace de la {\it vitesse} avec laquel les intersections transverses cr\'e\'ees apr\`es mise en position g\'en\'erale se rapprochent lorsque l'on revient \`a la position sp\'eciale. L'\'etude du cas des branches planes sugg\`ere que le passage \`a $\overline C_{X,Y}$ revient \`a faire la th\'eories des intersections apr\`es avoir tout plong\'e  dans un espace anisotrope dont les poids des coordonn\'ees correspondent aux {\it vitesses} de rapprochement possibles qui sont ''primitives''  c'est \`a dire engendrent toutes les autres. Si l'on veut \'etudier les vitesses de rapprochement dans l'espace ambient de d\'epart, alors la section 4.2 de [T 2] et les r\'esultats de [T3], [T6] section 5.15 sugg\`erent que l'objet qui mesure la distribution des vitesses est le polygone de Newton du compl\'ement 2 ci-dessus. 
\par\noindent\vskip1truecm
\centerline{BIBLIOGRAPHIE}\vskip.5truecm\noindent     
[B\"o1] Erwin B\"oger, \it Einige Bemerkungen zur theorie der ganzalgebraischen
Abh\"angigkeit in Idealen, \rm Math. Ann., {\bf 185} (1970), 303--308. \par\noindent
[Bon] Romain Bondil, \it Geometry of superficial elements, \rm  Ann. Fac. Sci. Toulouse Math. (6)  14  (2005),  no. 2, 185--200.\par\noindent
[Br] Gregory W. Brumfiel, \it Real valuation rings and ideals, \rm Springer L.N.M., No. 959, 1981.\par\noindent
[C] Pierrette Cassou-Nogu\`es, \it Courbes de semi-groupe donn\'e, \rm Rev. Mat. Univ. Complut. Madrid  4, no. 1 (1991), 13--44. \par\noindent
[C-E-S] C\v at\v alin Ciuper\v ca, Florian Enescu, Sandra Spiroff, \it Asymptotic growth of powers of ideals, \rm ArXiv: Math. AC/0610774.\par\noindent
[C-K1] Jacek Chadzy\'nski, Tadeusz Krasi\'nski,  \it The \L ojasiewicz exponent of an analytic function of two complex variables at an isolated zero, \rm Singularities 1985, Banach Center publications 20, PWN Varsovie, 1988. \par\noindent  
[C-K2] J. Chadzy\'nski, T. Krasi\'nski, \it A set on which the local \L ojasiewicz exponent is attained, \rm Annales Polon. Math., 67, (1997) 297-301 . \par\noindent  
[C-S] Craig Huneke, Irena Swanson, \it Integral closure of ideals, rings, and modules, \rm London Mathematical Society Lecture Note Series, 336. Cambridge University Press, Cambridge, 2006. xiv+431\par\noindent
[D] John P. D'Angelo, \it Real hypersurfaces, orders of contact, and applications, \rm Annals of Math., (2), 115(3) (1982), 615-637.\par\noindent
[E-L] Lawrence Ein, Robert Lazarsfeld, \it A geometric effective Nullstellensatz, \rm Invent. Math., 137(2) (1999), 427-448.\par\noindent
[F1] Azzeddin Fekak, \it Interpretation alg\'ebrique de l'exposant de \L ojasiewicz, \rm Annales Polonici Mathematici, LVI, 2 (1992), 123-131.\par\noindent
[F2] A. Fekak, \it Exposants de \L ojasiewicz pour les fonctions semi-alg\'ebriques, \rm C.R.A.S.
Paris, t. 310, S\'erie 1 (1990), 193-196.\par\noindent
[Fu] William Fulton, \it Intersection Theory, \rm Springer 1983.   \par\noindent
[G1] Terry Gaffney, \it  Integral closure of modules and Whitney equisingularity, \rm Invent. Math, 
{\bf 107} (1992), 301-322.\par\noindent
[G2] T. Gaffney, \it Polar multiplicities and equisingularity of map-germs, \rm Topology, Vol.
32, No.1 (1993), 185-223.\par\noindent
[G3] T. Gaffney, \it Multiplicities and equisingularity of ICIS germs, \rm Inventiones Math., {\bf 123} (1996), 209-220.\par\noindent
[G4] T. Gaffney,\it The theory of integral closure of ideals and modules: applications and new developments \rm With an appendix by Steven Kleiman and Anders Thorup. NATO Sci. Ser. II Math. Phys. Chem., 21,  New developments in singularity theory (Cambridge, 2000),  379--404, Kluwer Acad. Publ., Dordrecht, 2001.\par\noindent
[Ge] Murray Gerstenhaber, \it On the deformation of rings and algebras, II. \rm   Ann. of Math.  {\bf 84} (1966), 1-19.\par\noindent
[Gw] Janusz Gwo\'zdziewicz, \it The Lojasiewicz exponent  at an
isolated zero, \rm Commentarii Math. Helvetici, 74, (1999), 364-375.\par\noindent
[G-K 1] T. Gaffney, Steven L. Kleiman, \it Specialization of integral dependence for modules, \rm Inv. Math., {\bf 137} (1999), no. 3, 541-574.
\par\noindent
[G-K 2] T. Gaffney, S. Kleiman, \it $W_f$ and integral dependence, \rm Real and Complex singularities (Sao Carlos, 1998) Chapman and Hall//CRC Res. Notes in Math., 412, 33-45, Chapman and Hall//CRC Boca Raton, Florida, 2000.
\par\noindent
[GB] Evelia Garc\'\i a Barroso, \it Sur les courbes polaires d'une courbe plane r\'eduite, \rm  Proc. London Math. Soc. (3) 81 (2000) 1-28.\par\noindent
[GB-G] E. Garc\'\i a Barroso, Janusz Gwo\'zdziewicz, \it Characterization of jacobian Newton polygons of branches, \rm Manuscrit, 2007.\par\noindent
[GB-P] E. Garc\'\i a Barroso, Arkadiusz P\l oski, \it On the \L ojasiewicz numbers, \rm C. R. Acad. Sci. Paris, Ser. I., 336 (2003), 585-588.\par\noindent
[GB-K-P 1] E. Garc\'\i a Barroso, Tadeusz Krasi\'nski, A. P\l oski, \it On the \L ojasiewicz numbers, II, \rm C. R. Acad. Sci. Paris, Ser. I., 341 (2005), 357-360.\par\noindent
[GB-K-P 2] E. Garc\'\i a Barroso, T. Krasi\'nski, A. P\l oski, \it The \L ojasiewicz numbers and plane curve singularities, \rm Ann. Pol. Math., 87 (2005), 127-150.\par\noindent
[G-T] Rebecca Goldin, B. Teissier, \it Resolving plane branch singularities with
one toric morphism, \rm in ``Resolution of Singularities, a research textbook in tribute to Oscar Zariski'',
Birkh\"auser, Progress in Math. No. 18 (2000), 315-340.\par\noindent
[H] Gordon Heier, \it Finite type and the effective Nullstellensatz, \rm ArXiv: Math/AG 0603666.
\par\noindent
[Hi] Michel Hickel, \it Solution d'une conjecture de C. Berenstein-A. Yger et invariants de contact \`a l'infini, \rm  Ann. Inst. Fourier (Grenoble),  {\bf 51}, No.3 (2001), 707-744.\par\noindent
[Hu] Craig Huneke, \it Tight closure and its applications, \rm C.B.M.S. Lecture Notes 88 (1996), A.M.S., Providence\par\noindent
[Hu-S] Craig Huneke, Irena Swanson, \it Integral closure of ideals, rings, and modules, \rm 
London Mathematical Society Lecture Note Series, 336. Cambridge University Press, Cambridge, 2006.\par\noindent
[H-I-O] Manfred Hermann, Shin Ikeda, Ulrich Orbanz, \it Equimultiplicity and blowing up, an algebraic study, with an appendix by Boudewin Moonen, \rm  Springer Verlag, 1988.\par\noindent
[H-L] G. Heier, Robert Lazarsfeld, \it Curve selection for finite type ideals, \rm  ArXiv: Math/CV0506557.
\par\noindent
[I 1] Shuzo Izumi, \it A measure of integrity for
local analytic algebras, \rm Publ. R.I.M.S., Kyoto University, {\bf 21}, 4(1985), 719--735.\par\noindent
[I 2] S. Izumi, \it  Gabrielov's rank condition is equivalent to an inequality of reduced orders, \rm
Math.  Annalen, {\bf 276} (1986), 81--87.\par\noindent
[I 3] S. Izumi, \it  Fundamental properties of germs of analytic mappings of analytic sets and related topics, \rm  Real and Complex singularities, Proceedings of the Australian-Japanese Workshop, University of Sidney 2005, L. Paunescu, A. Harris, T. Fukui, S. Koike, Editors.
World Scientific, 2007, 109-123.\par\noindent
[K] Steven L.  Kleiman, \it Equisingularity, multiplicity, and dependance, \rm Commutative algebra and algebraic geometry (Ferrara), 211-225, Lecture Notes in Pure and Applied Math., 206, Dekker, New York, 1999.\par\noindent
[K-T] S. Kleiman, Anders Thorup, \it A geometric theory of the Buchsbaum-Rim multiplicity,   \rm J. Algebra, {\bf 167}, (1) (1994), 168-231.\par\noindent
[Kn] Allen Knutson, \it Balanced normal cones and Fulton-MacPherson's intersection theory, \rm   Pure Appl. Math. Q. 2, no. 4 (2006), 1103--1130. \par\noindent
[L] Robert Lazarsfeld, \it Positivity in algebraic geometry II, \rm Ergebnisse der Mathematik vol. 49, Springer Verlag 2004.\par\noindent 
[Le] Andrzej Lenarcik, \it On the jacobian Newton polygon of plane curve singularities, \rm Soumis.\par\noindent 
[M1] Marcelo Morales,  \it Le polyn\^ome de Hilbert-Samuel associ\'e \`a la filtration par les cl\^otures int\'egrales des puissances de l'id\'eal maximal pour une courbe plane. \rm  C. R. Acad. Sci. Paris S\'er. A-B  289  (1979), no. 6, A401--A404.\par\noindent
[M2] M. Morales, \it Polyn\^ome d'Hilbert-Samuel des cl\^otures int\'egrales des puissances d'un id\'eal $m$-primaire. \rm  Bull. Soc. Math. France  112  (1984),  no. 3, 343--358.\par\noindent
[M3] M. Morales, \it Cl\^oture int\'egrale d'id\'eaux et anneaux gradu\'es Cohen-Macaulay. \rm  G\'eom\'etrie alg\'ebrique et applications, I (La R\'abida, 1984),  151--171, Travaux en Cours, 22, Hermann, Paris, 1987.\par\noindent
[Mi] John Milnor, \it Singular points of complex hypersurfaces, \rm  Annals of Mathematics Studies, No. 61, Princeton U.P., 1968.\par\noindent
[M-T-V] M. Morales, Ng\^o Vi\^et Trung, Orlando Villamayor \it Sur la fonction de Hilbert-Samuel des cl\^otures int\'egrales des puissances d'id\'eaux engendr\'es par un syst\`eme de param\`etres. \rm   J. Algebra  129  (1990),  no. 1, 96--102.\par\noindent
[Mc-N] Jeffery D. McNeal, Andr\'as N\'emethi, \it The order of contact of a holomorphic ideal in $\C^2$, \rm Math. Z.,  250  (2005),  no. 4, 873--883.\par\noindent
[N] Masayoshi Nagata, \it Note on a paper of Samuel concerning asymptotic properties of powers of ideals, \rm  Mem. Coll. Sci. Univ. Kyoto, Series A, Math., {\bf 30} (1957), 165-175.\par\noindent
[No] D.G. Northcott, \it Lessons on Rings, Modules, and Multiplicities, \rm University Press,
Cambridge,1968.\par\noindent
[N-R] D.G. Northcott, D. Rees, \it Reductions of ideals in local rings, \rm Proc. Camb. Phil. Soc., {\bf 50}
(1954), 145--158.\par\noindent
[P] Patrice Philippon, \it D\'enominateurs dans le th\'eor\`eme des z\'eros de Hilbert, \rm Acta Arithm., {\bf 58},1 (1991), 1-25.\par\noindent
[P\l 1] Arkadiusz P\l oski, \it
On the growth of proper polynomial mappings, \rm 
Annales Polonici Math., XLV (1985), 297-309.\par\noindent
[P\l 2] A. P\l oski, \it Remarque sur la multiplicit\'e d'intersection des branches planes, \rm Bull. Pol. Acad. Sci. Math., 33 (1985) No. 11-12, 601-605.\par\noindent
[P\l 3] A. P\l oski, \it Multiplicity and the \L ojasiewicz exponent, \rm  in: "Singularities", Banach Center publications, 1988,  353-364.\par\noindent
[Po] Vladimir L. Popov, \it Contraction of the action of reductive algebraic groups, \rm Math. USSR Sbornik, {\bf 58} (1987), no. 2, 311-335.\par\noindent
[P-T] F. Pham et B. Teissier, \it Saturation Lipschitzienne d'une alg\`ebre analytique complexe et saturation de Zariski, \rm  Preprint 1969. Fichier .pdf disponible sur http://www.math.jussieu.fr/~teissier/old-papers.html\par\noindent
[P-U-V]  Claudia Polini, Bernd Ulrich, Wolmer V. Vasconcelos, \it  Normalization of ideals and Brian\c con-Skoda numbers, \rm  Math. Res. Lett.  12  (2005),  no. 5-6, 827--842.\par\noindent
[R1] David Rees, \it Valuations associated with ideals, \rm  Proc. London Math. Soc. (3)  6  (1956), 161--174.\par\noindent
[R2] D. Rees, \it valuations associated with ideals, II, \rm  J. London Math. Soc.  31 (1956), 221--228\par\noindent
[R3] D. Rees, \it Degree functions in local rings, \rm Proc. Camb. Phil. Soc., 57 (1961),
1-7.\par\noindent
[R4] D. Rees, \it A-transforms of local rings and a theorem on multiplicities of ideals, \rm Proc. Camb. Phil.
Soc., 57 (1961), 8--17.\par\noindent
[R5] D. Rees, \it Local birational Geometry, \rm Actas del Coloquio internacional sobre Geometr\'\i a
algebraica, Madrid, Sept. 1965.\par\noindent
[R6] D. Rees, \it Multiplicities, Hilbert functions and degree functions, \rm
Commutative Algebra-Durham 1981, London Math. Soc. Lectures Notes 72 (Ed. R.Y. Sharp,
University Press, Cambridge 1983), pp 170-178.\par\noindent 
[R7] D. Rees, \it Hilbert functions and pseudo-rational local rings of dimension two, \rm J. London Math.
Soc. (2), 24 (1981), 467-479.\par\noindent
[R8] D. Rees, \it Rings associated with ideals and analytic spread, \rm Math. Proc. Camb. Phil. Soc., 
{\bf 89} (1981), 423-432.\par\noindent
[R9] D. Rees, \it Generalizations of reductions and mixed
multiplicities, \rm J. London Math. Soc., (2), {\bf 29} (1984), 397-414.\par\noindent
[R10] D. Rees, \it The general extension of a local ring and mixed multiplicities, \rm Springer Lecture Notes in mathematics No. 1183 (1986), 339-360.\par\noindent
[R11] D. Rees, \it Asymptotic properties of ideals, \rm London Math. Soc. Lecture Series, 113, 1988.\par\noindent
[R12] D. Rees, \it Izumi's theorem, \rm Commutative Algebra (Berkeley, CA., 1987), Math. Sci. Res. Inst. Publications, 15, Springer New-York, 1989.\par\noindent
[R-S] D. Rees, Rodney Y. Sharp, \it On a Theorem of B. Teissier on multiplicities of ideals in local rings, \rm J. London Math. Soc., (2), 18 (1978), 449-463.\par\noindent
[Sa] Pierre Samuel, {\it Some asymptotic properties of powers of ideals,} Annals of Math., (2), {\bf 56} (1952), 11-21 .\par\noindent
[T1] Bernard Teissier, \it  Cycles \'evanescents, sections planes, et conditions de Whitney, \rm Singularit\'es \`a Carg\`ese, Ast\'erisque No. 7-8, S.M.F., (1973), 285-362.\par\noindent
[T2]  B. Teissier, \it Jacobian Newton polyhedra and equisingularity, \rm Proceedings
R.I.M.S. Conference on singularities, Kyoto,  April
1978. (Publ. R.I.M.S. 1978) et traduction
dans: \it S\'eminaire sur les singularit\'es, \rm Publ. Math. Universit\'e  Paris
VII no.7, (1980), 193-211. Fichier .pdf disponible sur http://www.math.jussieu.fr/~teissier/articles-Teissier.html\par\noindent
[T3] B. Teissier, \it Vari\'et\'es polaires I; invariants polaires des singularit\'es d'hypersurfaces, \rm  Inventiones Math. 40 (1977), 267-292. \par\noindent
[T4] B. Teissier, \it Vari\'et\'es polaires II; multiplicit\'es polaires, sections planes, et conditions de Whitney, \rm Proc. Conf. Algebraic Geometry, La R\'abida, Springer Lecture Notes in Math., no. 961, 314-491.\par\noindent
[T5] B. Teissier, \it Appendice: la courbe monomiale et ses d\'eformations, \rm in: Oscar Zariski, "Le probl\`eme des modules pour les branches planes", Publ. Ecole Polytechnique, Paris 1975, reprinted by Hermann ed., Paris, 1986, English translation by Ben Lichtin in \it The moduli problem for plane branches, \rm  University Lecture Series, Vol. 39, A.M.S., 2006.\par\noindent
[T6] B. Teissier, \it The Hunting of invariants in the Geometry of discriminants, \rm in: Real and complex singularities, Oslo 1976, Per Holm editeur, Sijthoff \& Noordhoff 1977, pages 565-677.\par\noindent
[T7] B. Teissier, \it R\'esolution simultan\'ee II, \rm in S\'eminaire sur les Singularit\'es des Surfaces, Lecture Notes in Mathematics, No. 777. Springer, Berlin, 1980. 82-146. Fichier .pdf disponible sur http://www.math.jussieu.fr/~teissier/articles-Teissier.html\par\noindent
[V] Vasconcelos, Wolmer, \it  Integral closure, Rees algebras, multiplicities, algorithms, \rm Springer Monographs in Mathematics. Springer-Verlag, Berlin, 2005.\par\noindent
[V-S] Vui, H\`a Huy, Ph\d am Tien So'n, \it Newton-Puiseux approximation and \L ojasiewicz exponents, \rm  Kodai Math. J.  26  (2003),  no. 1, 1--15.\par\noindent
[Z] Oscar Zariski, \it Le probl\`eme des modules pour les branches planes, \rm  Publ. Ecole Polytechnique, Paris 1975, reprinted by Hermann ed., Paris, 1986, English translation by Ben Lichtin in \it The moduli problem for plane branches, \rm  University Lecture Series, Vol. 39, A.M.S., 2006.\par\noindent

\end